\documentclass[11pt]{amsart}
\usepackage{latexsym,amssymb,amsxtra,amsmath,amsfonts}
\usepackage{verbatim,appendix}
\usepackage[pdftex]{graphicx,color,graphics} 
\usepackage[top=1in, bottom=1in, left=1in, right=1in]{geometry}
\allowdisplaybreaks

\newtheorem{theorem}{Theorem}
\newtheorem{lemma}[theorem]{Lemma}
\newtheorem{proposition}[theorem]{Proposition}
\newtheorem{corollary}[theorem]{Corollary}
\newtheorem{conjecture}[theorem]{Conjecture}
\newtheorem{remark}[theorem]{Remark}

\theoremstyle{definition}

\newtheorem{example}[theorem]{Example}

\newcommand{\beq}{\begin{equation}}
\newcommand{\eeq}{\end{equation}}
\newcommand{\beqa}{\begin{eqnarray}}
\newcommand{\eeqa}{\end{eqnarray}}
\newcommand{\beaa}{\begin{eqnarray*}}
\newcommand{\ben}{\begin{eqnarray*}}
\newcommand{\eaa}{\end{eqnarray*}}
\newcommand{\een}{\end{eqnarray*}}

\newcommand{\leftexp}[2]{{\vphantom{#2}}^{#1}{#2}}

\def\al{\alpha}
\def\be{\beta}

\def\de{\delta}
\def\ep{\varepsilon}

\def\la{\lambda}
\def\si{\sigma}
\def\Si{\Sigma}
\def\Om{\Omega}

\def\d{\partial}

\def\tI{\widetilde{I}}
\def\tGamma{\widetilde{\Gamma}}
\def\tf{\widetilde{\bf f}}

\def\cD{\mathcal{D}}

\def\cO{\mathcal{O}}
\def\cF{\mathcal{F}}
\def\cT{\mathcal{T}}
\def\cA{\mathcal{A}}
\def\cH{\mathcal{H}}
\def\cW{\mathcal{W}}
\def\tDelta{\widetilde{\Delta}}
\def\tal{\widetilde{\alpha}}
\def\tphi{\widetilde{\phi}}
\def\tomega{\widetilde{\omega}}
\def\tlieh{\widetilde{\mathfrak{h}}}

\def\Z{\mathbb{Z}}
\def\R{\mathbb{R}}
\def\C{\mathbb{C}}
\def\f{{\bf f}}
\def\q{{\bf q}}
\def\t{{\bf t}}
\def\one{{\bf 1}}

\def\y{{\bf y}}
\def\x{{\bf x}}

\def\oDelta{\Delta^{(0)}}
\def\olieg{\mathfrak{g}^{(0)}}
\def\olieh{\mathfrak{h}^{(0)}}
\def\osi{\sigma_b}
\def\oLambda{\Lambda^{(0)}}

\def\lieh{{\mathfrak{h}}}
\def\lieg{{\mathfrak{g}}}
\def\lies{{\mathfrak{s}}}
\def\liea{{\mathfrak{a}}}
\def\lieI{{\mathfrak{I}}}

\begin{document}

\setlength{\unitlength}{1mm}

\title[GW of Fano orbifold curves, Gamma integral structures and ADE-Toda Hierarchies]{Gromov--Witten theory of Fano orbifold curves, Gamma integral structures and ADE-Toda Hierarchies}

\author{Todor Milanov}
\address{Kavli IPMU\\
University of Tokyo (WPI)\\
Japan}
\email{todor.milanov@ipmu.jp}

\author{Yefeng Shen}
\address{Kavli IPMU\\
University of Tokyo (WPI)\\
Japan}
\email{yefeng.shen@ipmu.jp}

\author{Hsian-Hua Tseng}
\address{
Department of Mathematics\\
Ohio State University\\
USA}
\email{hhtseng@math.ohio-state.edu}

\begin{abstract}
We construct an integrable hierarchy in the form of Hirota quadratic equations (HQE) that governs the Gromov--Witten (GW) invariants of the Fano orbifold projective curve $\mathbb{P}^1_{a_1,a_2,a_3}$. The vertex operators in our construction are given in terms of the $K$-theory of $\mathbb{P}^1_{a_1,a_2,a_3}$ via Iritani's $\Gamma$-class modification of the Chern character map.  We also identify our HQEs with an appropriate Kac--Wakimoto hierarchy of ADE type. In particular, we obtain a generalization of the famous Toda conjecture about the GW invariants of $\mathbb{P}^1$ . 
\end{abstract}

\date{\today}

\maketitle

\tableofcontents

\section{Introduction}
Witten's conjecture \cite{W1}, proven by Kontsevich \cite{Ko1} states
that the GW theory of $X=\text{pt}$ is governed by the KdV
hierarchy. Although Witten was cautious in proposing that there should
be an integrable hierarchy for every target $X$, several groups of
physicists and mathematicians, including Witten himself \cite{W2}, have tried to find a generalization of Witten's conjecture. The next important discovery was the Toda conjecture \cite{EY,Ge,Z}, proven by \cite{Ge, DZ3, M}. It states that the GW theory of $X=\mathbb{P}^1$ is governed by the extended Toda hierarchy (see \cite{CDZ}  for the definition in terms of a Lax operator and \cite{Ge,Z} for the bi-Hamiltonian definition). The Toda conjecture was further generalized by Milanov and Tseng \cite{MT} (see also \cite{J, CvdL}) by allowing the target $X$ to be a projective line with two orbifold points. The corresponding integrable hierarchy is the extended bigraded Toda hierarchy, which was introduced and studied by Carlet \cite{Ca}. The relationship between topological field theories and integrable hierarchies is studied in other examples, such as \cite{G1, GM, FSZ, FJR, FJR2, LRZ}. 

Motivated by GW theory, Dubrovin--Zhang \cite{DZ} proposed a general construction based on the theory of semi-simple Frobenius manifolds. While their construction produces flows that are rational functions on the jet variables, it was expected that for the important classes of semi-simple Frobenius manifold, such as quantum cohomology, the flows are in fact polynomial and  that the hierarchy can be used to
compute uniquely the higher genus invariants. The polynomiality of the flows for a semi-simple Frobenius manifold associated with a cohomological field theory (this includes the case of GW theory) was proved recently by Buryak--Posthuma--Shadrin \cite{BPS, BPS2} using the higher genus reconstruction of Givental. In particular, Witten's conjecture generalizes for all targets $X$ that have semi-simple quantum cohomology. The discovery of this new class of integrable hierarchies is a major breakthrough in the theory of integrable systems. It is natural to study further their properties and to look for applications to other areas of Mathematics and even beyond.

The higher genus reconstruction of Givental which was mentioned above is one of the major achievements in GW theory. The reconstruction was discovered and proved by Givental in the equivariant settings when $X$ is equipped with a torus action with isolated fixed points \cite{G2}. Based on his work \cite{G2}, Givental conjectured a certain higher genus reconstruction formula for the total ancestor potential of $X$ with semi-simple quantum cohomology. Givental's conjecture was proved in various cases in \cite{G3, JK, I, BCFK}, and in full generality by C. Teleman \cite{Te}. Givental's reconstruction inspires an approach to study the relation between GW theory, representation theory of vertex algebras, and integrable systems. In this approach one aims at constructing an integrable hierarchy in the form of {\em Hirota quadratic equations} (HQE)\footnote{The word ``quadratic'' in HQE was used by Givental in \cite{G1}. The equations are also known as ``Hirota bilinear equations.''} and show that the generating function of GW invariants is a tau-function of the hierarchy (i.e. it satisfies the HQEs). This approach has been successfully worked out for GW theory of $X$ when $X=\mathbb{P}^1$ \cite{M, M2} and $X=\mathbb{P}^1_{a,b}$ \cite{MT}. See also \cite{G1, GM, FGM} for instances of this approach in the setting of singularity theory.

While the construction of Dubrovin and Zhang is general, the approach with HQEs is not so easy to generalize. The main difficulty is that we have to deal with vanishing cycles and period integrals whose properties are still not very well understood. This is probably one of the main motivation to pursue the HQEs approach. It gives us a new motivation and a new view point in the theory of vanishing cycles and period integrals. Let us point out that there are no examples of targets $X$ of dimension $>1$ for which the HQEs are known to exist, although there are some indications that such examples exists (see \cite{Br}).  Even for orbifolds of dimension 1 (with semi-simple quantum cohomology) the HQEs are not known in general. In this paper, we would like to solve this problem for Fano orbifolds of dimension 1, i.e., $\mathbb{P}^1$-orbifolds $\mathbb{P}^1_{a_1,a_2,a_3}$ (with 3 orbifold points), s.t., $1/a_1+1/a_2+1/a_3>1$. It was already noticed in \cite{CDZ} that the extended Toda hierarchy is equivalent to an extended Kac--Wakimoto hierarchy of type $A_1$. While KdV is the so called principal Kac--Wakimoto hierarchy of type $A_1$, the extended Toda hierarchy is obtained by extending the homogeneous Kac--Wakimoto hierarchy of type $A_1$. Our main result is that for the remaining Fano  $\mathbb{P}^1_{a_1,a_2,a_3}$-orbifolds the corresponding integrable hierarchy is an extension of a Kac--Wakimoto hierarchy as well, but this time it is neither the homogeneous, nor the principal realization, but something in between.

\subsubsection{GW theory of Fano orbifold curves}\label{sec:gw}

Let
\ben
{\bf a}=\{a_1,a_2,a_3\},\quad a_1\leq a_2\leq a_3,
\een 
be a triple of positive integers. Let  $\mathbb{P}^1_{\bf a}$ be the orbifold projective line obtained from $\mathbb{P}^1$ by adding\footnote{For example, by root constructions \cite{AbGV}, \cite{Cadman}.} $\mathbb{Z}_{a_1}$-, $\mathbb{Z}_{a_2}$-, and $\mathbb{Z}_{a_3}$-orbifold points. The nature of the problem of constructing HQEs depends on the orbifold Euler characteristic of $\mathbb{P}^1_{\bf a}$:
\ben
\chi:=\frac{1}{a_1}+\frac{1}{a_2}+\frac{1}{a_3}-1.
\een 
In this paper we will study the Fano case $\chi>0$, leaving the other two cases $\chi=0$ (elliptic) and $\chi<0$ (hyperbolic) for a future investigation.

We consider the Chen-Ruan orbifold cohomology of $\mathbb{P}^1_{\bf a}$,
$$
H:=H_{\rm CR}(\mathbb{P}^1_{\bf a}, \C).
$$ 
As a vector space $H$ is just $H^*({\rm I} \mathbb{P}^1_{\bf a}, \C)$,
where ${\rm I} \mathbb{P}^1_{\bf a}$ is the {\em inertia orbifold} of $\mathbb{P}^1_{\bf a}$,
\ben
{\rm I} \mathbb{P}^1_{\bf a}=\{(x,g)\ |\ x\in \mathbb{P}^1_{\bf a},\ g\in \operatorname{Aut}(x)\}.
\een
We can fix a homogeneous basis $\{\phi_i\}_{i\in\mathfrak{I}}$ of $H$, where the index set is defined by 
\beq\label{ind-set}
\mathfrak{I}:=\mathfrak{I}_{\rm tw}\cup\{(01),(02)\}:=\left\{(k,p)\mid 1\leq k\leq 3, 1\leq p\leq a_{k}-1\right\}\cup\{(01),(02)\}.
\eeq
The index set $\mathfrak{I}$ reflects the structure of the forgetful map ${\rm I} \mathbb{P}^1_{\bf
  a}\to \mathbb{P}^1_{\bf a},$ $(x,g)\mapsto x$. The connected
components of ${\rm I} \mathbb{P}^1_{\bf a}$ split into several types
depending on their fate under the forgetful map. The entry $k$
enumerates the different types, while $p$ enumerates the cohomology
classes supported on the connected components of type $k$.
Motivated by GW theory, Chen and Ruan (see \cite{CR}) have introduced
a new product, called {\em Chen--Ruan} or {\em orbifold cup product}. It is defined
as the degree-$0$ component of the quantum cup product. It is
graded homogeneous with respect to a new grading denote by
$\operatorname{deg}_{\rm CR}$. In our notation, $\phi_{01}={\bf 1}$ is the unit, $\deg_{\rm CR}\phi_{02}=1$, and $\deg_{\rm CR}\phi_{k,p}=p/a_k.$

The main objects in the orbifold GW theory of $\mathbb{P}^1_{\bf a}$ are the moduli spaces  $\overline{\mathcal{M}}_{g,n}(\mathbb{P}^1_{\bf a},d)$ of orbifold stable maps $f$ from a domain orbifold genus $g$ curve $\Sigma$ with $n$ marked points, to the target orbifold $\mathbb{P}^1_{\bf a}$, such that the homology class of the image of $f$ is $d$ times the fundamental class of the underlying curve of $\mathbb{P}^1_{\bf a}$. The descendant GW invariants (see \eqref{des-inv}) are intersection numbers on the moduli space of stable maps, denoted by
\ben
\langle \phi_{1}\,\psi_1^{k_1},\dots,\phi_{n}\psi_n^{k_n}\rangle_{g,n,d},
\een
where $\phi_j\in H$ and $\psi_{j}$ is the $j$-th $\psi$-class on the moduli space of stable maps.

Our main interest is in the so-called {\em total descendant  potential}, defined by the following generating 
series of GW invariants:
\beq\label{total-desc}
\cD_{\bf a}(\hbar;\t) = \exp\Big( \sum_{g,n,d} \hbar^{g-1}\,\frac{Q^d}{n!}
\langle \t(\psi_1),\dots,\t(\psi_n)\rangle_{g,n,d}\Big),
\eeq
where $Q$ is a non-zero complex number called the Novikov variable, $\t(z):=t_0+t_1z+t_2
z^2+\cdots$,
with $t_0,t_1,\ldots\in H $ and $\hbar$ are formal variables. Using the so called {\em dilaton shift} $q_m=t_m-\delta_{m,1}\one$ we denote $\cD_{\bf a}(\hbar;\t)$ by $\cD_{\bf a}(\hbar;\q)$ and identify it with a vector in a certain {\em Fock space} (see \eqref{eq:fock}).

The construction of the HQEs is given in Section \ref{sec:hqe}. It
relies on the theory of vanishing cycles and period integrals
associated to a Landau-Ginzburg (LG) mirror model of $\mathbb{P}^1_{\bf
  a}$. An appropriate mirror was constructed in \cite{MT} in the case
$a_1=1$ and in general by P. Rossi \cite{Ro}, who managed to compute
the quantum cohomology of $\mathbb{P}^1_{\bf a}.$ We also need to know
how to solve the quantum differential equations in terms of period
integrals. This was achieved recently by
Ishibashi--Shiraishi--Takahashi \cite{IST}, where the mirror model was constructed from the miniversal deformation space $M$ of the {\em affine cusp polynomial} 
$$f_{\bf a}(x_1,x_2,x_3):=x_1^{a_1}+x_2^{a_2}+x_3^{a_3}-\frac{1}{Q}x_1x_2x_3.$$
With such a mirror model
at hands we can pursue the same idea as in \cite{G1, FGM, GM, M, MT} to construct HQEs for GW theory of $\mathbb{P}^1_{\bf a}$ by studying periods of the mirror model $f_{\bf a}$. 
   
\subsubsection{$\Gamma$-conjecture for the Milnor lattice}\label{subsec:gamma}
Compared to the earlier works, one novelty of this paper is that we
made use of Iritani's integral structure \cite{Ir} (see also \cite{KKP}), which allows us to
express the vertex operators in our construction in terms of
K-theory. This observation seems to be quite general, so we formulated
a conjecture for the general case (see Conjecture \ref{Gamma-conj}
below), which we refer to as the {\em $\Gamma$-conjecture for the
  Milnor lattice}. 

The homology space $\lieh$ of the Milnor fiber at a reference point $(0,1)\in M\times\C$ has a lattice structure on the vanishing cycles, called {\em the Milnor lattice}. Conjecture \ref{Gamma-conj} in the case of $\mathbb{P}^1_{\bf a}$ says that for each element of the $K$-group $K(\mathbb{P}^1_{\bf a})$, there exists a corresponding integral cycle in the Milnor lattice of the mirror model, such that both integrals structures match.  
We will give a proof of this conjecture for the Fano orbifold curves $\mathbb{P}^1_{\bf a}$, based on Iritani's proof for the $\Gamma$-conjecture of toric orbifolds $\mathbb{P}^2_{\bf a}$. After inverse Laplace transformations, this allows us to get explicit formulas for the calibrated periods over the Milnor lattice, in terms of integral structures in the A-model quantum cohomology, see formula \eqref{calp-im}. Then we can embed the root system of vanishing cycles into the quantum cohomology via period maps. The vanishing cycles form an affine root system of type $A, D,$ or $E$, and we can identify the classical monodromy of the Milnor lattice with an affine Coxeter transformation, see Proposition \ref{affine}. The calibrated periods over the vanishing cycles are very important for the construction of the vertex operators later.

\subsubsection{The Kac--Wakimoto hierarchy}\label{subsec:basic-rep}
The triplets ${\bf a}=\{a_1,a_2,a_3\}$ with $\chi>0$ are classified by the Dynkin diagrams
of type $ADE$ together with a choice of a {\em branching node}. In the
$D$ and $E$ cases there is a unique choice of a branching node, while
in the $A$-case any node can be chosen. By removing the branching node we
obtain 3 diagrams of type\footnote{if $a_k=1$ then the corresponding diagram is empty.} $A_{a_k-1}$, $k=1,2,3$. Let us denote by $\lieh^{(0)}$ the Cartan
subalgebra of the corresponding simple Lie algebra $\lieg^{(0)}$ and define (cf. eqn. \eqref{si-0})
\beq\label{si-0-intro}
\si_b=\prod_{k=1}^3 \Big( s^{(0)}_{k,a_k-1}\cdots s^{(0)}_{k,2}s^{(0)}_{k,1} \Big),
\eeq
where $s^{(0)}_{k,p}:\lieh^{(0)}\to\lieh^{(0)}$ is the reflection through the hyperplanes orthogonal to $\gamma_{k,p}^{(0)}$, which is the $p$-th simple root on the $k$-th branch of the Dynkin diagram. 
The automorphism $\si_b$ can be extended to a Lie algebra automorphism of $\lieg^{(0)}$. Let us denote by 
$\kappa$ the order of $\si_b$ as an automorphism of $\lieg^{(0)}$. Due to a mirror symmetry phenomenon the spectrum of $\si_b$ is given by the degrees of the cohomology classes $\phi_i$.
More precisely, there exists a $\osi$-eigenbasis $\{H_i\}_{i\in \mathfrak{I}}$ 
of $\lieh^{(0)}$, s.t.,  $\si_b(H_i)=e^{-2\pi\sqrt{-1} d_i}H_i$, where 
$d_i=1-\operatorname{deg}_{\rm CR}(\phi_{i}).$ The index set $\mathfrak{J}$ admits a natural involution $*$ compatible with the
Poincar\'e pairing: 
\ben
d_i+d_{i^*}=1. 
\een
The $\osi$-eigenbasis can be normalized so that $(H_i|H_{j^*})=\kappa\,\delta_{i,j}$. 

The Kac-Wakimoto hierarchy corresponding to the conjugacy class of
$\si_b$ in the Weyl group can be described as follows. Let $\C[y]$ be
the algebra of polynomials on $y=(y_{i,\ell})$,  $i\in
\mathfrak{I}\backslash\{(01)\}$ and $\ell\geq 0.$
The vector space\footnote{This is a direct product of copies of $\mathbb{C}[y]$ indexed by $n\in \mathbb{Z}$.} $\C[y]^\Z$ is equipped with the structure of a module
over the algebra of differential operators in $e^\omega$ by setting
\ben
(e^\omega\cdot \tau)_n = \tau_{n-1},\quad (\d_\omega\cdot \tau)_n =
n\tau_n, \quad \tau=(\tau_n)_{n\in \mathbb{Z}}\in \mathbb{C}[y]^\mathbb{Z}.
\een

For every root $\al\in\Delta^{(0)}$ of $\lieg^{(0)}$ we define {\em vertex operators} $E^{(0)}_\al(\zeta)$ (see \eqref{E0}) and $E^{*}_\al(\zeta)$ (see \eqref{vop:twisted}) in Section \ref{kac-peterson}, both acting on $\C[y]^\Z$. Let $E_\al(\zeta)=E^{(0)}_\al(\zeta) E^{*}_\al(\zeta)$.
The HQE of the $\osi$-twisted Kac--Wakimoto hierarchy are given by the following
bilinear equation for $\tau=(\tau_n(y))_{n\in \Z}$: 
\begin{equation}\label{HQE_KW}
\begin{split}
\operatorname{Res}_{\zeta=0}\frac{d\zeta}{\zeta}\Big( \sum_{\al\in \oDelta} a_\al(\zeta) E_\al(\zeta)\otimes E_{-\al}(\zeta) \Big)\,
\tau \otimes \tau
=\Big(  \frac{1}{12}\sum_{k=1}^3 
\frac{a_k^2-1}{a_k}+
\frac{\chi}{2}\, (\d_\omega\otimes 1-1\otimes \d_\omega)^2+\\
+\, \sum_{i\in\mathfrak{I}\backslash\{(01)\}}  \sum_{\ell\geq0}^{\infty} ( d_{i^*}+\ell) (y_{i,\ell}\otimes
1-1\otimes y_{i,\ell} )(\d_{y_{i,\ell}}\otimes 1-1\otimes \d_{y_{i,\ell}})
\Big) \tau\otimes \tau,
\end{split}
\end{equation}
with the coefficients $a_\al(\zeta)$ defined by \eqref{a-al} in Section \ref{sec:KW}. 

\subsubsection{The main result}
We can write the Kac--Wakimoto HQE in terms of the descendant
variables $\{q_m\}_{m\geq 0}$, using the change of
variables between $y_{i,\ell}$ and $q_\ell^i$ (see \eqref{dv:change1}--\eqref{dv:change2}). 
Our main result can be stated as follows.
\begin{theorem}\label{t1}
Let $\cD_{\bf a}(\hbar;\q)$ (with ${\bf a}=\{a_1,a_2,a_3\}$) be the total descendant potential \eqref{total-desc}
of an orbifold projective line $\mathbb{P}^1_{\bf a}$ with a positive
orbifold Euler characteristic, then the sequence
$\left(\tau_n(\hbar;\q)\right)_{n\in\Z}$ of formal power series defined by
\ben
\tau_n(\hbar;\q) = (\kappa^\chi Q)^{\frac{1}{2}n^2}\cD_{\bf a}(\hbar;\q+n\sqrt\hbar\one ),\quad
n\in \Z.
\een
is a solution to the $\si_b$-twisted Kac--Wakimoto HQE \eqref{HQE_KW}, where $\si_b$ is the element 
\eqref{si-0-intro} of the Weyl group of the corresponding finite root system.
\end{theorem}
In other words, Theorem \ref{t1} shows that the GW theory of $\mathbb{P}^1_{\bf a}$ is governed by the Kac-Wakimoto hierarchy associated to the triple ${\bf a}$. Let us emphasize that the variables $q_1^{01},q_2^{01},\dots$ appear as parameters in the differential equations for $\tau$.  It is natural to expect that the $\si_b$-twisted Kac--Wakimoto HQE can be extended in order to include differential equations in $q_1^{01},q_2^{01},\dots$ as well. We hope that our work will motivate the specialists in integrable systems and representation theory to investigate more systematically the possibility of extending the Kac--Wakimoto hierarchies. For example, in the case of Dynkin diagrams of type $A$, our hierarchy should agree with a certain reduction of the 2D Toda hierarchy and the required extension was constructed by G. Carlet \cite{Ca} based on the ideas of \cite{CDZ}.  For the type $D$ and $E$ cases, the extension can be constructed with the same idea as in \cite{M2} with a slight necessary modification. The details will be presented elsewhere. We suggest to call the $\si_b$-twisted Kac--Wakimoto hierarchy appearing in Theorem \ref{t1} the {\em ADE-Toda hierarchy}, and call the corresponding extension the {\em Extended ADE-Toda hierarchy}. 

Our approach to Theorem \ref{t1} systematically explores representation theoretic properties of the Landau-Ginzburg mirror of $\mathbb{P}^1_{\bf a}$ and realizes these properties in quantum cohomology of $\mathbb{P}^1_{\bf a}$ using the period maps. A new observation is that we can also use K-theory to obtain explicit formulas for the leading terms of the period mapping. In particular, this simplifies the analysis of the monodromy representation. Such an approach should be helpful for more general target spaces as well. 

To our knowledge, Theorem \ref{t1} is the first case where the problem of constructing HQEs governing GW theory of a target $X$ is solved for a {\em non-toric} $X$.

Finally, it is very interesting also to investigate the relation between the integrable hierarchies obtained by applying Dubrovin and Zhang's construction \cite{DZ2} to the quantum cohomology of $\mathbb{P}^1_{\bf a}$ and the integrable hierarchies in Theorem \ref{t1}. It is natural to expect that the two approaches yield the same integrable hierarchy. We hope to return to this problem in the near future. 

\subsubsection{Outline of the proof of Theorem \ref{t1}}
First, the hierarchy (\ref{HQE_KW}) is shown to be equivalent (via a
Laplace transform) to another hierarchy (\ref{eth-ade}) defined for
affine cusp polynomials, see Theorem \ref{t2}. Then by Proposition
\ref{anc-desc}, the descendant potential $\mathcal{D}_{\bf a}$
satisfies the hierarchy (\ref{eth-ade}) if and only if the ancestor
potential $\mathcal{A}_t$ (see equation (\ref{ancestor})) satisfies
another hierarchy (\ref{eth-anc}). Finally, the most difficult step is
to prove (Theorem \ref{thm:ancestor}) that $\mathcal{A}_t$ indeed satisfies (\ref{eth-anc}). 

Let us point out that although our proof of Theorem \ref{thm:ancestor} follows closely the argument of \cite{GM}, we managed to simplify one of the crucial steps in \cite{GM}. Namely, there is a certain analyticity property (c.f. Section \ref{sec:phase-factors}) of the so called phase factors that was previously established via the theory of finite reflection groups and their relation to Artin groups. This is one of the main obstacles to generalize the result of \cite{GM} to other singularities. Our argument now seems to apply in much more general settings, since it relies only on the fact that the Gauss--Manin connection has regular singularities and that the vertex operators are local to each other (in the sense of the theory of vertex operator algebras). 

The rest of this paper is organized as follows. In Section
\ref{sec:ogw}, we recall the orbifold GW theory for Fano projective
curves $\mathbb{P}^1_{\bf a}$ and the corresponding LG mirror model. 
In Section \ref{sec:int_GW} we recall Iritani's integral structure
(see \cite{Ir}) in the quantum cohomology of a smooth projective
orbifold $X$. Furthermore, we prove that for $X=\mathbb{P}^1_{\bf a}$ the integral
structure corresponds to the Milnor lattice under mirror
symmetry. Finally, using the period mapping we identify the root system
arising from the set of vanishing cycles with an affine root system in the quantum cohomology of $\mathbb{P}^1_{\bf a}$. The integral structures allows us to obtain an explicit description of the leading order terms of the period mapping in terms of finite root systems. In Section \ref{sec:basic-rep}, using the results from Section \ref{sec:int_GW}, we give a Fock-space realization of the basic representations of the affine Lie algebras of ADE type. Then we recall the Kac-Wakimoto hierarchies and construct integrable hierarchies for affine cusp polynomials and show that these hierarchies are related by a Laplace transform (Theorem \ref{t2}). In Section \ref{sec:eth-anc} we construct another hierarchy (\ref{eth-anc}) and describe its relation with the hierarchies from previous sections, see Proposition \ref{anc-desc}. Then we show that the ancestor potential of $\mathbb{P}^1_{\bf a}$ satisfies the integrable hierarchy (\ref{eth-anc}) and deduce Theorem \ref{t1}. In Section \ref{ex} we consider the example ${\bf a}=\{2,2,2\}$. In the appendix, we give an alternative proof for the higher genus reconstruction of total ancestor potential. \\

\noindent 
{\bf Acknowledgments}. We thank Yongbin Ruan for his interests in this work, and for his comments and suggestions. T. M.  would like to thank Tadashi Ishibe for useful discussions on the affine Artin group, Atsushi Takahashi and Yuuki Shiraishi for useful discussion on mirror symmetry for affine cusp polynomials, and especially to S. Galkin and H. Iritani for very useful discussions on the $\Gamma$-conjecture. Y. S. would like to thank Shifra Reif for helpful discussions on affine Kac--Moody algebras. H.-H. T. thanks H. Iritani for discussions on mirror symmetry. 

T. M. and Y. S. acknowledge the World Premiere International Research Center Initiative (WPI Initiative), Mext, Japan. H.-H. T. thanks the hospitality and support for his visits to IPMU. 

T. M. is supported in part by JSPS Grant-in-Aid 26800003. H.-H. T. is supported in part by Simons Foundation Collaboration Grant.

\section{Orbifold GW theory of Fano orbifold curves $\mathbb{P}^1_{\bf a}$ and their mirror symmetry}\label{sec:ogw}

\subsection{Orbifold GW theory of $\mathbb{P}^1_{\bf a}$}
Fano orbifold curves are closed orbifold curves with positive orbifold Euler characteristics. They are classified by triplets of positive integers 
${\bf a}=\{a_1,a_2,a_3\}$
where $a_1\leq a_2\leq a_3$ and 
$\chi:=\frac{1}{a_1}+\frac{1}{a_2}+\frac{1}{a_3}-1>0.$
Each Fano orbifold curve is an orbifold curve with an underlying curve $\mathbb{P}^1$ and has at most three orbifold points $p_k$ ($k=1,2,3$) with local isotropy groups $\mathbb{Z}_{a_k}$. We denote such an Fano orbifold curve by $\mathbb{P}^1_{\bf a}$. Note that such notation also includes the smooth curve $\mathbb{P}^1$ with $a_1=a_2=a_3=1$. It is easy to see that $\chi$ is  the orbifold Euler characteristic of $\mathbb{P}^1_{\bf a}$.

We use the index set $\mathfrak{I}$ (see \eqref{ind-set}) to label a
fixed basis of the Chen-Ruan orbifold cohomology $H:=H_{\rm CR}(\mathbb{P}^1_{\bf a};\C)$ as follows:
\ben
\phi_{01}=\one,\quad \phi_{02}= P
\een
are the unit and the hyperplane class of the underlying $\mathbb{P}^1$ respectively and
\ben
\phi_i=\phi_{k,p},\quad i:=(k,p)\in\mathfrak{I}_{\rm tw}.
\een
are the units of the corresponding twisted sectors of $\mathbb{P}^1_{\bf a}$.
The cohomology degrees of the classes are:
\ben
\deg_{\rm CR}\phi_{01}=0, \quad \deg_{\rm CR}\phi_{02}=1, \quad
\deg_{\rm CR}\phi_i=\frac{p}{a_{k}}, \quad i=(k,p)\in\mathfrak{I}_{\rm tw},
\een
where slightly violating the standard conventions we work with complex
degree, i.e., half of the usual real degrees. 
There is a natural involution $*$ on $\mathfrak{I}$ induced by orbifold Poincar\'e duality
\beq\label{eqn:involution}
(01)^*=(02),\quad (k,p)^*=(k,a_{k}-p).
\eeq
The orbifold Poincar\'e pairing $(-,-)$ on $H$ is non-zero only for the following cases:
\ben
(\phi_{01},\phi_{02})=1,\quad (\phi_{i},\phi_{j})=
\frac{1}{a_{i}}\delta_{i,j^*},
\een 
where $i,j\in \mathfrak{I}_{\rm tw}$ correspond to twisted classes, and we set
$a_i:=a_{k}$ for $i=(k,p)\in\mathfrak{I}_{\rm tw}.$
 
GW theory studies integrals over moduli spaces of stable maps. In this paper, we will use both the descendant invariants and the ancestor invariants. Let us introduce their definitions for Fano orbifold curves $\mathbb{P}^1_{\bf a}$. For more details on orbifold GW theory we refer to \cite{CR} for the analytic approach and to \cite{AbGV}  for the algebraic geometry approach. 
Let $d\in {\rm Eff}(\mathbb{P}^1_{\bf a})\subset H_2(\mathbb{P}^1_{\bf a};\Z)\cong
\Z$ be an effective curve class. By choosing the homology class
$[\mathbb{P}^1_{\bf a}]$ as a $\Z$-basis of $H_2(\mathbb{P}^1_{\bf a};\Z)$ we may
identify $d$ with a non-negative integer. 
Let $\overline{\mathcal{M}}_{g,n}(\mathbb{P}^1_{\bf a},d)$ be the moduli space of stable orbifold maps $f$ from a genus-$g$ nodal orbifold Riemann surface $\Sigma$ to
$\mathbb{P}^1_{\bf a}$, such that $f_*[\Sigma]=d$. In addition, $\Sigma$ is
equipped with $n$ marked points $z_1,\dots,z_n$ that are pairwise distinct and not
nodal and the orbifold structure of $\Sigma$ is non-trivial only at
the marked points and the nodes.  The moduli space $\overline{\mathcal{M}}_{g,n}(\mathbb{P}^1_{\bf a},d)$ has a virtual fundamental cycle $[\overline{\mathcal{M}}_{g,n}(\mathbb{P}^1_{\bf a},d)]^{\rm virt}$. Its homology degree is 
\beq\label{virt-deg}
2\left((3-\dim \mathbb{P}^1_{\bf a})(g-1)+\chi\cdot d+n\right).
\eeq

The moduli space is naturally
equipped with line bundles $\mathcal{L}_j$ formed by the cotangent
lines\footnote{Here $\bar{\Sigma}$ is the nodal Riemann surface underlying $\Sigma$ and $\bar{z}_j\in\bar{\Sigma}$ is the $i$-th marked point on $\bar{\Sigma}$.} $T^*_{\bar{z}_j}\bar{\Sigma}/{\rm Aut}(\Si,z_1,\dots,z_n;f)$ and with evaluation map 
\ben
{\rm ev}: \overline{\mathcal{M}}_{g,n}(\mathbb{P}^1_{\bf a},d)\rightarrow
\underbrace{{\rm I}\, \mathbb{P}^1_{\bf a}\times\cdots\times
  {\rm I}\,\mathbb{P}^1_{\bf a}}_\text{$n$},
\een
obtained by evaluating $f$ at the (orbifold) marked points $z_1, ..., z_n$ and
landing at the connected component of the inertia orbifold ${\rm I}\,\mathbb{P}^1_{\bf a}$
corresponding to the generator of the automorphism group of the
orbifold point $z_j$ (c.f. \cite{CR}).

The {\em descendant orbifold GW invariants} of $\mathbb{P}^1_{\bf a}$ are intersection numbers 
\beq\label{des-inv}
\langle \phi_{1}\,\psi_1^{k_1},\dots,\phi_{n}\psi_n^{k_n}\rangle_{g,n,d}
:=\int_{[\overline{\mathcal{M}}_{g,n}(\mathbb{P}^1_{\bf a},d)]^{\rm virt}} {\rm ev}^*(\phi_{1}\otimes\cdots \otimes \phi_{n})\,\psi_1^{k_1}\cdots \psi_n^{k_n},
\eeq
where $\phi_{j}\in H:=H_{\rm CR}(\mathbb{P}^1_{\bf a};\C) $, $\psi_{j}=c_1(\mathcal{L}_{j})$. The  {\em total descendant  potential} is
\ben
\cD_{\bf a}(\hbar;\t) = \exp\Big( \sum_{g,n,d} \hbar^{g-1}\,\frac{Q^d}{n!}
\langle \t(\psi_1),\dots,\t(\psi_n)\rangle_{g,n,d}\Big),
\een
where $Q$ is a non-zero complex number called the {\em Novikov variable}, $\hbar$,
$t_0,t_1,\ldots\in H $ are formal variables and $\t(z):=t_0+t_1z+t_2
z^2+\cdots$.

Let $\pi:\overline{\mathcal{M}}_{g,n}(\mathbb{P}^1_{\bf a},d)\to\overline{\mathcal{M}}_{g,n}$ be the forgetful morphism and
\ben
\Lambda_{g,n,d}\left(\phi_{1},\cdots,\phi_{n}\right):=
\pi_*\left([\overline{\mathcal{M}}_{g,n}(\mathbb{P}^1_{\bf a},d)]^{\rm virt}\cap {\rm ev}^*(\phi_{1}\otimes\cdots \otimes \phi_{n})\right).
\een
The {\em ancestor orbifold GW invariants} of $\mathbb{P}^1_{\bf a}$ are intersections numbers over the moduli space of stable curves $\overline{\mathcal{M}}_{g,n}$ ($2g-2+n>0$):
\beq\label{anc-cor}
\langle \phi_{1}\,\bar{\psi}_1^{k_1},\dots,\phi_{n}\bar{\psi}_n^{k_n}\rangle_{g,n,d}
:=\int_{\overline{\mathcal{M}}_{g,n}}\Lambda_{g,n,d}\left(\phi_{1},\cdots,\phi_{n}\right)\,\bar{\psi}_1^{k_1}\cdots \bar{\psi}_n^{k_n}, 
\eeq
where $\bar{\psi}_{j}$ is the $j$-th $\psi$-class over $\overline{\mathcal{M}}_{g,n}$. 
We define the total ancestor potential of $\mathbb{P}^1_{\bf a}$ as follows
\beq\label{total-anc}
\mathcal{A}_{\bf a}(\hbar;\t):= \exp\Big( \sum_{g,n,d} \hbar^{g-1}\,\frac{Q^d}{n!}
\langle \t(\bar{\psi}_1),\dots,\t(\bar{\psi}_n)\rangle_{g,n,d}\Big).
\eeq

For each element $t\in H$, it is useful to introduce the double bracket notation:
\ben
\langle\langle \phi_{1}\,\bar{\psi}_1^{k_1},\dots,\phi_{n}\bar{\psi}_n^{k_n}\rangle\rangle_{g,n}(t):
=\sum_{k,d}\frac{Q^{d}}{k!}\langle \phi_{1}\,\bar{\psi}_1^{k_1},\dots,\phi_{n}\bar{\psi}_n^{k_n}, t, \dots, t\rangle_{g,n+k,d}
\een
We define a total ancestor potential that depends on the choice of $t$,
 \beq\label{total-anc-t}
\mathcal{A}_{t}(\hbar;\t) = \exp\Big( \sum_{g,n} \hbar^{g-1}\,\frac{1}{n!}
\langle\langle\t(\bar{\psi}_1),\dots,\t(\bar{\psi}_n)\rangle\rangle_{g,n}(t)\Big).
\eeq
According to \cite{G3} the total ancestor potential $\mathcal{A}_{t}(\hbar;\t)$ and the total descendant potential $\cD_{\bf a}(\hbar;\t)$ are related by the quantization of a calibration operator $S_{t}(z)$ in Section \ref{sec:cali}. We will explain the details of the quantization in Section \ref{sec:higher}.

The {\em quantum cup product} is a family of associative commutative
multiplications $\bullet_t$ (or just $\bullet$ if the reference point $t$ is mentioned) in $H$ defined for each $t\in H$ via the correlators
\ben
(\phi_i\bullet_t\phi_j,\phi_k) =
\langle\langle \phi_i,\phi_j,\phi_k\rangle\rangle(t).
\een
The degree-0 part of $\bullet_t$ at $t=0$ is called the {\em Chen-Ruan cup product}. We denote it by 
\ben
\cup_{\rm CR}=\bullet_{t=0}|_{Q=0}
\een
Let $t_i$, $i\in \mathfrak{I}$ be the corresponding coordinates of $\phi_i$. 
The quantum cup product induces on $H$ a Frobenius
structure of conformal dimension 1 with respect to the {\em Euler vector field}
\ben
E=\sum_{i\in \mathfrak{I}} d_i t_i\frac{\d}{\d t_i} +
\chi\, \frac{\d}{\d t_{02}}
\een
where $d_i=1-\deg_{\rm CR }(\phi_i)$, i.e., 
\ben
d_{01}= 1,\quad d_{02}=0,\quad 
d_i =1-\frac{p}{a_{k}},\quad i=(k,p)\in\mathfrak{I}_{\rm tw}.
\een

\subsection{Mirror symmetry for the quantum cohomology}
The Frobenius structure on $H$ arising from quantum cohomology can be identified with the Frobenius structure
on a certain deformation space of the {\em affine cusp polynomial}  
\beq\label{aff:poly}
f_{\bf a}(x) = x_1^{a_1}+x_2^{a_2}+x_3^{a_3} - \frac{1}{Q} x_1x_2x_3,\quad x=(x_1,x_2,x_3).
\eeq
where $Q\in\C^*$ is the Novikov variable. The isomorphism in the case
$a_1=1$ was established in \cite{MT} and the general case can be found
in \cite{Ro}. According to Ishibashi--Shiraishi--Takahashi (see \cite{IST}), the Frobenius structure can be
described also in the general framework of K. Saito's theory of
primitive forms.  This is precisely the point of view suitable for our
purposes.

Denote the Milnor number of $f_{\bf a}$ (i.e., the number of critical points of a Morsification of $f_{\bf a}$) by $$N+1=a_1+a_2+a_3-1.$$ Denote the space of a miniversal deformation of the polynomial $f_{\bf a}$ by $$M=\C^{N+1}.$$
Note that the cardinality of the set $\mathfrak{I}$ is $N+1$,
so we can enumerate the coordinates on $M$ via $s=(s_i)_{i\in \mathfrak{I} }$. Recall $\mathfrak{I}_{\rm tw}=\mathfrak{I}\setminus\{(01), (02)\}$. Given $s\in M$, we put
\ben
F(x,s) = x_1^{a_1}+x_2^{a_2}+x_3^{a_3}-
\frac{1}{Qe^{s_{02}}}\,x_1x_2x_3+s_{01} + \sum_{i=(k,p)\in
  \mathfrak{I}_{\rm tw}} s_i\, x_{k}^{p}.
\een
 Let $C\subset M\times \C^3$ be the analytic subvariety with structure sheaf 
\ben
\cO_C = \cO_{M\times \C^3}/(\partial_{x_1}F, \partial_{x_2}F, \partial_{x_3}F);
\een
then the {\em Kodaira-Spencer map} 
\beq\label{KS-iso}
\cT_M\to p_* \cO_C,\quad \frac{\d}{\d s_i} \mapsto \frac{\d F}{\d s_i}\,
{\rm mod}\ (\partial_{x_1}F, \partial_{x_2}F, \partial_{x_3}F),
\eeq
where $p:M\times \C^3\to M$ is the projection onto the first factor, 
is an isomorphism which allows us to define an associative, commutative
multiplication $\bullet$ on $\cT_M$. 
The main result in \cite{IST} is that 
\ben
\omega=\frac{\sqrt{-1}}{Qe^{s_{02}}}\, dx_1\wedge dx_2\wedge dx_3
\een
is a {\em primitive form} in the sense of K. Saito (see \cite{Sa}),
which allows us to construct a Frobenius structure on $M$ (see \cite{SaT}).  More precisely, the form $\omega$ gives rise to a residue
pairing on $\cO_C$
\ben
(\phi_1,\phi_2)=-\frac{1}{Q^2e^{2\, s_{02}} }\, {\rm Res}_{M\times \C^3/M}\, \frac{\phi_1\phi_2 \,
  dx_1\wedge dx_2\wedge dx_3}{ \partial_{x_1}F \partial_{x_2}F \partial_{x_3}F}\, ,
\een
which via the Kodaira--Spencer isomorphism \eqref{KS-iso} induces a non-degenerate bilinear form on $\cT_M$. Let us form the following family of connections  on $\cT_M$
\ben
\nabla = \nabla^{\rm L.C.} - \frac{1}{z} \sum_{i\in \mathfrak{I}}
(\d_{s_i}\bullet )\, ds_i,
\een
where $\nabla^{\rm L.C.} $ is the Levi-Cevita connection associated
with the residue pairing and $\d_{s_i}\bullet$ is the operator of
multiplication by the vector field $\d/\d s_i$. 
Let us also introduce the {\em oscillatory integrals}
\ben
J_\cA(s,z) = (-2\pi z)^{-3/2}\,z d_M\, \int _{\cA_{s,z}} e^{F(x,s)/z}\, \omega \in T^*_s M,
\een
where $d_M$ is the de Rham differential on $M$, and $\cA$ is a flat section
of the bundle on $M\times \C^*$, whose fiber over a point $(s,z)$ is
given by the space of semi-infinite homology cycles
\ben
H_3 (\C^3,\{x|{\rm Re}(F(x,s)/z)\ll 0\};\C)\cong \C^{N+1}.
\een
The fact that $\omega$ is primitive means that the connection $\nabla$
is flat for all $z\neq 0$ and that after identifying $\cT_M\cong
\cT_M^*$ via the residue pairing, the oscillatory integrals $J_\cA$ give rise
to flat sections of $\nabla$. Moreover, since the oscillatory integrals are
weighted-homogeneous functions if one assigns weights $d_i$ ($i\in \mathfrak{I}$), $1/a_j$ ($1\leq j\leq 3$),
and $\chi$ to $s_i$, $x_j$, and $Q$ respectively, they satisfy an additional differential equation with respect to $z$. Let $E\in \cT_M$
be the {\em Euler vector field}
\ben
E=\sum_{i\in \mathfrak{I}} d_i s_i \frac{\d}{\d s_i} + \chi
\frac{\d}{\d s_{02}}.
\een
Note that under the Kodaira--Spencer isomorphism $E$ corresponds to
the equivalence class of $F$ in $p_*\cO_C$. The oscillatory integrals
satisfy the following differential equation:
\beq\label{deq-z}
(z\d_z + E) \, J_\cA(t,z) = \theta\,  J_\cA(t,z),
\eeq
where $\theta:\cT_M\to \cT_M$ is the {\em Hodge grading operator} defined via 
\beq\label{eq:theta}
\theta (X) = \nabla^{L.C.}_X(E) - \frac{1}{2}\, X
\eeq
where the constant $\frac{1}{2}$ is chosen in such a way that $\theta$
is {\em anti-symmetric} with respect to the residue pairing:
$(\theta(X),Y)=-(X,\theta(Y))$. 

The quantum cohomology computed at $t=0$ is isomorphic as a Frobenius algebra with $T_0M$ (see \cite{IST, Ro}). The
identification has the following form  
$$
\phi_i = x_{k}^{p} +\cdots,\quad  
\phi_{01}=1,\quad \phi_{02}= \frac{1}{Q} x_1x_2x_3+\cdots.
$$
where $i=(k,p)$ is the index of a twisted class and the dots stand for some 
polynomials that involve higher-order powers of $Q$. More precisely, using
the Kodaira-Spencer isomorphism we have
\ben
\phi_i= \d_{s_i}+\cdots,\quad \phi_{01}=\d_{s_{01}},\quad \phi_{02}=\d_{s_{02}}+\cdots,
\een  
where the dots stand for some vector fields depending holomorphically on $Q$ near $Q=0$
and vanishing at $Q=0.$ These additional terms are uniquely fixed by the requirement
that the vector fields $\phi_i$ $(i\in \mathfrak{I})$ are flat, i.e., the residue pairing is constant independent of $Q$. 
On the other hand the flatness of $\nabla$ implies that the residue
pairing is flat, therefore we can extend uniquely the isomorphism
$H\cong T_0M$ to an isomorphism $$TH\cong TM$$ such that the residue
pairing coincides with the Poincar\'e pairing. In other words, the
linear coordinates $t_i$, $i\in \mathfrak{I}$ on $H$ are functions on
$M$ such that $t_i(0)=0$, the vector field $\d/\d t_i$ is flat with respect to the
Levi--Civita connection, and at $s=0$ it coincides with
$\phi_i.$ The mirror symmetry for quantum cohomology
can be stated as follows.

\begin{theorem}[\cite{IST}, Theorem 4.1]\label{ms}
The isomorphism $M\cong H,\ s\mapsto t(s)$ is an isomorphism of Frobenius manifolds, i.e., $T_sM\cong T_{t(s)} H$ as Frobenius
algebras.   
\end{theorem}
\begin{remark}
Theorem \ref{ms} can be proved also by using the extended $J$-function of $\mathbb{P}^1_{\bf a}$ (see Section \ref{sec:JX}). Namely, it is not hard to derive an identification between the quantum cohomology $D$-module of $\mathbb{P}^1_{\bf a}$ and the $D$-module defined by $f_{\bf a}(x)$.
\end{remark}

From now on we will make use of the residue pairing to identify
$T^*M\cong TM$. Also the flat Levi--Civita connection $\nabla^{\rm
  L.C.}$ allows us to construct a trivialization $$TM\cong M\times T_0M,$$ and finally, the Kodaira--Spencer map \eqref{KS-iso} together
with the mirror symmetry isomorphism gives $T_0M\cong H$. In other
words, we have natural trivializations
\beq\label{triv}
T^*M\cong TM\cong M\times H. 
\eeq

\subsection{The period integrals and the calibration operator}\label{sec:periods}
Givental noticed that certain period integrals (c.f. formula \eqref{per} below) in singularity theory
play a crucial role in the theory of integrable systems. In this
section, we recall Givental's construction as well as some of its
basic properties. See \cite{G1} for more details.

Put $X=M\times \C^3$ and let 
\ben
\varphi: X \to M\times \C,\quad (s,x)\mapsto (s,F(x,s)).
\een
Let $$X_{s,\la} = \varphi^{-1}(s,\la)$$ be the fibers of $\varphi$. The
  set of all $(s,\la)\in M\times \C$ such that the fiber $X_{s,\la}$ is
  singular is an analytic hypersurface, called {\em discriminant}. Its complement in $M\times
  \C$ will be denoted by $(M\times \C)'.$ The homology and
  cohomology groups $H_2(X_{s,\la};\C)$ and $H^2(X_{s,\la};\C)$,
  $(s,\la)\in (M\times \C)'$ form vector
  bundles over the base $(M\times \C)'$. Moreover, the integral
  structure in the fibers allows us to define a flat connection known
  as the {\em Gauss--Manin} connection. 

Let us fix the point $(0,1)\in (M\times \C)'$ (for $Q\ll 1$) to be our reference point. The vector space $$\lieh=H_2(X_{0,1};\C)$$ has a very
rich structure, which we would like to recall. Let $$\Delta\subset \lieh$$ be the set of {\em vanishing cycles}, and $(\cdot |\cdot )$
be the {\em negative} of the intersection pairing. The negative sign
is chosen so that $(\alpha|\alpha)=2$ for all $\alpha\in \Delta$. The
parallel transport with respect to the Gauss--Manin
connection induces a monodromy representation
\ben
\pi_1((M\times \C)') \to \operatorname{GL}(\lieh).
\een 
The image $$W\subset \operatorname{GL}(\lieh)$$ of the fundamental group under this representation is a subgroup of the group of
linear transformations of $\lieh$ that preserve the intersection
form. The Picard--Lefschetz theory can be applied in our setting as
well and $W$ is in fact a reflection group generated by the
reflections
\ben
s_\al(x) = x-(\alpha|x)\alpha,\quad \alpha\in \Delta.
\een
The reflection $s_\al$ is the monodromy transformation along a simple loop that goes
around a generic point on the discriminant over which the cycle
$\alpha$ vanishes.  
Finally, recall that the {\em classical} monodromy $\si\in W$ 
is the monodromy transformation along a big loop around the discriminant. 
For more details on vanishing homology and
cohomology and the Picard--Lefschetz theory we refer to the book
\cite{AGV}. We will see in Proposition \ref{affine} below that $\Delta$ is an affine root system.

The main objects in our construction are the following multi-valued
analytic functions:
\beq\label{per}
I^{(n)}_\alpha(t,\la) = -\frac{1}{2\pi}\,\d_\la^{n+1}\, d_M\
\int_{\alpha_{t,\la}} d^{-1}\omega,
\eeq
where the value of the RHS depends on the choice of a path avoiding
the discriminant, connecting the reference point with $(t,\la)$. The
cycle $\alpha_{t,\la}$ is obtained from $\alpha\in \lieh$ via a
parallel transport (along the chosen path), $d^{-1}\omega$ is any
holomorphic 2-form $\eta$ on $\C^3$ such that $\omega=d\eta$, and $d_M$
is the de Rham differential on $M$. The RHS in \eqref{per} defines
naturally a cotangent vector in $T^*_tM$, which via the trivialization
\eqref{triv} is identified with a  vector in $H$. 

The period vectors \eqref{per} are uniquely defined for all $n\geq
-1$. For $n\leq -2$ there is an ambiguity in choosing integration
constants, which can be removed by means of the following
differential equations:
\beqa\label{2str-con-1}
\d_{t_i} \, I^{(n)}_\al(t,\la) & = & -\phi_i\bullet
I^{(n+1)}_\al(t,\la) ,\quad i\in \mathfrak{I},  
\\ \label{2str-con-2}
\d_\la \, I^{(n)}_\al(t,\la) & = & 
I^{(n+1)}_\al(t,\la),  
\\ \label{2str-con-3}
(\la-E\bullet)\d_\la I^{(n)}_\al(t,\la) & = & \Big(\theta-n-1/2 \Big) I^{(n)}_\al(t,\la) .
\eeqa
Finally, note that the unit vector
$\one\in H\cong M$ has coordinates $t_{01}=1$, $t_i=0$ for $i\neq (01)$ and that the period vectors have the 
following translation symmetry:
\ben
I^{(n)}_\al(t,\la) = I^{(n)}_\al(t-\la\one,0),\quad \forall n\in \Z,\quad \forall \al\in \lieh.
\een
The oscillatory integrals are related to the period integrals via a
Laplace transform along an appropriately chosen path:
\beq\label{lt}
J_\cA(t,z)=(-2\pi z)^{-1/2} \int_{u_j}^\infty e^{\la/z} I^{(0)}_\al(t,\la)d\la, 
\eeq
where $u_j(t)$ is such that $(t,u_j(t))$ is a point on the
discriminant over which the cycle $\alpha$ vanishes. The differential
equations \eqref{2str-con-1} are the Laplace transform of $\nabla
J_\cA=0$, while the equation \eqref{2str-con-3} is the Laplace
transform of the differential equation \eqref{deq-z}. Using equations
\eqref{2str-con-2} and \eqref{2str-con-3} we can express $I^{(n)}$ in
terms of $I^{(n+1)}$ as long as the operator $\theta-n-1/2$ is
invertible. This is the case for $n\leq -2$, which allows us to extend
the definition of $I^{(n)}$ to all $n\in \Z$.


\subsubsection{Stationary phase asymptotic}\label{sec:spa}
Let $u_j(t)$, $1\leq j\leq N+1$ be the critical values of
$F(x,t)$. The set $$M_{\rm ss}\subset M$$ of all points $t\in M$ such that the
critical values $u_j(t)$ form locally
near $t$ a coordinate system is open and dense. Let us fix some
$t_0\in M_{\rm ss}$; then in a neighborhood of $t_0$ the critical
values give rise to a coordinate system in which the pairing and the
product $\bullet$ are diagonal, i.e.,
\ben
\d/\d u_j \bullet \d/\d u_{j'} = \delta_{j,j'} \d/\d u_j,\quad (\d/\d u_j,\d/\d u_{j'})=\delta_{j,j'}/\Delta_j,
\een
where $\Delta_j$ are some multi-valued analytic functions on $M_{\rm
  ss}.$ Following Dubrovin's terminology (see \cite{Du}), we refer to $u_j$
as {\em canonical coordinates}.
\begin{remark}
It is easy to see that the critical variety $C$ of the function $F$ is
non-singular, i.e., it is a manifold. It can be proved that
the projection map $p:C\subset M\times \C^3 \to M$ is a finite
branched covering of degree $N+1$. The branching points are precisely
$M\setminus{M_{ss}}$.
\end{remark}

Using the canonical coordinates we can construct a trivialization of
the tangent bundle 
\ben
\Psi: M_0\times \C^{N+1}\cong TM_0,\quad (t,e_j)\mapsto (t,\sqrt{\Delta_j} \,
\frac{\d}{\d u_j}).
\een
Here $M_0\subset M_{\rm ss}$ is an open contractible neighborhood of
$t_0$ and $\{e_j\}$ is the standard basis of $\mathbb{C}^{N+1}$, where the $j$-th component of $e_j$ is 1, while the remaining ones are 0. According to Givental (see \cite{G2}), there exists a 
unique formal asymptotic series $\Psi_t R_t(z) e^{U_t/z}$
that satisfies the same
differential equations as the oscillatory integrals $J_\cA$, where 
\beq\label{R}
R_t(z)=1+\sum_{\ell=1}^{\infty}R_\ell(t) z^\ell, \quad R_\ell(t)\in {\rm End}(\C^{N+1}).
\eeq

We will make use of the following formal series
\beq\label{falpha}
\f_\al(t,\la;z) = \sum_{n\in \Z} I^{(n)}_\al(t,\la) \, (-z)^n \,
,\al\in \lieh.
\eeq
\begin{example}
Note that for $A_1$-singularity $F(t,x)=x^2/2+t$ we have
$u:=u_1(t)=t.$ Up to a sign there is a unique vanishing cycle. The
series \eqref{falpha}  will be
denoted simply by $\f_{A_1}(t,\la;z).$ The corresponding
period vectors can be computed explicitly:
\ben
\begin{aligned}
I^{(n)}_{A_1}(u,\la) & = (-1)^n\, \frac{(2n-1)!!}{2^{n-1/2}}\,
(\la-u)^{-n-1/2},\quad n\geq 0 \\
I^{(-n-1)}_{A_1}(u,\la) & = 2\, \frac{2^{n+1/2}}{(2n+1)!!}\,
(\la-u)^{n+1/2}, \quad n\geq 0.
\end{aligned}
\een
\end{example}

The key lemma (see \cite{G1}) is the following.
\begin{lemma}\label{vanishing_a1}
Let\/ $t\in M_{\rm ss}$ and\/ $\be$ be a vanishing cycle vanishing
over the point\/ $(t,u_j(t))$. Then for all\/ $\la$ near\/ $u_j:=u_j(t)$, we have
\ben
\f_{\be}(t,\la;z) = \Psi_t R_t(z)\,  e_j\,  \f_{A_1}(u_j,\la;z)\,.
\een
\end{lemma}
An important corollary of Lemma \ref{vanishing_a1} is the following
remarkable formula due to K. Saito (\cite{Sa}):
\beq\label{int-form}
(\al|\be)=(I^{(0)}_\al(t,\la),(\la-E\bullet)I^{(0)}_\beta(t,\la)).
\eeq
To prove this formula, first note that the the differential equations
\eqref{2str-con-1}--\eqref{2str-con-3} imply that the RHS is independent of $t$
and $\la.$ In order to compute the RHS, let us fix $t\in M_{\rm ss}$
and let $\la$ approach one of the critical values $u_j(t)$ in such a
way that the cycle $\beta$ vanishes over $(t,u_j(t))$. According to
Lemma \ref{vanishing_a1} we have
\ben
I^{(0)}_\be(t,\la) = 2(2(\la-u_j))^{-1/2} e_j + O((\la-u_j)^{1/2}).
\een
Similarly, decomposing $\al=\al'+(\al|\be)\be/2$, where $\al'$ is
invariant with respect to the local monodromy, we get
\ben
I^{(0)}_\al(t,\la) = (\al|\beta)\, (2(\la-u_j))^{-1/2} e_j+ O((\la-u_j)^{1/2}).
\een
It is well known (see \cite{Du}) that in canonical coordinates the
Euler vector field has the form $E=\sum u_j\d_{u_j}$. Now it is easy
to see that the RHS of  \eqref{int-form}, up to higher order terms in
$(\la-u_j)$ is $(\al|\beta)$ and since the latter must be independent
of $\la$ the higher-order terms must vanish.   

\subsubsection{The calibration operator}\label{sec:cali}
The {\em calibration} of the Frobenius structure on $H$ is by definition a
gauge transformation $S$ of the form
\beq\label{gauge}
S_t(z)=1+\sum_{\ell=1}^\infty S_\ell(t)z^{-\ell},\quad S_\ell(t)\in {\rm End}(H),
\eeq
such that $\nabla =S d S^{-1}$. In GW theory there is a canonical choice of calibration given by genus-0 descendant invariants as follows (see \cite{G3}): 
\beq\label{S:def}
(S_t(z)\phi_i,\phi_j) = (\phi_i,\phi_j) + \sum_{\ell=0}^\infty \langle\langle \phi_i\psi^\ell,\phi_j\rangle\rangle_{0,2}(t)z^{-\ell-1}. 
\eeq
Here 
\begin{equation*}
\langle\langle \phi_i\psi^\ell,\phi_j\rangle\rangle_{0,2}(t)=\sum_{m\geq0}\sum_{d\geq0}\frac{Q^d}{m!} \langle \phi_i\psi^\ell,\phi_j, t,...,t\rangle_{0,2+m,d}.
\end{equation*}

It is a general fact in GW theory (see \cite{G3}) that 
\beq\label{cal-z}
S_t(z)^{-1} \Big(\d_z - z^{-1} \theta +z^{-2} E\bullet\Big) S_t(z) =\d_z - z^{-1}\theta + z^{-2} \rho, 
\eeq
where $\rho=\chi\,P\cup_{\rm CR}$. By definition the operator $\rho$ acts on $H$ as follows
\beq\label{eq:rho}
\rho(\phi_{01}) = \chi \phi_{02},\quad \rho(\phi_i) = 0,\quad \mbox{for } i\in\mathfrak{I}\backslash\{(01)\}.
\eeq

We define a new series 
\beq\label{eqn:calib_lim}
\widetilde{\f}_\al(\la;z) := S_t(z)^{-1} \, \f_\al(t,\la;z).
\eeq
Note that the RHS is independent of $t$. Put 
\beq\label{eqn:I^n}
\widetilde{\f}_\al(\la;z) = \sum_{n\in \Z} \widetilde{I}^{(n)}_\al(\la)\,(-z)^n.
\eeq
We will refer to $\widetilde{I}^{(n)}_\al(\la)$ as
the {\em calibrated} limit of the period vector $I^{(n)}_\al(t,\la).$

In our general set up the Novikov variable $Q$ is a fixed non-zero
constant. However, it will be useful also to allow $Q$ to vary in a
small contractible neighborhood and to study the dependence
of the periods and their calibrated limits on $Q$. By definition
$I^{(n)}_\al(t,\la)$  depend on $Qe^{t_{02}}$, so we simply have
\ben
Q\d_Q\, I^{(n)}_\al(t,\la) = \d_{t_{02}}\, I^{(n)}_\al(t,\la). 
\een
Using the {\em divisor equation} in GW theory, it is easy
to prove (c.f. \cite{G3}) that the gauge transformation $S_t(z)$ satisfies the following
differential equation:
\ben
zQ\d_Q\, S_t(z) = z \d_{t_{02}} S_t(z) - S_t(z)\, (P\,\cup_{\rm CR}\,).
\een
Finally, the gauge identity $\nabla = SdS^{-1}$ and the differential equations
\eqref{2str-con-1}--\eqref{2str-con-3} imply that 
the calibrated limit of the period vectors satisfy the following
system of differential equations: 
\beqa \label{calp-div}
Q\d_Q\tI^{(n)}_\al(\la) & = &  -P\,\cup_{\rm CR}\, \tI^{(n+1)}_\al(\la) 
\\ \label{calp-1}
\d_\la \, \tI^{(n)}_\al(\la) & = & 
\widetilde{I}^{(n+1)}_\al(\la),  
\\ \label{calp-2}
(\la-\rho)\d_\la \widetilde{I}^{(n)}_\al(\la) & = & \Big(\theta-n-1/2 \Big) \widetilde{I}^{(n)}_\al(\la) .
\eeqa

\begin{lemma}\label{per-1}
a) Let $\{B_i\}_{i\in \mathfrak{I}}$ be a basis of $\lieh^\vee:=H^2(X_{0,1};\C)$, then the following formula holds
\ben
\tI^{(-1)}_\al(\la)= \langle B_{01},\al\rangle\,
\Big( \la \one+  (\chi\log\la-\log Q)\,P \Big)  + 
\langle B_{02},\al\rangle\, P + 
\sum_{i\in \mathfrak{I}_{\rm tw}} \langle B_i,\al\rangle\, \la^{d_i}\,\phi_i,
\een

b) The analytic continuation of 
$\tI^{(n)}_\al(\la)$ along a closed loop around $0$ is $\tI^{(n)}_{\si(\al)}(\la)$, where $\si$ is the classical monodromy.
\end{lemma}
\proof
a)
Recall $\rho$ acts on $H$ by \eqref{eq:rho},
while the operator $\theta$ defined in \eqref{eq:theta} has the form (via \eqref{triv})
\ben
\theta(\phi_i) = (d_i-1/2)\,\phi_i ,\quad i\in \mathfrak{I}. 
\een
Note that the $H$-valued functions that follow the
pairings $\langle B_i,\al\rangle$ are solutions to the system
\eqref{calp-div}--\eqref{calp-2} with $n=-1$.  These solutions are
linearly independent, therefore they must give a basis in the space of
all solutions. 

b) Now the statement follows, because it is true for $I^{(n)}_\al(t,\la)$, for $|\la|\gg1$, where
\ben
I^{(n)}_\al(t,\la)=\widetilde{I}^{(n)}_\al(\la) + \sum_{\ell=1}^\infty (-1)^\ell\,S_\ell(t)\widetilde{I}^{(n+\ell)}_\al(\la).
\een
 
\qed

\subsection{Mirror symmetry at higher genus}\label{sec:higher}
A Frobenius manifold is called {\em semi-simple} if the multiplication has a semi-simple basis.
The Frobenius manifold $(H,(\ ,\ ),\bullet_t,\phi_{01},E)$ is isomorphic to the Frobenius manifold constructed from the mirror model of $\mathbb{P}^1_{\bf a}$ \cite{MT, Ro, IST}, see Theorem \ref{ms}. Using the mirror model, it is easy to see that $\bullet_t$ is semi-simple for generic $t$. 

For any semi-simple Frobenius manifold, Givental introduced a higher genus reconstruction formula \cite{G2} using the symplectic loop space formalism \cite{G3}. Furthermore, he conjectured that the higher genus GW ancestor invariants are uniquely determined from its semi-simple quantum cohomology.
Teleman \cite{Te} has proved this conjecture.  Let us recall the construction.
\subsubsection{Canonical quantization}\label{sec:can_quantize}
Equip the space $$\cH:=H(\!(z^{-1})\!)$$ of formal Laurent series in $z^{-1}$ with
coefficients in $H$  with the following \emph{symplectic form}: 
\ben
\Omega(\phi_1(z),\phi_2(z)):={\rm Res}_z \left( \phi_1(-z),\phi_2(z)\right) \,,
\qquad \phi_1(z),\phi_2(z)\in\cH \,,
\een 
where, as before, $(,)$ denotes the residue pairing on $H$
and the formal residue ${\rm Res}_z$ gives the coefficient in front of $z^{-1}$.

Let $\{\phi_i\}_{i\in \mathfrak{I}}$ and $\{\phi^i\}_{i\in\mathfrak{I}}$ be dual bases of $H$ with respect to the residue pairing.
Then
\ben
\Om(\phi^i(-z)^{-\ell-1}, \phi_j z^m) = \de_{ij} \de_{\ell m} \,.
\een
Hence, a Darboux coordinate system is provided by the linear functions $q_\ell^i$, $p_{\ell,i}$ on $\cH$ given by:
\ben
q_\ell^i = \Om(\phi^i(-z)^{-\ell-1}, \cdot) \,, \qquad
p_{\ell,i} = \Om(\cdot, \phi_i z^\ell) \,.
\een
In other words,
\ben
\phi(z) = \sum_{\ell=0}^\infty \sum_{i\in\mathfrak{I}} q_\ell^i(\phi(z)) \phi_i z^\ell +  
\sum_{\ell=0}^\infty \sum_{i\in\mathfrak{I}} p_{\ell,i}(\phi(z)) \phi^i(-z)^{-\ell-1} \,,
\qquad \phi(z)\in\cH \,.
\een  
The first of the above sums will be denoted by $\phi(z)_+$ and the second by $\phi(z)_-$.

The \emph{quantization} of linear functions on $\cH$ is given by the rules
\ben
\widehat q_\ell^{\,i} = \hbar^{-1/2} q_\ell^i \,, \qquad
\widehat p_{\ell,i} = \hbar^{1/2} \frac{\d}{\d q_\ell^i} \,,
\een
where the RHSs of the above definitions are operators acting on the Fock space
\beq\label{eq:fock}
\C_\hbar[\![\q]\!]:=\C_{\hbar}[\![q_0,q_1+\one,q_2,\cdots]\!],
\quad
\mbox{where}
\quad
\C_{\hbar} := \C(\!(\hbar)\!)
\quad q_\ell:=(q_\ell^i)_{i\in \lieI}.
\eeq 
Every $\phi(z)\in\cH$ gives rise to the linear function $\Om(\phi(z),\cdot)$ on $\cH$,
so we can define the quantization $\widehat{\phi(z)}$. Explicitly,
\beq\label{phihat}
(\phi_i z^\ell)\sphat = -\hbar^{1/2}\frac{\d}{\d q_\ell^i} ,\quad 
(\phi^i(-z)^{-\ell-1})\sphat = \hbar^{-1/2} \,q_\ell^i.
\eeq
The quantization also makes sense for $\phi(z)\in H[[z,z^{-1}]]$ if we interpret $\widehat{\phi(z)}$
as a formal differential operator in the variables $q_\ell^i$ with
coefficients in $\C_{\hbar}$.
\begin{lemma}\label{lphihat}
For all\/ $\phi_1(z),\phi_2(z)\in\cH$, we have\/ $[\widehat{\phi_1(z)}, \widehat{\phi_2(z)}] = \Om(\phi_1(z),\phi_2(z))$.
\end{lemma}
\proof
It is enough to check this for the basis vectors $\phi^i(-z)^{-\ell-1}$, $\phi_i z^\ell$, in which case it is true by definition.
\qed

\subsubsection{Quantization of symplectic transformations}\label{hgr}

It is known that both series $S_t(z)$ and $R_t(z)$ described in Sections\ \ref{sec:spa} and \ref{sec:cali} are symplectic
transformations on $(\cH, \Omega)$. Moreover, they both have the form $e^{A(z)},$ where
$A(z)$ is an infinitesimal symplectic transformation.  

A linear operator $A(z)$ on $\cH:=H(\!(z^{-1})\!)$ is infinitesimal
symplectic if and only if the map $\cH\ni \phi(z) \mapsto A(\phi(z))\in \cH$ is
a Hamiltonian vector field with a Hamiltonian given by the quadratic
function $$h_A(\phi(z)) = \frac{1}{2}\Omega(A(\phi(z)),\phi(z)).$$ 
By definition, the \emph{quantization} of $e^{A(z)}$ is given by the
differential operator $e^{\widehat{h}_A},$ where the quadratic
Hamiltonians are quantized according to the following rules: 
\ben
(p_{\ell,i}p_{m,j})\sphat = \hbar\frac{\d^2}{\d q_\ell^i\d q_m^j} \,,\quad 
(p_{\ell,i}q_m^j)\sphat = (q_m^jp_{\ell,i})\sphat = q_m^j\frac{\d}{\d q_\ell^i} \,,\quad
(q_\ell^iq_m^j)\sphat = \frac1{\hbar} q_\ell^iq_m^j \, .
\een     

In the case of the orbifold $\mathbb{P}^1_{\bf a}$, the Frobenius manifold is semi-simple at a generic point $t\in H$. Teleman's higher genus reconstruction theorem \cite{Te} implies that the total ancestor potential defined in \eqref{total-anc-t} can be identified with Givental's higher genus reconstruction formula \cite{G3}
\beq\label{ancestor}
\cA_t(\hbar;\q(z)) =\widehat{\Psi_t  R_t\Psi_t^{-1}}\ \prod_{j=1}^{N+1}\,
\cD_{\rm pt}(\hbar\Delta_j;\leftexp{j}{\q}(z)\sqrt{\Delta_j}) \in
\C_{\hbar,Q}[[q_0,q_1+\one,q_2\dots]] 
\eeq
and the total descendant potentials defined in \eqref{total-desc} can be identified with
\beq\label{dxna}
\cD_{\bf a}(\hbar;\q(z))=e^{F^{(1)}(t)} \widehat{S}_t^{-1} \cA_t(\hbar;\q(z)) \,,
\eeq
where $\leftexp{j}{\q}(z):=\sum_{\ell=0}^\infty \leftexp{j}{q}_\ell z^\ell$ and
the coefficents $\leftexp{j}{q_\ell}$ are defined by
\ben
\sum_{j=1}^{N+1}\leftexp{j}{q_\ell}\Psi(e_j) =
\sum_{i\in \lieI}  q_\ell^i\phi_i  \,.
\een
Recall that $\cD_{\rm pt}$ is the total descendant potential of a point and the factor 
\ben
F^{(1)}(t)= \sum_{d,n=0}^\infty \frac{Q^d}{n!}\langle t,\ldots,t\rangle_{1,n,d}\,
\een
is the genus-$1$ primary (i.e. no descendants) potential. 
Let us examine more carefully the quantized action of the operators in
formula \eqref{ancestor} and \eqref{dxna}.

\subsubsection{The action of the asymptotical operator}\label{asop}
The operator ${\widehat{U_t/z}}$ is known to annihilate the
Witten--Kontsevich tau-function. Therefore, 
$e^{\widehat{U_t/z}}$ is redundant and it can be dropped from the
formula. The action of the operator $\widehat{R}_t$ on formal functions, whenever it makes sense, is given as follows.
\begin{lemma}[Givental \cite{G2}]\label{lrfock}
We have
\ben
\widehat{R}^{-1}_t\,F(\q) = \left.\Big(e^{\frac{\hbar}{2}V_t(\d,\d)}F(\q)\Big)\right|_{\q\mapsto R_t\q} \,,
\een
where\/ $V_t(\d, \d)$ is the quadratic differential operator 
\ben
V_t(\d,\d) = \sum_{\ell, m=0}^\infty \sum_{i,j\in \lieI} (\phi^i,V_{\ell m}(t)\phi^j) \, \frac{\d^2}{ \d q_\ell^i \d q_m^j}
\een
whose coefficients\/ $V_{\ell m}(t)$ are given by 
\ben
\sum_{\ell, m=0}^\infty V_{\ell m}(t) z^\ell w^m = \frac{1-R_t(z)(\leftexp{T}{R}_t(w))}{z+w} 
\een
and\/ $\leftexp{T}{R}_t(w)$ denotes the transpose of\/ $R_t(w)$ with respect to the Poincare pairing.
\end{lemma}

The substitution $\q\mapsto R_t\q $ can be written more explicitly as follows:
\ben
q_0\mapsto q_0,\quad q_1\mapsto R_1(t)q_0 + q_1,\quad q_2\mapsto R_2(t)q_0 +R_1(t)q_1 + q_2 \,,\dots.
\een
The above substitution is not a well-defined
operation on the space of formal functions. This complication,
however, is offset by a certain property of the
Witten--Kontsevich tau-function, which we now explain.
By definition, an {\em asymptotical function} is a formal function  of the type:
\ben
\cA(\q) = \exp \Big(\sum_{g=0}^\infty\hbar^{g-1} F^{(g)}(\q)\Big) \,.
\een
Such a function is called {\em tame} if the following {\em $(3g-3+n)$-jet constraints} are satisfied:
\ben
\frac{\d^{n} F^{(g)}}{\d q_{k_1}^{i_1}\cdots \d
  q_{k_n}^{i_n}}\Bigg|_{\q = 0} = 0 \quad \text{ if } \quad k_1+\cdots
+ k_n > 3g-3+n \,. 
\een   
The Witten--Kontsevich tau-function (up to the shift $q_1\mapsto
q_1+1$) is tame for dimensional reasons: $\dim \overline{\mathcal{M}}_{g,n}=3g-3+n$. The total ancestor potential $\cA_t$
is also tame, as it can be seen from its geometric definition
(cf. \cite{G3}) or by using the fact that  the action of the operator
$\widehat{R}_t$ on tame functions is well defined and it
preserves the tameness property (\cite{G1}).

\subsubsection{The action of the calibration}\label{cali}
The quantized symplectic transformation $\widehat{S}^{-1}_t$ acts on formal functions as follows.
\begin{lemma}[Givental \cite{G2}]\label{lsfock}
We have
\beq\label{S:fock}
\widehat{S}^{-1}_t\, F(\q) = e^{\frac1{2\hbar} W_t(\q,\q)}F\bigl((S_t\q)_+\bigr) \,,
\eeq
where\/ $W_t(\q,\q)$ is the quadratic form 
$$
W_t(\q,\q) = \sum_{\ell, m=0}^\infty (W_{\ell m}(t)q_m,q_\ell)
$$ 
whose coefficients are defined by
\ben
\sum_{\ell, m=0}^\infty W_{\ell m}(t) z^{-\ell}w^{-m}=\frac{\leftexp{T}{S}_t(z)S_t(w)-1}{z^{-1}+w^{-1}} \,.
\een
\end{lemma}
The subscript $+$ in \eqref{S:fock} means truncation of all negative powers
of $z$, i.e., in $F(\q)$ we have to substitute (cf.\ \eqref{gauge}): 
\ben
q_\ell\mapsto q_\ell + S_1(t) q_{\ell+1}+S_2(t) q_{\ell+2}+\cdots \,,\quad \ell=0,1,2,\dots\ .
\een
This operation is well-defined on the space of formal power series.

\section{$\Gamma$-integral structures and the root system}\label{sec:int_GW}

If $X$ is a compact complex orbifold, then using the $K$-ring $K(X)$
of orbifold vector bundles on $X$ and a certain $\Gamma$-modification
of the Chern character map, Iritani has introduced an integral lattice
in the Chen-Ruan cohomology group $H_{\rm CR}(X;\C)$ (see \cite{Ir}
and also \cite{KKP}). If $X$ has semi-simple quantum cohomology, then
it is expected that $X$ has a LG mirror model and it is natural to
conjecture that Iritani's embedding of the $K$-theoretic lattice
coincide with the image of the Milnor lattice via an appropriate
period map. In our case, when $X=\mathbb{P}^1_{\bf a}$, we prove the
above conjecture by using the same argument as in \cite{Ir}, where the
toric case was proved.  Moreover, we obtain an explicit identification
of the set of vanishing cycles with a certain $K$-theoretic affine
root system.

\subsection{Iritani's integral structure and mirror symmetry}\label{ir-lattice}
Let us recall Iritani's construction in the most general case when $X$ is a compact complex orbifold. Let $IX$ be the {\em inertia orbifold} of $X$, i.e., as a groupoid the points of $IX$ are 
\ben
(IX)_0=\{(x,g)\ |\ x\in X_0,\  g\in \operatorname{Aut}(x)\}
\een
while the arrows from $(x',g')$ to $(x'',g'')$ consists of all arrows $g\in X_1$ from $x'$ to $x''$, s.t., 
$g''\circ g = g\circ g'$. It is known that $IX$ is an orbifold consisting of several connected components $X_v$, $v\in T:=\pi_0(|IX|)$. Following Iritani, we define a linear map 
\ben
\Psi: K(X)\to H^*(IX;\C)=\bigoplus_{v\in T} H^*(X_v;\C)
\een
via 
\beq\label{psi-map}
\Psi(V)=(2\pi)^{-\operatorname{dim}_\C X/2}\ 
\widehat{\Gamma}(X) \cup (2\pi\sqrt{-1})^{\operatorname{deg}}
\operatorname{inv}^* \widetilde{\operatorname{ch}}(V).
\eeq
Here $\cup$ is the usual cup product in $H^*(IX;\C)$. Let us recall the notation. The linear operator 
\ben
\operatorname{deg}: H^*(IX;\C)\to H^*(IX;\C)
\een 
is defined by $\operatorname{deg} (\phi)=r \phi$ if $\phi \in H^{2r}(IX;\C)$. 
The involution $\operatorname{inv}:IX\to IX$ inverts all arrows while on the points it acts as $(x,g)\mapsto (x,g^{-1}).$ If $V$ is an orbifold vector bundle, then we have an eigenbasis decomposition
\ben
\operatorname{pr}^*(V) = \bigoplus_{v\in T} V_v=\bigoplus_{v\in T}\bigoplus_{0\leq f<1} V_{v,f},
\een
where $\operatorname{pr}:IX\to X$ is the forgetful map $(x,g)\mapsto x$ and $V_{v,f}$ is the subbundle of $V_v:=\operatorname{pr}^*(V)|_{X_v}$ whose fiber over a point $(x,g)\in (IX)_0$ is the eigenspace of $g$ corresponding to the eigenvalue $e^{2\pi\sqrt{-1} f}$. Let us denote by $\delta_{v,f,j}$ $(1\leq j\leq l_{v,f}:={\rm rk}(V_{v,f}))$ the Chern roots of $V_{v,f}$, then the Chern character and the $\Gamma$-class of $V$ are defined by
\ben
\widetilde{\operatorname{ch}}(V) = \sum_{v\in T} \sum_{0\leq f<1} e^{2\pi\sqrt{-1} f}
\operatorname{ch}(V_{v,f}),
\een
\ben
\widehat{\Gamma}(V) = \sum_{v\in T} \prod_{0\leq f<1} \prod_{j=1}^{l_{v,f}} \Gamma(1-f+\delta_{v,f,j}),
\een
where the value of the $\Gamma$-function $\Gamma(1-f+y)$  at $y=\delta_{v,f,j}$ is obtained by first expanding in Taylor's series at $y=0$ and then formally substituting $y=\delta_{v,f,j}$. By definition $\widehat{\Gamma}(X):=\widehat{\Gamma}(TX).$

\subsubsection{The $\Gamma$-conjecture for the Milnor lattice}
We denote by $H_{\rm CR}(X;\C)$ the vector space  $H^*(IX;\C)$ equipped with the Chen--Ruan cup product $\cup_{\rm CR}$. We define a shift function $\iota: T \to \mathbb{Q}$ by
\ben
\iota(v) = \sum_{0\leq f<1} f\operatorname{dim}_\C(TX)_{v,f}.
\een
The Chen--Ruan product is graded homogeneous with respect to the following grading 
\ben
\operatorname{deg}_{\rm CR}(\phi) = (r+\iota(v))\phi,\quad \phi\in H^{2r}(X_v;\C). 
\een
The vector space $H^*(IX;\C)$ is equipped with a Poincar\'e pairing, i.e.
\ben
(\phi_1,\phi_2)=\int_{IX} \phi_1\cup {\rm inv}^*(\phi_2).
\een 
This pairing turns both algebras $H^*(IX;\C)$ and $H_{\rm CR}(X;\C)$ into Frobenius algebras. 
Let us point out also that by using the Kawasaki Riemann--Roch formula we can also prove that the map $\Psi$ is compatible (up to a sign) with the natural pairing on $K(X)$ and the Poincar\'e pairing
\beq\label{K-pairing}
\chi(V_1\otimes V_2^\vee) = (e^{\pi \sqrt{-1} \theta_X} e^{\pi\sqrt{-1} \rho_X} \Psi(V_1),\Psi(V_2)),
\eeq
where $\rho_X=c_1(TX)\cup_{\rm CR}$ and $\theta_X$ is the {\em Hodge grading operator} of $X$,
\ben
\theta_X = \frac{1}{2}\,\operatorname{dim}_\C X -\operatorname{deg}_{\rm CR}.
\een

On the other hand, if $X$ has a LG-mirror model, then we can define the calibrated periods $\widetilde{I}_\alpha^{(-\ell)}(\la)$ in the same way as in formulas \eqref{per} and \eqref{eqn:calib_lim}. The main motivation for the above construction is the following conjecture, which is motivated by Iritani's mirror symmetry theorem in \cite{Ir}. To simplify the formulation we set all Novikov variables to be $1$. Using the divisor equation one can recover easily the Novikov variables. 
\begin{conjecture}[$\Gamma$-Conjecture for the Milnor lattice]\label{Gamma-conj}
Given an integral cycle $\alpha$ there exists a class $V_\alpha\in K(X)$ in the $K$-theory of vector bundles, s.t. for all $\ell\gg 0$,
\ben
\frac{1}{\sqrt{2\pi}}\int_0^\infty e^{-\la s}\widetilde{I}^{(-\ell)}_{\alpha}(\la) d\la = 
s^{-\theta_X -\ell-1/2}s^{-\rho_X} \Psi(V_\alpha).
\een
\end{conjecture}
The conjecture can be refined even further, by saying that if $\alpha$ is a vanishing cycle then 
$V_\alpha$ can be represented by an exceptional object in the derived category $\mathcal{D}^b(X)$ and that the monodromy transformations of $\alpha$ correspond to certain mutation operations in $\mathcal{D}^b(X)$. See \cite{GGI} for more discussions. 
Now we describe Conjecture \ref{Gamma-conj} in the case of $\mathbb{P}^1_{\bf a}$.

\subsubsection{The $K$-ring of $\mathbb{P}^1_{\bf a}$}\label{p1a:K-ring}
Let ${\bf a}=(a_1,a_2,a_3)$ be a triple of non-negative integers and put $X=\mathbb{P}^1_{\bf a}$. The orbifold $\mathbb{P}^1_{\bf a}$ can be constructed as follows. Put 
\ben
G=\{t=(t_1,t_2,t_3)\in (\C^*)^3\ |\ t_1^{a_1}=t_2^{a_2}=t_3^{a_3} \}.
\een 
We have 
\ben
\mathbb{P}^1_{\bf a} = [Y_{\bf a}/G],\quad Y_{\bf a}=\{y=(y_1,y_2,y_3)\in \C^3\setminus{\{0\}}\ |\ y_1^{a_1}+y_2^{a_2}+y_3^{a_3} = 0\},
\een 
where the quotient is taken in the category of orbifolds, i.e., it should be viewed as an orbifold groupoid. The $K$-ring of orbifold vector bundles on $\mathbb{P}^1_{\bf a}$ can be presented as a quotient of the polynomial ring $\C[L_1,L_2,L_3]$ by the following relations
\ben
L=L_1^{a_1}=L_2^{a_2}=L_3^{a_3},\quad
(1-L_k)(1-L_{k'}) = 0 \ (1\leq k<k'\leq 3). 
\een
Here $L$ is the pullback of $\mathcal{O}_{\mathbb{P}^1}(1)$ under the natural map $\mathbb{P}^1_{\bf a}\to \mathbb{P}^1$, and the product is given by tensor product of vector bundles. The orbifold vector bundle $L_k$ is the trivial line bundle $Y_{\bf a}\times \C$ equipped with the following $G$-action 
\ben
G\times L_k\to L_k,\quad (t,y,v)\mapsto (t y, t_k v).
\een
It is easy to see that the $K$-ring is generated by $L_1, L_2, L_3, L$. The first set of relations follows from the definition of $G$. 
To see the remaining ones, note that the coordinate function $y_k$ on $Y_{\bf a}$ gives rise to a section of $L_k.$ The Koszul complex associated with the sections $(y_k,y_{k'})$ is $G$-equivariaint and it gives rise to the exact sequence 
\ben
0\to L_k^\vee\otimes L_{k'}^\vee \to  L_k^\vee\oplus L_{k'}^\vee\to \mathcal{O}_{\mathbb{P}^1_{\bf a}} \to 0.
\een 
This proves that $(1-L_k)(1-L_{k'})=0$.

\subsubsection{The image of $K(X)$}
The connected components of $\mathbb{P}^1_{\bf a}$ are indexed by $\{(0,0)\}\cup \mathfrak{I}_{\rm tw}$.
Let us denote by $P=c_1(L)$, then $c_1(L_k)=P/a_k$. By the adjunction formula $TX=L_1L_2L_3L^{-1}$, we get 
\ben
c_1(TX)=\chi\, P,\quad \chi=\frac{1}{a_1}+\frac{1}{a_2}+\frac{1}{a_3}-1.
\een
Furthermore, note that 
\ben
(L_k)_{k',p,f} = 
\begin{cases}
0 & \mbox{ if } k\neq k' \mbox{ and } f\neq 0 \\
0&   \mbox{ if } k= k' \mbox{ and } f\neq p/a_{k'} \\
\C & \mbox{ otherwise}.
\end{cases}
\een
From here we get that the eigenspace decomposition of $TX$ is
\ben
(TX)_{k,p,f} = 
\begin{cases}
TX  & \mbox{ if } k=0,\quad f=0, \\
\C & \mbox{ if } k\neq 0,\quad f=p/a_k, \\
0 & \mbox{ otherwise }.
\end{cases}
\een
Recall that for $i=(k,p)\in\mathfrak{I}_{\rm tw}$, $d_i=d_{k,p}=1-p/a_k$, we get the following formulas
\ben
\widehat{\Gamma}(X)&=&\Gamma(1+\chi P) + \sum_{i\in\mathfrak{I}_{\rm tw}} \Gamma(d_i) \phi_{i},\\
\widetilde{\rm ch}(L_k^m) &=& \one+\frac{m}{a_k}\,P + \sum_{(j,p)\in\mathfrak{I}_{\rm tw}}
\zeta_j^{mp\delta_{k,j}}\phi_{j,p}, \quad \zeta_j:=e^{2\pi\sqrt{-1}/a_j}.
\een
Let us point out that in the above formulas $\one,P\in H^*(X_{0,0})$, while $\phi_{k,p}\in H^0(X_{k,p})$ is the standard generator for the twisted sector. Note that the unit of the algebra ($H^*(IX;\C), \cup$) is
\ben
\widetilde{\rm ch}(\mathcal{O})={\bf 1} + \sum_{i\in\mathfrak{I}_{\rm tw}} \phi_{i}.
\een
Finally, since
\ben
(2\pi\sqrt{-1})^{\rm deg} {\rm inv}^* \widetilde{\rm ch}(L_k^m) = 
\one+\frac{2\pi\sqrt{-1}m}{a_k}\, P + \sum_{(j,p)\in\mathfrak{I}_{\rm tw}}
\zeta_{j}^{-mp\delta_{k,j}}\phi_{j,p},
\een
we get the following formula 
\beq\label{gamma-integral}
(2\pi)^{1/2}\Psi(L_k^m)=
\one+\left(-\gamma\chi+\frac{2\pi\sqrt{-1}m}{a_k}\right)P+ \sum_{(j,p)\in\mathfrak{I}_{\rm tw}}\frac{\Gamma(d_{j,p})}{\zeta_j^{mp\delta_{k,j}}}\phi_{j,p}.
\eeq

\subsection{$\Gamma$-conjecture for Fano orbifold curves}\label{sec:gamma}
Now we give a proof\footnote{Note that $\mathbb{P}^1_{\bf a}$ is not
  covered by results in \cite{Ir, Iri}.} of the $\Gamma$-conjecture
for the Milnor lattice for $\mathbb{P}^1_{\bf a}$. The proof is
obtained by applying Iritani's argument of the proof of \cite[Theorem
4.11]{Iri} and \cite[Theorem 5.7]{Iri} and relies on the
$\Gamma$-conjecture for the Milnor lattice for the Fano toric orbifold (proven in \cite{Ir})
\ben
Y:=\mathbb{P}^2_{\bf a}=[(\C^3\setminus{\{0\}})/G]
\een
and the explicit formulas for the $J$-functions of $X:=\mathbb{P}^1_{\bf a}$ and $Y$. 
Note that $X$ is a suborbifold of $Y$. 
\begin{remark}
There is a natural map $p: \mathbb{P}^2_{\bf a}\to \mathbb{P}^2$. The above description of $X=\mathbb{P}^1_{\bf a}$ realizes $X$ as the locus of zero of a section of the line bundle $p^*\mathcal{O}_{\mathbb{P}^2}(1)$ on $\mathbb{P}^2_{\bf a}$. Applying the recipe of constructing mirrors of complete intersections in \cite{G0}, we obtain $f_{\bf a}$ as the mirror of $X$. 
\end{remark}

Notice that the line bundles
$L_k$ are restrictions of line bundles on $Y$ and the
$K$-ring of $Y$ is the quotient of  the polynomial ring $\C[L_1,L_2,L_3]$ by the following relations
\ben
L=L_1^{a_1}=L_2^{a_2}=L_3^{a_3},\quad (1-L_1)(1-L_2)(1-L_3)=0.
\een
Put $L=p^*\mathcal{O}_{\mathbb{P}^2}(1)$ and $P=c_1(L)$. We have isomorphisms
\ben
\mathbb{Q} \cong  H_2(X;\mathbb{Q})\cong H_2(Y;\mathbb{Q}),\quad
d\mapsto d[\mathbb{P}^1_a]
\een
and 
\ben
H^2(Y;\mathbb{Q})\cong  H^2(X;\mathbb{Q})\cong \mathbb{Q},\quad
\alpha \mapsto \langle \alpha, [\mathbb{P}^1_a]\rangle .
\een

The $J$-function of an orbifold $X$ used by Iritani is
$$J_X(\tau,z)=L(\tau, z)^{-1}\one,$$ 
where $\tau\in H_{\rm CR}(X)$, 
$$L(\tau,z) := S_\tau(-z)  e^{-P\log Q/z}$$ 
and $S$ is the calibration operator \eqref{S:def}. 
Note that this definition differs from Givental's one by a sign and by the exponential factor. 

\subsubsection{Combinatorics of the inertia orbifolds} 
The orbifold $Y$ is toric. We describe its stacky fan as follows. Put 
\ben
b_1=(a_1,0),\quad b_2=(0,a_2),\quad b_3=(-a_3,-a_3) \in   \Z^2.
\een
The fan of $Y$ is 
\ben
\Sigma\cong \{\emptyset,\{k\}, \{k,k'\}\ |\ 1\leq k,k'\leq 3\}
\een
where each set $I$ on the RHS determines a cone in $\mathbb{R}^2$ spanned
by $b_k$, $k\in I$. Note that $\Sigma$ is the fan for $\mathbb{P}^2$. The fan map for $Y$ sends the standard basis $\{e_1, e_2, e_3\}$ of $\Z^3$ to $\Z^2$ by $$\mathbb{Z}^3\to \Z^2, \quad e_k\mapsto b_k.$$

The connected components of $IY$ are parametrized
by 
\ben
{\rm Box}(\Si) = \{(c_1,c_2,c_3)\ |\ 0\leq c_k< 1,\quad \sum_k c_k b_k
\in \Z^2\cap \si \mbox{ for some cones } \si\in \Si\},
\een
where $c\in {\rm Box}(\Si) $ determines the twisted sector
\ben
Y_c=[\{ y\in \C^3\ |\ y_k=0 \mbox{ if } c_k\neq 0  \}/G]
\een
which has a generic stabilizer given by the cyclic subgroup of $G$
generated by 
\ben
(e^{2\pi\sqrt{-1}c_1},e^{2\pi\sqrt{-1}c_2},e^{2\pi\sqrt{-1}c_3})\in G. 
\een
The inertia orbifold $IX$ is a suborbifold of $IY$ and the twisted
sectors of $IX$ are parametrized by those $c\in {\rm Box}(\Si)$ for
which $\operatorname{dim}(Y_c)>0$, i.e., at most one component of $c$
is non-zero. 
\subsubsection{The $J$-function of $Y$}\label{JX}
Let $\one_c\in H^0(Y_c)$ be the dual of the fundamental class for $c\in {\rm Box}(\Si)$ and
\ben
\tau=\tau_1 \one_{(1/a_1,0,0)}+\tau_2\one_{(0,1/a_2,0)}+\tau_3\one_{(0,0,1/a_3)}.
\een
According to the mirror theorem of \cite{CCIT2}, the $J$-function $J_Y(\tau, z)$ depending on $\tau$ is equal to the $S$-extended $I$-function \cite[Definition 28]{CCIT2} with $S=\{(1,0), (0,1), (-1,-1)\}$. This gives 
\ben
J_Y(\tau,z) = e^{P\log Q/z}\Big(\sum_{d=0}^\infty\
\sum_{n_1,n_2,n_3=0}^\infty \ 
\frac{Q^d}{z^{\operatorname{deg}_Y(Q^d) }} \, \frac{t^n}{n!z^{\operatorname{deg}_Y(t^n)}}\,  J^Y_{d,n}(\tau,z)\Big),
\een
where we introduced homogeneous parameters $t=(t_1,t_2,t_3)$, whose
dependance on $\tau$ and $Q$ can be determined from the expansion
$J_Y=1+\tau/z+\cdots$, the degrees of $Q$ and $t$ are
\ben
\operatorname{deg}_Y(Q^d):=\int_d c_1(TY)=d\Big(\frac{1}{a_1}+\frac{1}{a_2}+\frac{1}{a_3}\Big),\quad 
\operatorname{deg}_Y(t_k):= \operatorname{deg}_Y(\tau_k)= 1-1/a_k.
\een 
Finally, we denoted $n=(n_1,n_2,n_3)$ and we used the standard multi-index notations
\ben
t^n=t_1^{n_1}t_2^{n_2}t_3^{n_3},\quad n!=n_1!n_2!n_3!.
\een
In order to define
the component $J^Y_{d,n}$ let us define $m_k \in \Z$ and $ c_k\in \mathbb{Q}$ by 
\ben
\frac{n_k-d}{a_k}= -m_k + c_k, \quad 0\leq c_k<1.
\een
Then we have
\ben
J^Y_{d,n}(\tau,z) = 
\frac{\one_c}{z^{\operatorname{deg}_Y(\one_c)}}
\prod_{k=1}^3 \frac{\Gamma(1-c_k+(P/a_k)z^{-1}) }{\Gamma(1-c_k+m_k+(P/a_k)z^{-1}) },
\een
where if $c\notin {\rm Box}(\Sigma)$ then we set $\one_c=0$. In other words we sum over all $(d,n)$, s.t., at least one of the numbers $c_k$ is $0$.

\subsubsection{The $J$-function of $X$}\label{sec:JX}
Since $p^*\mathcal{O}_{\mathbb{P}^2}(1)$ is a convex line bundle in the sense of \cite[Example B]{CCIT1}, the $J$-function of $\mathbb{P}^1_{\bf a}$ can be computed from that of $\mathbb{P}^2_{\bf a}$ using the quantum Lefschetz theorem of \cite{tseng} and \cite{CCIT1}. 

Using the embedding $j:IX\to IY$ we restrict $\tau$ and $\one_c$ to $H^*(IX)$. Slightly abusing the notation we use the same notation for the restrictions. Note that now $\one_c=0$ if $c$ has more than one non-zero component. 
The formula for $J_X$ has the same form
\ben
J_X(\tau,z) = e^{P\log Q/z}\Big(\sum_{d=0}^\infty\
\sum_{n_1,n_2,n_3=0}^\infty \ 
\frac{Q^d}{z^{\operatorname{deg}_X(Q^d) }} \, \frac{t^n}{n!z^{\operatorname{deg}_X(t^n)}}\,  J^X_{d,n}(\tau,z)\Big),
\een
where 
\ben
J^X_{d,n}(\tau,z) = 
\frac{\one_c}{z^{\operatorname{deg}_X(\one_c)}}
\frac{\Gamma(1+d+Pz^{-1}) }{\Gamma(1+Pz^{-1}) }
\prod_{k=1}^3 \frac{\Gamma(1-c_k+(P/a_k)z^{-1}) }{\Gamma(1-c_k+m_k+(P/a_k)z^{-1}) }
\een
Note that the grading takes the form
\ben
\operatorname{deg}_X(Q^d):=\int_d c_1(TX)=d\Big(\frac{1}{a_1}+\frac{1}{a_2}+\frac{1}{a_3}-1\Big)
\een
while the degrees of $t$ and $\one_c$ do not change, because the restriction map preserves the grading. 

\subsubsection{The Galois action}
The Picard group $\operatorname{Pic}(X)$  of isomorphism classes of (topological) orbifold line bundles on $X$ can be presented as a quotient
\ben
\Z^3\to \operatorname{Pic}(X),\quad (r_1,r_2,r_3)\mapsto L_1^{r_1}L_2^{r_2}L_3^{r_3}
\een
with kernel given by the relations 
\ben
a_1e_1=a_2e_2=a_3e_3,
\een
where $\{e_1,e_2,e_3\}$ is the standard basis of $\Z^3$.  
The group $\operatorname{Pic}(X)$ acts naturally on the Milnor fibration via
\ben
\nu\cdot (x,t) = (\nu\cdot x, \nu\cdot t), \quad \nu=(r_1,r_2,r_3)\in \operatorname{Pic}(X),
\een
where
\ben
(\nu \cdot x)_k  & = &  e^{2\pi\sqrt{-1} r_k/a_k }x_k, 
\een
and the action on the remaining components is defined in such a way that 
\ben
F(\nu\cdot x,\nu\cdot t)=F(x,t),
\een 
i.e., 
\ben
(\nu\cdot t)_{k,p} & = &  e^{-2\pi\sqrt{-1}  r_k p/a_k}\,t_{k,p},\quad 1\leq k\leq 3,\quad 1\leq p\leq a_k-1,\\
(\nu\cdot t)_{02} & =& t_{02} +2\pi\sqrt{-1} \sum_{k=1}^3 \frac{r_k}{a_k}, \\
(\nu\cdot t)_{01} & =& t_{01} .
\een
Let us fix some $(t,\la)\in M\times \C$ with $\la$ sufficiently large, then for every $\nu=(r_1,r_2,r_3)$ we can construct a path from $(t,\la)$ to $(\nu\cdot t,\la)$ as follows. Using the above formulas we let $c\in \mathbb{R}^3$ act on $M$. As $c$ varies along the straight segment from $0$ to $(r_1,r_2,r_3)\in \Z^3\subset \mathbb{R}^3$ we get a path in $M$ connecting $t$ and $\nu\cdot t$. The parallel transport along this path with respect to the Gauss-Manin connection gives an identification $H_2(X_{\nu\cdot t,\la};\Z) \cong H_2(X_{t,\la};\Z).$ Combined with the $\operatorname{Pic}(X)$-action on $\C^3$ we get an action\ben 
\operatorname{Pic}(X)\times H_2(X_{t,\la};\Z)\to H_2(X_{t,\la};\Z),\quad (\nu,\alpha)\mapsto \nu(\alpha).
\een  
Following Iritani, we refer to the above action as Galois action of $\operatorname{Pic}(X)$ on the Milnor lattice.
\begin{lemma}\label{galois}
If the $\Gamma$-conjecture for the Milnor lattice is true for some cycle $\alpha$ and $V_\alpha\in K(X)$ is the corresponding $K$-theoretic vector bundle, then the conjecture is true for all $\nu(\alpha)$, $\nu=(r_1,r_2,r_3)\in \operatorname{Pic}(X)$. Moreover,
\ben
V_{\nu(\alpha)}=V_{\alpha}\otimes L_\nu,\quad L_\nu=L_1^{r_1}L_2^{r_2}L_3^{r_3}.
\een
\end{lemma}
\proof
Using the vector space decomposition 
\ben
H_{\rm CR}(X) = H^*(X) \bigoplus\left(\bigoplus_{(k,p)\in\mathfrak{I}_{\rm tw}} H^{p/a_k}_{\rm CR}(X)\right),
\een
we define a linear operator
\ben
\theta_\nu :H_{\rm CR}(X) \to H_{\rm CR}(X),\quad \theta_\nu = 
\sum_{k=1}^3\sum_{p=1}^{a_k-1} \frac{r_k p}{a_k} {\rm pr}_{k,p},
\een
where ${\rm pr}_{k,p}$ is the projection onto the subspace $ H^{p/a_k}_{\rm CR}(X)$. 
By changing the variables $y=\nu\cdot x$ in the period integrals we get 
\ben
I^{(\ell)}_{\nu(\alpha)}(t,\la) = e^{-2\pi\sqrt{-1}\theta_\nu} I^{(\ell)}_\al(\nu^{-1}\cdot t,\la),\quad \forall \ell\in \Z. 
\een
On the other hand the calibration operator satisfies 
\ben
S_{\nu^{-1}(t)}(z)=e^{2\pi\sqrt{-1}\theta_\nu} S_t(z) e^{-2\pi\sqrt{-1}\theta_\nu} e^{-2\pi\sqrt{-1} c_1(L_\nu)/z},
\een
which can be seen easily by using that  if the correlator 
\ben
\langle \alpha_1\psi_1^{k_1},\dots,\alpha_n\psi_n^{k_n}\rangle_{0,n,d}
\een
is non-zero then, since we have at least one stable map $f: C\to X$, we have
\ben
\chi(f^*L_\nu)=\int_d c_1(L_\nu) -\sum_{j=1}^n \theta_\nu(\alpha_j)\in \Z.
\een
Since by definition 
\ben
I^{(\ell)}_{\al}(t,\la) = S_{ t}(-\d_\la^{-1}) \widetilde{I}_{\al}^{\ell}(\la), 
\een
the above formulas imply that 
\ben
\widetilde{I}^{(\ell)}_{\nu(\al)}(\la) = e^{-2\pi\sqrt{-1}\theta_\nu} e^{2\pi\sqrt{-1} c_1(L_\nu)\d_\la}
\widetilde{I}_{\alpha}^{\ell}(\la) .
\een
In particular, after taking a Laplace transform, we get
\ben
\frac{1}{\sqrt{2\pi}} \int_0^\infty e^{-\la s} \widetilde{I}^{(-\ell)}_{\nu(\al)}(\la) =  
e^{-2\pi\sqrt{-1}\theta_\nu} e^{2\pi\sqrt{-1} c_1(L_\nu)s}
\frac{1}{\sqrt{2\pi}} \int_0^\infty e^{-\la s} \widetilde{I}^{(-\ell)}_{\al }(\la) .
\een
On the other hand, using the definition of $\Psi$ we get
\ben
\Psi(V\otimes L_\nu) = e^{-2\pi\sqrt{-1}\theta_\nu} e^{2\pi\sqrt{-1} c_1(L_\nu)} \Psi(V).
\een
It remains only to notice that $s^{-\theta} P = (P s) s^{-\theta}$ and that $\theta_\nu$ commutes with both $\theta$ and the Chen-Ruan product multiplication operators.
\qed

The Milnor lattice is known to be unimodular with respect to the $K$-theoretic bilinear form 
\ben
(\ ,\ ): K(X)\otimes_\Z K(X) \to \Z,\quad (L_1,L_2) = \chi(L_1\otimes L_2^\vee)
\een
 (see \cite{Iri}, Section 2). The above Lemma implies that it is enough to prove that the $\Gamma$-conjecture holds for the structure sheaf. Indeed, if this is true, then since $K(X)$ is generated by $\operatorname{Pic}(X)$, the $\Gamma$ conjecture correspondence will embed $K(X)$ into a sublattice of the Milnor lattice. Since both lattices are unimodular,  they must coincide.  

\subsubsection{The central charge}
Iritani's $\Gamma$-conjecture for the Milnor lattice looks different
since he works with Lefschetz thimbles. Nevertheless, our formulation
is completely equivalent. The reason is the following. Let us take a
Lefschetz thimble $\mathcal{A}$ corresponding to a vanishing cycle
$\alpha$, i.e., for fixed $(t,z)\in M\times \C^*$ we fix a path $C$ in
$\C$ from $u_j$ to $\infty$, s.t., $\operatorname{Re}(\la/z)>0$ for
all $\la\in C$ and the cycle $\alpha_{t,\la}$ vanishes when $\la$
approaches $u_j$. In this way we can identify the Milnor lattice with a lattice of Lefschetz thimbles. 

We claim that 
\ben
L(t,z)^{-1} \, \int_{u_j}^\infty e^{-\la/z} I^{(-\ell)}_\alpha(t,\la) d\la = 
e^{z^{-1}P\log Q} \int_0^\infty e^{-\la/z} \widetilde{I}^{(-\ell)}_\alpha(\la)d\la,
\een
where $L(t,z)=S_t(-z) e^{-z^{-1}P\log Q}$. Indeed, one can check easily using the quantum differential equations that the LHS is independent of $t$ and $Q$. On the other hand we have 
\ben
L(t,z)=1-z^{-1}P\log Q +\cdots,
\quad 
u_j=0 +\cdots 
\quad
\mbox{and} 
\quad
I^{(-\ell)}_\al(t,\la) = \widetilde{I}^{(-\ell)}_\alpha(\la) + \cdots,
\een
where the dots stand for terms that vanish at $t=Q=0$. So modulo terms that vanish at $t=Q=0$ the LHS coincides with the RHS. Our claim follows.

We define the {\em central charge} of $V_\alpha\in K(X)$ by
\ben
Z^{(0)}_X(V_\alpha)(t,z):=\left( L(t,z) z^{\theta } z^\rho \Psi(V_\alpha),\one \right).
\een
Since we will use the result of Iritani, let us clarify the relation
between our notations. In Iritani's notation, the central charge is defined to be 
\ben
Z^{(n)}_X(V)(t,z):= (2\pi z)^{n/2} (2\pi \sqrt{-1})^{-n} \Big( L(t,z) z^{\theta } z^\rho \Psi(V_\alpha),\one \Big), \quad n=\operatorname{dim}_\C(X).
\een 
For the LG models studied in
\cite{Ir} the $\Gamma$-conjecture for the central charge is stated as
\ben
(2\pi \sqrt{-1})^{-n} \int_\mathcal{A} e^{-F(x,t)/z}\omega = Z_X^{(n)}(V_\al)(t,z).
\een 
As we see from the LG model that we use in general one should choose
$n$ to be the number of variables in the LG potentials. For the LG models in \cite{Ir}  the number of variables coincides with the dimension of the orbifold, so this difference does not matter.

The identity in the $\Gamma$-conjecture for the Milnor lattice is equivalent to 
\beq\label{gc-param}
\frac{1}{\sqrt{2\pi}}\int_{u_j}^\infty e^{-\la/z} I^{(-\ell)}_\alpha(t,\la) d\la  = 
L(t,z) z^{\theta +\ell+1/2} z^{\rho} \Psi(V_\alpha).
\eeq
The number $\ell$ must be chosen sufficiently large. We will see that in our case $\ell=1$ works.
In general $\ell$ can be chosen, s.t., the number of variables in the LG potential is $2\ell+1$. 
Recalling the definition of the period integrals, we transform the LHS into 
\ben
(-z d_M) (2\pi)^{-3/2} \int_\mathcal{A} e^{-F(x,t)/z} \omega,
\een 
where $d_M$ is the de Rham differential on $M$. In particular, since the Poincar\'e pairing of the RHS with $\one$ corresponds to contracting the LHS with $\d_{01}$ we get
\beq\label{cc}
(2\pi z)^{-3/2} \int_\mathcal{A} e^{-F(x,t)/z} \omega = Z^{(0)}_X(V_\alpha)(t,z).
\eeq

In order to prove the $\Gamma$-conjecture for the Milnor lattice it is enough to prove that if $V=\mathcal{O}_X$, then we can find an integral cycle $\mathcal{A}$, s.t., the identity \eqref{cc} holds for all parameters $t$ of the form
\ben
t=t_{1,1} \one_{1,1}+t_{2,1}\one_{2,1}+t_{3,1}\one_{3,1}.
\een
One can check that the partial derivatives of the LHS and the RHS of \eqref{cc} with respect to any other parameter $t_{k,p}$ can be expressed in terms the same differential operator involving only $t_{k,1}$, $1\leq k\leq 3$. Therefore if \eqref{cc} holds for all $t$ of the above form, then \eqref{gc-param} holds also for all such $t$. As it was explained above if the identity \eqref{gc-param} holds for a single point $t=t_0$ then it holds for all $t$ and it is equivalent to the identity in our $\Gamma$-conjecture. In other words, the $\Gamma$-conjecture holds for the structure sheaf. Recalling Lemma \ref{galois} we get that the $\Gamma$-conjecture holds for the entire Milnor lattice.  

\subsubsection{The central charge as an oscillatory integral}
It remains only to prove \eqref{cc}. Following Iritani, it is
convenient to rewrite the RHS of \eqref{cc} in terms of the so-called
$H$-function 
$$H^{(0)}_X(t,z)=\widetilde{\operatorname{ch}}(H^{(0)}_K(t,z)),$$ 
where the $K(X)$-valued function
\ben
H^{(0)}_K: M\times \C^*\to K(X)
\een
is defined by the equation
\ben
\one=L(t,z)z^\theta z^\rho \Psi(H^{(0)}_K(t,z)).
\een
For the central charge $Z^{(0)}_X(V) $ we have 
\ben
(L(t,z)z^\theta z^\rho\Psi(V),L(t,-z) (-z)^\theta (-z)^\rho \Psi(H^{(0)}_K(t,-z)))=
(\Psi(V), e^{\pi\sqrt{-1} \theta} e^{\pi\sqrt{-1} \rho} \Psi(H^{(0)}_K)),
\een
where we define $(-1)^R := (e^{\pi\sqrt{-1}R})$ for all linear operators $R$. 
Recalling \eqref{K-pairing} and the Kawasaki Riemann--Roch formula we get
\ben
Z^{(0)}_X(V) = \chi(H^{(0)}_K\otimes V^\vee) = 
\int_{IX} H^{(0)}_X(t,-z)\cup \widetilde{\operatorname{ch}}(V^\vee)\cup \widetilde{\operatorname{Td}}(TX),
\een
where in the notation of Section \ref{ir-lattice} the Todd class of an orbifold vector bundle is a multiplicative characteristic class defined by 
\ben
\widetilde{\operatorname{Td}}(V) = 
\sum_{v\in T} \ 
\prod_{j=1}^{l_{v,0}} \frac{\delta_{v,0,j}}{1-e^{-\delta_{v,0,j}}} \,
\prod_{0<f<1} \prod_{j=1}^{l_{v,f}}\frac{1}{1-e^{-2\pi\sqrt{-1} f}e^{-\delta_{v,f,j}}}.
\een
The proof of formula \eqref{cc} requires a simple lemma. The main ingredient is a slight modification of the usual Laplace transform defined as follows. Let $f(t,Q;z)$ be any function, then we define 
\ben
\mathfrak{L}\left(f\right)(t,Q;z) = \int_0^\infty e^{-\eta} f(t,-\eta z Q; z) d\eta.
\een
The integral is convergent if for example $f$ depends polynomially on $Q$ and $\log Q$, which is the case that we have. Note that this Laplace transform does not commute with the involution $z\mapsto -z$. 
\begin{lemma} \label{H:LT}
Let $j: IX\to IY$ be the natural embedding, then 
\ben
j_*H^{(0)}_X(t,Q;-z)= (-z/2\pi)^{1/2}\, \widetilde{e}(L) \cup \mathfrak{L} (H^{(0)}_Y) (t,Q;-z), 
\een 
where $\widetilde{e}(L)=\sum_{v\in T} e(L_v)$ is the orbifold Euler class of $L$. 
\end{lemma}
\proof
Since $j_*(j^*\alpha)=\widetilde{e}(L)\cup \alpha$ for every $\al\in H^*(IY)$, it is enough to prove that 
\beq\label{jHY}
\mathfrak{L}\left(j^* H^{(0)}_Y\right) (t,Q;-z) =  (-z/2\pi)^{-1/2}\  H^{(0)}_X(t,Q;-z).
\eeq
We have
\ben
(2\pi)^{-1/2}\widehat{\Gamma}(X)\cup (2\pi\sqrt{-1})^{\rm deg}{\rm inv}^* H^{(0)}_X(t,Q,z) = 
z^{-\rho_X}z^{-\theta_X} J_X(t,Q;z)
\een
and 
\ben
(2\pi)^{-1}\widehat{\Gamma}(Y)\cup (2\pi\sqrt{-1})^{\rm deg}{\rm inv}^* H^{(0)}_Y(t,Q,z) = 
z^{-\rho_Y}z^{-\theta_Y} J_Y(t,Q;z)
\een
On the other hand, using the explicit formulas for the $J$-functions it is easy to check that 
\ben
\mathfrak{L}(j^*J_Y)(t,Q;-z) = (-z)^{-P/z}\Gamma(1-P/z)\cup J_X(t,Q;-z).
\een
In order to prove formula \eqref{jHY}, it is enough only to recall the following identities
\ben
j^* (-z)^{-\rho_Y}(-z)^{-\theta_Y} = (-z)^{-P-1/2} (-z)^{-\rho_X}(-z)^{-\theta_X}  j^*,
\een
\ben
 (-z)^{-\rho_X}(-z)^{-\theta_X} (-z)^{-P/z}\Gamma(1-P/z) = 
(-z)^P\Gamma(1+P) (-z)^{-\rho_X}(-z)^{-\theta_X}, 
\een
and
\ben
j^*\widehat{\Gamma}(Y) = \widehat{\Gamma}(L) \widehat{\Gamma}(X).
\een
\qed 

Lemma \ref{H:LT} yields the following relation between the central charges of sheaves on $X$ and $Y$. Let $V\in K(Y)$, then 
\ben
Z^{(0)}_X(j^*V) = (-z/2\pi)^{1/2}\mathfrak{L}\left(Z^{(0)}_Y(V-V\otimes L)\right).
\een 
In particular
\beq\label{cc:LT}
Z^{(0)}_X(1) =  (-z/2\pi)^{1/2}\mathfrak{L}\left(Z^{(0)}_Y(1- L)\right).
\eeq

\begin{theorem}\label{thm-gamma}
For a Fano orbifold curve $X=\mathbb{P}^1_{\bf a}$, 
given a class $L_k^m\in K(X)$ in the $K$-theory of vector bundles, there exists an integral cycle $\alpha_{k,m}\in\lieh$, s.t. for all $\ell\gg 0$,
\ben
\frac{1}{\sqrt{2\pi}}\int_0^\infty e^{-\la s}\widetilde{I}^{(-\ell)}_{\alpha_{k,m}}(\la) d\la = 
s^{-\theta_X-\ell-1/2}s^{-\rho_X} \Psi(L_k^m).
\een
\end{theorem}
\proof
It is enough to prove \eqref{cc}.
Let us look at the corresponding oscillatory integrals. Recall that the LG model of $Y$ is given by the restriction of
\ben
F_{\mathbb{P}^2}(x,t) = \sum_{k=1}^3 (x_k^{a_k} + t_{k,1} x_k),
\een
to the complex torus $x_1x_2x_3=Q$, while the corresponding primitive form is 
\ben
\omega_{\mathbb{P}^2} = \frac{dx_1dx_2dx_3}{d(x_1x_2x_3)}.
\een
Let us assume now that $z$ and $Q$ are real numbers, s.t., $z>0$ and $Q<0$. Let $\mathcal{C}\subset \C^3$ be the chain 
\ben
\mathcal{C}=\{x\in \mathbb{R}^3\ |\ x_k\geq 0, k=1,2,3\}.
\een
The oscillatory integral 
\ben
(2\pi z)^{-3/2} \int_\mathcal{C} e^{-F(x,t)/z} \omega = 
(2\pi z)^{-3/2} (-1)^{1/2}\int_0^\infty e^{-\eta} \int_{\Gamma_{-zQ\eta}}e^{-F_{\mathbb{P}^2}(x,t)/z} 
\frac{dx_1dx_2dx_3}{d(x_1x_2x_3)}  (-z)d\eta,
\een
where we presented the chain $\mathcal{C}$ as a family of cycles 
\ben
\Gamma_{-z\eta Q}=\{ x\in \mathcal{C}\ |\ x_1x_2x_3 = -z\eta Q\}.
\een 
and used the Fubini theorem. The $\Gamma$-conjecture for $Y$ was
proved by Iritani \cite{Ir}. Moreover, the real cycle $ \Gamma_{-z\eta Q}$ corresponds to the structure sheaf $\mathcal{O}_Y$, so the above integral coincides with 
\ben
(-1)^{3/2}z (2\pi z)^{-3/2}\int_0^\infty e^{-\eta}  (2\pi z) Z^{(0)}_Y(1) (t,-z\eta Q;z)d\eta = (-1)^{3/2} (z/2\pi )^{1/2} \mathfrak{L}\left(Z^{(0)}_Y(1)\right). 
\een
Recalling the argument in Lemma \ref{galois} it is easy to see that the analytic continuation around
$Q=0$ in clockwise direction of $\mathfrak{L}\left(Z^{(0)}_Y(1)\right)$ is $\mathfrak{L}\left(Z^{(0)}_Y(L)\right)$, therefore the cycle that we are looking
for is $\widetilde{\mathcal{C}}-\mathcal{C}$, where $\widetilde{\mathcal{C}}$ is the chain obtained
from $\mathcal{C}$ via the monodromy transformation around $Q=0$ in the
clockwise direction. More precisely, $\widetilde{\mathcal{C}}$ is the family of
cycles $\widetilde{\Gamma}_{-z\eta Q}$ obtained from $ \Gamma_{-z\eta
  Q}$  by the monodromy transformation around $Q=0$. It remains only
to notice that the boundaries of $\widetilde{\mathcal{C}}$ and $\mathcal{C}$ are the same. Together with \eqref{cc:LT}, this proves \eqref{cc}.  
\qed

\subsection{Affine root systems and vanishing cycles}\label{sec:pmaps} 

According to Theorem \ref{thm-gamma} (recall that we have to put $Q=1$) and formula \eqref{gamma-integral}, we have 
\ben
\int_0^\infty e^{-\la s}\widetilde{I}^{(-\ell)}_{\alpha_{k,m}}(\la)d\la & = & 
\frac{\one}{s^{\ell+1}} + \Big( \frac{2\pi\sqrt{-1}m}{a_k} -\gamma\chi-\chi\log s\Big) \frac{P}{s^\ell}
+\sum_{(j,p)\in\mathfrak{I}_{\rm tw}}
\frac{\Gamma(d_j)}{e^{2\pi\sqrt{-1}mp\delta_{k,j}/a_j}}\frac{\phi_{j,p}}{s^{\ell+d_{j,p}}}.
\een
where $d_{j,p}=1-p/a_j$, $\gamma$ is the Euler's gamma constant defined by
\ben
\gamma=\lim_{m\to \infty} H_m-\ln m,\quad H_m:=1+\frac{1}{2}+\cdots +\frac{1}{m}, 
\een
If $\ell\geq 1$, then we can recall the inverse Laplace transform and also the divisor equation \eqref{calp-div} 
to get 
\beqa\label{calp-im}
\begin{split}
\widetilde{I}^{(-\ell)}_{\alpha_{k,m}}(\la)   =&   \frac{\lambda^\ell}{\ell!} \one+\frac{\lambda^{\ell-1}}{(\ell-1)!} 
\left(
\frac{2\pi\sqrt{-1}m}{a_k} + \chi(\log\lambda - C_{\ell-1}) \right) P +
\\ 
&+\sum_{(j,p)\in\mathfrak{I}_{\rm tw}}\frac{\lambda^{d_{j,p}+\ell-1} e^{2\pi\sqrt{-1}m\delta_{k,j}d_{j,p}} }{\left(d_{j,p}+\ell-1\right) \cdots \left(d_{j,p}\right) } 
\phi_{j,p},
\end{split}
\eeqa
where if $\ell=1$ we set $C_0:=\frac{1}{\chi}\log Q$ and if $\ell>1$ then $C_\ell=C_{\ell-1}+\frac{1}{\ell}$.
\begin{proposition}\label{affine}
\hfill
\begin{enumerate}
\item
The set of vanishing cycles $\Delta\subset \lieh=H_2(X_{0,1};\C)$ is an affine root system of type
$X_N^{(1)}$, where $N= a_1+a_2+a_3-2$ and 
\ben
X=
\begin{cases}
A & \mbox{  if } a_1=1, \\
D & \mbox{ if }  a_1=a_2=2, \\
E & \mbox{ otherwise}.
\end{cases} 
\een
\item
There exists a basis of simple roots such that the
classical monodromy $\si$ is an affine Coxeter transformation. 
\end{enumerate}
\end{proposition}
Part (1) of Proposition \ref{affine} is due to A. Takahashi (see
\cite{Ta}). The proof is based on a standard method developed by
Gusein-Zade and A'Campo. We give a proof of Proposition \ref{affine} based on Iritani's integral structure.

We will be interested in the two maps from the sequence
\beq\label{per-mor}
\tI^{(n)}(1):\lieh \to H,\quad \al\mapsto \tI_\al^{(n)}(1)
\eeq
corresponding to $n=-1$ and $n=0$. According to Lemma \ref{per-1} we
have
\ben
\tI^{(-1)}(1) = B_{01}\, (\one-(\log Q) P) +B_{02}\, P + 
\sum_{i\in\mathfrak{I}_{\rm tw}}\, B_i\,\phi_i,
\een
which proves that the map for $n=-1$ is an isomorphism.
Using $\tI^{(-1)}(1)$ we equip $H$ with an
intersection pairing $(\cdot|\cdot)$, i.e., 
\ben
(\phi'|\phi''):=(\alpha'|\alpha''), \quad  \mbox{for}\quad 
\phi'=\tI^{(-1)}_{\al'}(1),\quad \phi''=\tI^{(-1)}_{\al''}(1).
\een
The period map \eqref{per-mor} with $n=0$ has a 1-dimensional
kernel because (using  \eqref{calp-2})
\ben
\tI^{(0)}(1)=(1-\rho)^{-1}(\theta+1/2) \tI^{(-1)}(1) = (1+\rho)
(1-\operatorname{deg}_{\rm CR}) \tI^{(-1)}(1) ,
\een
so the kernel is $\C\, P$. We denote the image of $\tI^{(0)}(1)$ by $H^{(0)}$. 
Let us denote by $r: H\to H^{(0)}$ the map defined by
$\tI^{(0)}(1)=r\circ \tI^{(-1)}(1)$, i.e.,
\ben
r(b) =  (1+\rho) (1-\operatorname{deg}_{\rm CR}) (b).
\een
According to Saito's formula \eqref{int-form} the intersection pairing on $H$ takes the form
\beq\label{int-form-QC}
(\phi'|\phi'') = (r(\phi'),(1-\rho)r(\phi'')),\quad \phi',\phi''\in H.
\eeq
It follows that we can pushforward the intersection form
to a non-degenerate bilinear pairing on $H^{(0)}$, which we denote
again by $(\cdot | \cdot)$. More precisely we define
\ben
(\phi'|\phi'') = (\phi',(1-\rho)\phi''),\quad \phi',\phi''\in H^{(0)}.
\een
Let us denote by  $\Delta^{(-1)}\subset H$ and $\Delta^{(0)}\subset H^{(0)}$ the images of the set of vanishing cycles, i.e., 
\ben
\Delta^{(-1)}=\{\tI^{(-1)}_\al(1)\ |\ \al\in \Delta\},\quad
\Delta^{(0)}=\{\tI^{(0)}_\al(1)\ |\ \al\in \Delta\}.
\een

A straightforward computation with formula \eqref{int-form-QC} implies
\begin{lemma}
Consider $\alpha_{k,m}$ as in \eqref{calp-im}.
Then the cycles $\alpha_{k,m}$ ($1\leq k\leq 3$, $m\in \Z$) satisfy
\ben
(\al_{k,m}|\al_{k,n})=
\begin{cases}
2 & \mbox{if } m= n ({\rm mod}\ a_k),\\
1 & \mbox{if } m\neq n ({\rm mod}\ a_k),\\
\end{cases}
\een
and for $k\neq k'$
\ben
(\al_{k,m}|\al_{k',n})=
\begin{cases}
2 & \mbox{if } m=0 ({\rm mod}\ a_k) \mbox{ and } n=0 ({\rm mod}\ a_{k'}) ,\\
0 & \mbox{if } m\neq 0 ({\rm mod}\ a_k) \mbox{ and } n\neq 0 ({\rm mod}\ a_{k'}) ,\\ 
1 & \mbox{otherwise} .
\end{cases}
\een
\end{lemma}
\subsubsection{The toroidal cycle}\label{sec:tor_cycle}
Let $\Gamma_\ep\subset \C^3$ be the torus $$\Gamma_\epsilon:=\{|x_1|=|x_2|=1,
|x_3|=\ep\}.$$ If $\ep$ is sufficiently large, $\Gamma_\ep$ does not
intersect the Milnor fiber $X_{0,1}$. Hence we have a well-defined cycle 
\ben
[\Gamma_\ep]\in H_3(\C^3\setminus X_{0,1};\Z) \cong
H_2(X_{0,1};\Z),
\een
where the isomorphism is given by the so called {\em tube mapping}
(for more details see \cite{Gri}). Let us denote by $\varphi$ the
image of $[\Gamma_\ep]$ under the above isomorphism. 
\begin{proposition}\label{per-tor}
We have $I^{(-1)}_\varphi(t,\la) = 2\pi\sqrt{-1}\, P$. 
\end{proposition}
\proof
Increasing $\ep$ does not change the homology class $[\Gamma_\ep]$,
therefore by choosing $\ep\gg 0$ we may arrange that $\Gamma_\ep$ does
not intersect the Milnor fiber $X_{t,\la}$ for all $(t,\la)$
sufficiently close to $(0,1)$. In particular, the cycle
$\varphi_{t,\la}$ obtained from $\varphi$ via a parallel transport
with respect to the Gauss--Manin connection coincides with the image of 
$[\Gamma_\ep]$ via the tube mapping. We have (c.f. \cite{Gri})
\beq\label{tor-int}
I(t,\la,Q):=\int_{[\Gamma_\ep]} \frac{\omega}{F(t,x)-\la} = 2\pi\sqrt{-1}
\int_{\varphi_{t,\la}} \frac{\omega}{dF} = 2\pi\sqrt{-1}\, \d_\la\,
\int d^{-1}\omega. 
\eeq
Comparing with the definition \eqref{per} we get 
\ben
I(t,\la,Q) = -(2\pi)^2\sqrt{-1}\, (I^{(-1)}_\varphi(t,\la),\one).
\een
Using the differential equation \eqref{2str-con-3}, we get 
\beq\label{grading}
(\la\d_\la + E) I(t,\la,Q) = 0.
\eeq
The integral $I(t,\la,Q)$ is analytic at $(t,\la,Q)=(0,0,0)$ because
it has the form
\ben
\sqrt{-1}\, \int_{[\Gamma_\ep]} 
\frac{dx_1dx_2dx_3}{Qe^{t_{02}} \, (G(t,x)-\la) - x_1 x_2 x_3}\, ,
\een
where $G(t,x)$ is a holomorphic function in $t$ and $x$. However,
equation \eqref{grading} means that $I(t,\la,Q)$ is homogeneous of
degree $0$ and since the weights of all variables are positive, the
integral must be a constant. In particular, we may set $t=Q=\la=0$,
which gives
\ben
I(t,\la,Q)= -\sqrt{-1}\, \int_{[\Gamma_\ep]} 
\frac{dx_1dx_2dx_3}{ x_1 x_2 x_3} = (2\pi)^3.
\een
Note that equation \eqref{tor-int} implies that
$I^{(0)}_\varphi(t,\la)=0$. Recalling again the differential equation
\eqref{2str-con-3}, we get 
\ben
I^{(-1)}_\varphi(t,\la) = (I^{(-1)}_\varphi(t,\la) ,\one)\, P =
(2\pi)^{-2} \sqrt{-1}\, I(t,\la,Q) \, P= 2\pi\sqrt{-1}\, P.
\een 
\qed

We immidiately have a corollary
\begin{corollary}
The cycle $\varphi$ corresponds to the skyscraper sheaf $\mathcal{O}_{\rm pt}:=L-\mathcal{O}$, i.e., 
\ben
\delta:=\tI^{(-1)}_{\alpha_{k,a_k}}(1)-\tI^{(-1)}_{\alpha_{k,0}}(1) = 2\pi\sqrt{-1} P.
\een
\end{corollary}

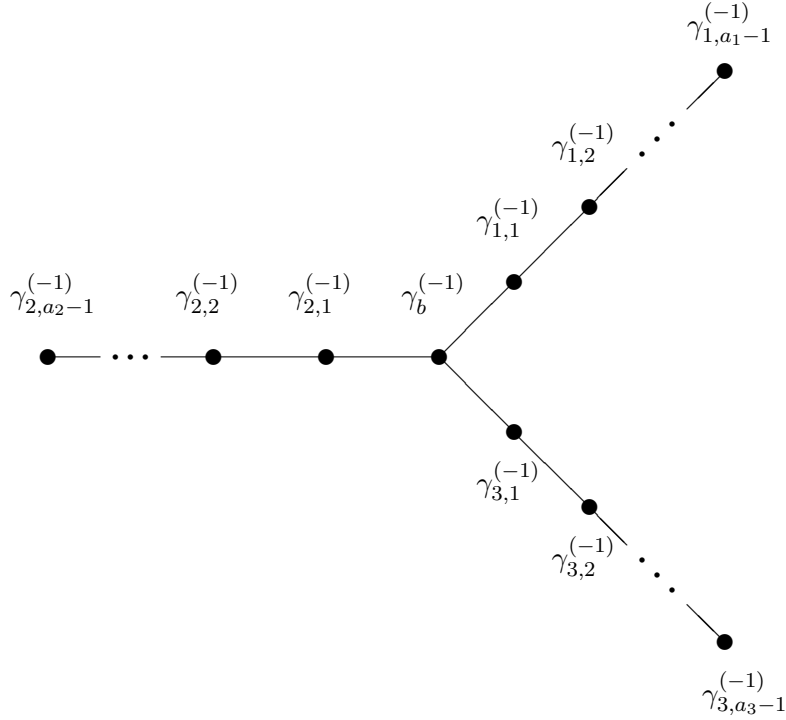
\begin{figure}[htbp]
\begin{center}
\begin{picture}(0,120)
\put(-5,57){$\gamma_b^{(-1)}$}
\put(0,50){\circle*{2}}
\put(0,50){\line(1,1){25}}
\multiput(27,77)(2,2){3}{\circle*{0.7}}
\put(33,83){\line(1,1){5}}
\put(5,67){$\gamma_{1,1}^{(-1)}$}
\put(10,60){\circle*{2}}
\put(15,77){$\gamma_{1,2}^{(-1)}$}
\put(20,70){\circle*{2}}
\put(33,93){$\gamma_{1,a_1-1}^{(-1)}$}
\put(38,88){\circle*{2}}
\put(0,50){\line(1,-1){25}}
\multiput(27,23)(2,-2){3}{\circle*{0.7}}
\put(33,17){\line(1,-1){5}}
\put(5,32){$\gamma_{3,1}^{(-1)}$}
\put(10,40){\circle*{2}}
\put(15,22){$\gamma_{3,2}^{(-1)}$}
\put(20,30){\circle*{2}}
\put(35,4){$\gamma^{(-1)}_{3,a_3-1}$}
\put(38,12){\circle*{2}}
\put(0,50){\line(-1,0){37}}
\multiput(-39,50)(-2,0){3}{\circle*{0.7}} 
\put(-45,50){\line(-1,0){7}}
\put(-20,57){$\gamma_{2,1}^{(-1)}$}
\put(-15,50){\circle*{2}}
\put(-35,57){$\gamma_{2,2}^{(-1)}$}
\put(-30,50){\circle*{2}}
\put(-57,57){$\gamma_{2,a_2-1}^{(-1)}$}
\put(-52,50){\circle*{2}}
\end{picture}
\caption{The branching node}
\label{fig:br}
\end{center}
\end{figure}

\proof[Proof of Proposition \ref{affine} (1)]
The image of the Milnor lattice in $H$ has a $\Z$-basis given by 
\ben
\delta,\quad \gamma^{(-1)}_b,\quad \gamma_{i}^{(-1)}, \quad i\in\mathfrak{I}_{\rm tw},
\een
where for $n=0$ or $-1$, and $i=(k,p)$, we get
\ben
\gamma^{(n)}_b & := & \tI^{(n)}_{\alpha_{k,0}}(1)\quad (\mbox{corresponds to } \mathcal{O}) \\
\gamma_{k,p}^{(n)} & := & \tI^{(n)}_{\al_{k,-p}} (1)-\tI^{(n)}_{\al_{k,-p+1}}(1).
\een
It is easy to check that the intersection diagram
of the set of cycles $\gamma_{b}^{(-1)},\gamma_{i}^{(-1)}, i\in\mathfrak{I}_{\rm tw}$ is given by the Dynkin diagram on Figure \ref{fig:br}. As usual, each node has self-intersection 2, each edge means that the intersection of the cycles corresponding to the nodes of the edge is $-1$, and no edge means that the intersection is $0$. It follows that the intersection form of the Milnor lattice is semi-positive definite with 1 dimensional kernel. This is possible only if $\Delta^{(-1)}$ is an affine root system. 
\qed

In particular we get also that $\delta$ is a $\Z$-basis for the imaginary roots and that 
$\Delta^{(0)}=r(\Delta^{(-1)})$ is a finite root system.

\subsubsection{Splitting of the affine root system}
It is convenient to enumerate the roots $\gamma_{b}^{(n)},\gamma_{i}^{(n)}, i\in\mathfrak{I}_{\rm tw}$ also by $\gamma_j^{(n)} (1\leq j\leq N)$.  The Dynkin diagram on Figure \ref{fig:br} is of type $X_N$, $X=ADE$. Let us denote by $\gamma_0^{(-1)}$ the {\em affine vertex}, i.e., the
extra node that we have to attach 
to $X_N$ in order to obtain the corresponding affine Dynkin diagram $X_N^{(1)}$. 

Vectors $\gamma_j^{(0)}$, $1\leq j\leq N$, form a basis
of simple roots of $\Delta^{(0)}.$ Let $W^{(0)}$ be the
reflection group generated by $\gamma_j^{(0)}$. It is well known that there exists a group embedding $W^{(0)}\to W$ which is induced by the map
\ben
s^{(0)}_j:=s_{\gamma_j^{(0)}} \mapsto s^{(-1)}_j:=s_{\gamma_j^{(-1)}},
\quad 1\leq j\leq N
\een
Given $\al\in \Delta^{(0)}$, let us define a {\em lift} $\tilde{\al}\in
\Delta^{(-1)}$ as follows
\ben
\al=\sum_{j=1}^N b_j \gamma_j^{(0)}\quad  \mapsto \quad 
\tilde{\al}:=\sum_{j=1}^N b_j \gamma_j^{(-1)}. 
\een
Then the root system $\Delta^{(-1)}$ coincides with the set
\ben
\Big\{\tilde{\al}+n\,\delta\ |\ \al\in \Delta^{(0)},\ n\in \Z\Big\}\ ,
\een
where $\delta=\gamma_0^{(-1)} +\theta^{(-1)}$ and $\theta\in
\Delta^{(0)}$ is the highest root with respect to the basis
$\{\gamma_j^{(0)}\}_{j=1}^N$ (see \cite{Kac}). Following Kac, we will
refer to $n\,\delta$ ($n\in \Z$) as {\em imaginary roots}.   
Finally, let us denote by $$\Lambda^{(-1)} :=H_2(X_{0,1};\Z)$$ the root lattice
of $\Delta^{(-1)}$. Given $\al\in \Lambda^{(-1)}$ such that $|\al|^2:=(\al | \al)\neq 0$, recall that the reflection with respect to $\al$ is defined by 
\ben
s_\al(x)=x-2\frac{(\al|x)}{(\al|\al)}\,\al.
\een
We also define the following translation:
\ben
T_\al(x) := s_{\al+\delta}s_\al(x)= x+2\frac{(\al|x)}{(\al|\al)}\,\delta.
\een
This definition induces a group embedding $T:
\Lambda^{(0)} \to W$. Recall that $w\,s_\al\, w^{-1} = s_{w(\al)}$
for all $w\in W$ and $\al\in \Lambda^{(-1)}$ such that $|\al|^2\neq 0$. 
Therefore, $\Lambda^{(0)}$ is a normal subgroup of $W$ and we
have an isomorphism $$W\cong \Lambda^{(0)} \rtimes W^{(0)}.$$ Let us
emphasize that the above isomorphism is not canonical -- it depends
on the choice of a basis of simple roots of $\Delta^{(-1)}$.  

\subsubsection{The Coxeter transformation}
Put $\si_b:=\si_b^{(0)}$, where
\beq\label{si-0}
\si_b^{(\ell)}=\prod_{k=1}^3 \Big(s^{(\ell)}_{k,a_k-1}\cdots
s^{(\ell)}_{k,2} s^{(\ell)}_{k,1}\Big)\quad\in\quad \operatorname{Aut}(\Delta^{(\ell)}),\quad \ell=-1,0.
\eeq
Note that while the order of the reflections that enter each factor of the above
product is important, the order in which the 3 factors 
are arranged is irrelevant since they pairwise commute. 
\begin{proposition}\label{Conjugacy}
The automorphism of $\Delta^{(0)}$ induced by the action of the
classical monodromy $\si$ coincides with $\si_b$.
\end{proposition}
\proof
The analytic continuation in $\la$ around $\la=0$ of the period $\tI^{(-1)}_{\al_{k,m}}(\la)$ is 
equivalent to tensoring the line bundle $L_k^m$ by $TX=L_1L_2L_3L^{-1}$ and then taking the corresponding periods. Using
\ben
(L_k-1)L_{k'} = L_k-1, \quad \mbox{for}\quad  k\neq k'
\een
it is easy to check that
\ben
(L_{k}^{-m}-L_{k}^{-m+1})TX & = & L_{k}^{-m+1}-L_{k}^{-m+2}\quad \forall m\in \Z, \\
TX^{-1}  & = &  1+ L_1^{-1}-1 + L_2^{-1}-1 +L_3^{-1}-1 +L-1, \\
(L_{k}^{-1}-1)TX & = &  1-L_{k}^{-(a_{k}-1)} + 1-L.
\een
According to the above remark, the classical monodromy acts as follows
\ben
\si(\gamma_{k,p}) & = & \gamma_{k,p-1},\quad (k,p)\in\mathfrak{I}_{\rm tw}, \\
\si^{-1}(\gamma_b) & = & \gamma_b +\gamma_{1,1}+\gamma_{2,1}+\gamma_{3,1} +\delta,\\
\si(\gamma_{k,1}) & = & -\gamma_{1,1}-\cdots -\gamma_{k,a_k-1}-\delta.
\een
It remains only to check that the action of $\si_b^{(-1)}$ is given by the same formulas modulo the imaginary root $\delta$.
\qed

It is known that up to a translation the affine Coxeter transformation coincides with $\si_b$ (see \cite{Ste, Stb}), so part (2) of Proposition \ref{affine} follows from Proposition \ref{Conjugacy}.

\subsection{Calibrated periods in terms of the finite root system}\label{cal-per}
Let $\omega_j^{(-1)}\in H^\vee$ $(0\leq j\leq N)$ be the fundamental weights of $\Delta^{(-1)}$, i.e., 
\ben
\langle \omega_j^{(-1)},\gamma_m^{(-1)}\rangle =\delta_{j,m}.
\een
Using the intersection form we identify $H^{(0)}$ and its dual. Let $\omega_j^{(0)}\in H^{(0)}$ $(1\leq j\leq N)$ be 
the fundamental weights of $\Delta^{(0)}$, i.e., 
\ben
(\omega_j^{(0)}|\gamma_m^{(0)} )= \delta_{j,m},\quad 1\leq j,m\leq N. 
\een
We have the following relation
\ben
\langle \omega_j^{(-1)} ,\tilde{\al}\rangle =(\omega_j^{(0)}|r(\tilde{\al})) - k_j \langle \omega_0^{(-1)} ,\tilde{\al}\rangle, \quad \tilde{\al}\in\Delta,
\een
where $k_j$ $(1\leq j\leq N)$ are the Kac labels defined by 
\ben
\delta = \gamma_0^{(-1)} + \sum_{j=1}^N k_j\gamma_j^{(-1)}.
\een
In terms of the fundamental weights, the splitting of the affine root system from the previous section can be 
stated also as the following isomorphism 
\ben
\Delta^{(-1)}\cong \Delta^{(0)}\times \Z,\quad \tilde{\al}\mapsto (\al,n),\quad 
\al=r(\tilde{\al}),\quad n=\langle\omega_0^{(-1)},\tilde{\al}\rangle.
\een
\begin{lemma}\label{ev-si}
The following identity holds
\ben
\omega_b^{(0)}+\sum_{m=1}^{a_k-1} (\zeta_k^{mp}-\zeta_k^{(m-1)p})\omega_{k,m}^{(0)} = a_k\phi_{k,p^*},
\quad 
1\leq k\leq 3.
\een
\end{lemma}
\proof
We have explicit formulas for the simple roots 
\ben
\gamma_b^{(0)} & = &  \one+\chi P +\sum_{k=1}^3\sum_{p=1}^{a_k-1} \phi_{k,p} \\
\gamma_{k,m}^{(0)} &= & \sum_{p=1}^{a_k-1} (\zeta_k^{mp}-\zeta_k^{(m-1)p})\phi_{k,p}.
\een
It remains only to check that the LHS and the RHS have the same intersection pairing with the above set of simple roots of $\Delta^{(0)}$.
\qed

Let $\kappa$ be a positive constant whose value will be specified later on. Put 
\beq\label{eigen-basis}
H_0:=H_{01}:=H_{02}:= (\kappa\chi)^{\frac{1}{2}}\, \omega_b^{(0)},\quad 
H_i:= (\kappa\, a_i)^{\frac{1}{2}}\, \phi_i,\quad i\in \mathfrak{I}_{\rm tw}.
\eeq
Note that according to Lemma \ref{per-1} $\{H_i\}_{i\in \mathfrak{I}}$ is a $\si_b$-eigenbasis of
$H^{(0)}$ with  $\si_b(H_i)=e^{-2\pi\sqrt{-1} d_i}H_i$ in which the intersection form takes the form
\beq\label{norm}
(H_i|H_j)=\kappa\, \delta_{i,j^*},\quad i,j\in \mathfrak{I},
\eeq
where for $i=j=01\in \mathfrak{I}$ we used that $\omega_b^{(0)}=\chi^{-1}\one+P$. Finally, put 
\ben
\rho_b = -\sum_{(k,p)\in\mathfrak{I}_{\rm tw}} \frac{1}{a_k} \omega_{k,p}^{(0)}.
\een
\begin{proposition}
Let $\tilde{\al}=(\al,n)\in \Delta^{(0)}\times \Z\cong \Delta$, be a vanishing cycle, then
the corresponding calibrated periods are given by the following formulas:
\ben
&&
\phantom{^{--1}}
\tI^{(\ell)}_{\tilde{\al}}(\la) = (-1)^\ell \ell! (\al | \omega^{(0)}_b)  \chi\la^{-\ell-1} P +
\sum_{i\in \mathfrak{I}_{\rm tw}}
(\al|H_{i^*})  (d_i-1)\cdots (d_i-\ell) \la^{d_i-\ell-1}
\sqrt{a_{i}/\kappa} \phi_i ,\\
&&
\phantom{^{--\ell}}
\tI^{(0)}_{\tilde{\al}}(\la) =(\al | \omega^{(0)}_b) \one+ (\al | \omega^{(0)}_b)  \chi\la^{-1} P +
\sum_{i\in \mathfrak{I}_{\rm tw}}
(\al|H_{i^*})  \la^{d_i-1}
\sqrt{a_{i}/\kappa} \phi_i , \\
&&
\tI^{(-1-\ell)}_{\tilde{\al}}(\la) = 
(\al|\omega^{(0)}_b) \frac{\la^{\ell+1}}{(\ell+1)!} \one+\Big(
(\al|\omega^{(0)}_b)  \chi\frac{\la^\ell}{\ell!}(\log \la - C_\ell)\, + 
2\pi\sqrt{-1} (n+(\rho_b|\al) )\, \frac{\la^\ell}{\ell!}\, \Big)\, P+ \\
&&
\phantom{\tI^{(-1-\ell)}_{\tilde{\al}}(\la) = }
\sum_{i\in \mathfrak{I}_{\rm tw}} (\al| H_{i^*})\,
\sqrt{a_{i}/\kappa}\, \frac{\la^{d_i+\ell}}{ d_i(d_i+1)\cdots
  (d_i+\ell)}\, \phi_i,
\een
where $\ell\geq 1$ and $C_\ell$ $(\ell\geq 1)$ are constants defined recursively by
\ben
C_\ell=C_{\ell-1}+\frac{1}{\ell},\quad C_0=\frac{1}{\chi}\log Q.
\een
\end{proposition}
\proof
It is enough to check the statement for the following basis of the Milnor lattice
\ben
\gamma_b^{(-1)},\quad \delta, \quad \gamma_{i}^{(-1)}, \quad i=(k,p)\in\mathfrak{I}_{\rm tw}.
\een
Let us check the last identity for $\tilde{\al}=\gamma_{k,p}^{(-1)}$, i.e., $n=0$ and $\al=\gamma_{k,p}^{(0)}$. 
Recalling the explicit formulas for 
$\tI^{(-1-\ell)}_{\alpha_{k,p}}(\la)$ and $\gamma_{k,p}=\alpha_{k,-p}-\alpha_{k,-p+1}$, we get 
(recall that $d_{k,m}=1-m/a_k$)
\ben
\tI^{(-1-\ell)}_{\gamma_{k,p}}(\la) = -\frac{2\pi\sqrt{-1}}{a_k}\frac{\lambda^\ell}{\ell!} P + \sum_{m=1}^{a_k-1}
\frac{\zeta_k^{mp}-\zeta_i^{m(p-1)}}{( \ell+d_{k,m})\cdots (1+ d_{k,m})} 
\lambda^{\ell+d_{k,m}}\phi_{k,m}.
\een
On the other hand, by definition $H_{i*}\sqrt{a_i/\kappa} = a_i \phi_{i*}$, so the identity follows from Lemma
\ref{ev-si}. The remaining two cases are proved in the same way.
\qed

\section{ADE-Toda hierarchies}\label{sec:basic-rep}

\subsection{Twisted realization of the affine Lie algebra}

Let $ \lieg^{(0)}$ be a simple Lie algebra of type $ADE$ with an invariant bilinear
form $(\ |\ ),$ normalized in such a way that all roots have length $\sqrt{2}$. By definition, the affine Kac--Moody algebra
corresponding to $ {\lieg}$ is the vector space 
$$
\lieg:= {\lieg^{(0)}} [t,t^{-1}] \bigoplus \C\,K \bigoplus \C\, d
$$ 
equipped with a Lie bracket defined by the following relations: 
\ben 
&& 
[X\ t^n,Y\ t^m]:=[X,Y]\ t^{n+m} + n\delta_{n,-m}(X\, |\, Y)  K,\\ 
&& 
[d,X\ t^n]:= n(X\ t^n),\quad [K,\lieg^{(0)}]:=0, 
\een
where $X,Y\in  {\lieg^{(0)}}$.

We fix a Cartan subalgebra  $\lieh^{(0)}\subset \lieg^{(0)}$ and a
basis $\gamma_b,\gamma_i$ $(i\in\lieI)$ of simple roots, s.t., the
corresponding Dynkin diagram has the standard shape with $\gamma_b$
corresponding to the branching node. 
If the root system is of type $A$; then we choose any of
the nodes to be a branching node and we have (at most) 2 instead of 3
branches. Let us define $\si_b:=\si_b^{(0)}$ by formula \eqref{si-0}. 

Let
$\Delta^{(0)}\subset \lieh^{(0)}$ be the root system of $\lieg^{(0)}$, i.e., 
$$
\lieg^{(0)}=\bigoplus_{\al\in \Delta^{(0)}} \lieg^{(0)}_\al.
$$ 
The Lie algebra $\olieg$ can be constructed in terms of the root
system via the so-called {\em Frenkel--Kac} construction \cite{FK}. 
Let $\oLambda\subset \olieh$ be the root lattice.
There exists a bimultiplicative function 
\ben
\epsilon: \Lambda^{(0)}\times  \Lambda^{(0)} \to \{\pm
1\}
\een
satisfying
\ben
\epsilon(\al,\beta)\epsilon(\beta,\al) = (-1)^{(\al|\beta) },\quad 
\epsilon(\al,\al) = (-1)^{|\al|^2/2},
\een
where $|\al|^2:=(\al|\al) .$ 
The map $(\al,\be)\mapsto \ep(\osi(\al),\osi(\be))$ is another
bimultiplicative function satisfying the above properties. It is known that
all bi-multiplicative functions of the above form are equivalent (see
\cite{Kac-v}, Corollary 5.5). Hence there exists a function
$\upsilon:\oLambda\to \{\pm 1\}$ such that
\beq\label{eq:zeta}
\upsilon(\al)\upsilon(\be)\ep(\al,\be) = \upsilon(\al+\be)\ep(\osi(\al),\osi(\be)).
\eeq   
There exists a set of root vectors 
\beq\label{A_alpha}
A_\al \in  \lieg^{(0)}_\al
\eeq 
such that
\begin{align}
\notag
[ A_\al,A_{-\al} ] & =  \epsilon(\al,-\al) \al  \\ 
\notag
[ A_\al,A_\beta ] & =  \epsilon(\al,\beta) A_{\al+\beta},\quad \mbox{if} \quad
(\al|\beta)  = -1 \\
\notag
[ A_\al,A_{\beta} ] & =  0 ,\quad \mbox{if} \quad
(\al|\beta)  \geq 0. 
\end{align}
We can extend $\osi$ to a Lie algebra automorphism of $\olieg$ as
follows
\ben
\osi(A_\al) = \upsilon(\al)^{-1} \, A_{\osi(\al)},\quad \al \in \oDelta.
\een
Let us denote by $\kappa$ the order of the extended automorphism $\osi:\olieg\to \olieg$. Clearly we have
$\kappa=|\osi|$ or $2|\osi|$.  
Since $(\ |\ ) $ is both $\olieg$-invariant (with respect to the
adjoint representation) and $ W^{(0)}$-invariant, we have
\ben
(A_\al|A_{-\al})  := \epsilon(\al,-\al),\quad 
(A_\al|A_\beta)  := (A_{\al}|H)  = 0,\quad \forall \beta\neq
-\al, \ H\in  \lieh^{(0)}.
\een
Put $\eta = e^{2\pi\sqrt{-1}/\kappa}$. We extend the action of $\osi$
to the affine Lie algebra $\lieg$ by 
\ben
 \osi\cdot (X\otimes t^n) =  \osi(X)\otimes
(\eta^{-1} t)^n,\quad
 \osi\cdot K = K,\quad
 \osi\cdot d = d.
\een 
Let $$\lieg^{ \osi}\subset \lieg$$ be the Lie subalgebra of
$ \osi$-fixed points. According to Kac (see \cite{Kac},
Theorem 8.6.), $\lieg^{ \osi} \cong \lieg.$ Let us recall the
isomorphism. The fixed points subspace $ (\lieg^{(0)})^{\osi}$
contains a Cartan subalgebra $ \tlieh^{(0)}$. We have a
corresponding decomposition into root subspaces
\ben
 \lieg^{(0)} = \tlieh^{(0)}\bigoplus\left(\bigoplus_{\tal\in  \tDelta^{(0)}}  \lieg^{(0)}_{\tal}\right),
\een
where $ \tDelta^{(0)}\subset  \tlieh^{(0)}$ are the
corresponding roots. Note that since
the root subspaces are 1 dimensional, they must be eigen-subspaces of
$ \osi$. Therefore, by choosing a set of simple roots
$\tal_j$, $j=1,2,\dots, N$ in $ \tDelta^{(0)}$ we can
uniquely define an integral vector $s=(s_1,\dots,s_{N})$, $0\leq
s_j<\kappa$ such that the eigenvalue of the eigensubspace
$ \olieg_{\tal_j}$ is $\eta^{s_j}$. 
Put
\ben
\rho_s:  \tlieh^{(0)} \to  \tlieh^{(0)},\quad
\rho_s=\sum_{j=1}^{N} s_j\tomega_j,
\een
where $\tomega_j\in \tlieh^{(0)} \, (1\leq j\leq N)$ are the fundamental weights corresponding to the
simple roots $\tal_j \, (1\leq j\leq N)$, i.e.,
$(\tomega_j|\tal_{j'})=\delta_{j,j'}$. The isomorphism 
\ben
\Phi: \lieg\longrightarrow \lieg^{ \osi}
\een
is defined as follows
\beqa\label{tw:iso1}
\Phi(X t^n) & = &   
t^{n\kappa +{\rm ad}_{\rho_s }} X + \delta_{n,0}\,
(\rho_s  | X)K\\
\notag
\Phi(K) & = & \kappa\, K \\
\label{tw:iso2}
\Phi(d) & = & \kappa^{-1}
\Big(d-\rho_s -\frac{1}{2}(\rho_s |\rho_s ) \, K\Big ),
\eeqa
where 
\ben
t^{{\rm ad}_{\rho_s }} X = \exp\Big(\log t\   {\rm ad}_{\rho_s }\Big) X.
\een
Note that the RHS is single-valued in $t$ and $\osi$-invariant in
$X$, because 
\ben
\exp \Big( 2\pi\sqrt{-1} {\rm ad}_{\rho_s /\kappa}\Big)=\osi.
\een

Finally,  we make a remark on $\kappa$. There is no a canonical way to extend $\si_b$
to a Lie algebra automorphism of $\lieg^{(0)}$. Therefore, the value of $\kappa$ depends on our choice of the cocycle
$\epsilon(\al,\be)$ and the corresponding sign-function $\upsilon(\alpha)$. We will see however that replacing $\kappa$ by 
$m\kappa$, where $m$ is a positive integer, does not change the HQEs, so we may assume that $\kappa$ is a sufficiently large integer, s.t., $\si_b^\kappa=1$. For the sake of completeness, let us fix an extension that seems natural for our purposes. 
Put $\omega_{k,0}=\omega_{b}$ and $\omega_{k,a_k}=0$ and define
\beq\label{seifert-form}
{\rm SF}(\al,\beta)=\sum_{k=1}^3 \sum_{p=0}^{a_k-1} (\omega_{k,p}|\al)(\omega_{k,p}-\omega_{k,p+1}|\be).
\eeq
Since $\operatorname{SF}(\al,\be)+ \operatorname{SF}(\be,\al)=(\al|\be)$, the bi-multiplicative function 
$\epsilon(\ ,\ )=(-1)^{\rm SF(\ ,\ )}$ is an acceptable choice for the Frenkel--Kac construction. Note that
\beq\label{zeta-def}
\upsilon(\al)=(-1)^{\sum_{k=1}^{3}(\omega_b|\al)(\omega_{k,1}|\al)}
\eeq
satisfies formula \eqref{eq:zeta}, so we get an explicit formula for an extension of $\si_b$ to a Lie algebra automorphism 
of $\lieg^{(0)}$. Moreover, since
\ben
\prod_{m=1}^{|\sigma_b|}\upsilon(\sigma_b^{m}(\al))=(-1)^{\chi|\sigma_b|},
\een
we get that $\kappa=|\sigma_b|$ if $\chi|\sigma_b|$ is even and $\kappa=2|\sigma_b|$ if $\chi|\sigma_b|$ is odd. Notice that $|\sigma_b|={\rm lcm}(a_1,a_2,a_3)$, the least common multiple of $a_1,a_2,a_3$.
\begin{remark}
The notation $\operatorname{SF}$ is motivated from the notion of a Seifert form in singularity theory (cf. \cite{AGV, BM}). 
We do not claim that \eqref{seifert-form} is a Seifert form, although it would be interesting to investigate whether 
definition \eqref{seifert-form} can be interpreted as a linking number between $\al$ and $\be$. 
\end{remark}

\subsection{The Kac--Peterson construction}\label{kac-peterson}
Following \cite{KP}, we would like to recall the realization of the basic level 1 representation of the affine
Lie algebra $\lieg$ corresponding to the automorphism $\osi$. The idea is to
construct a representation of the Lie algebra $\lieg^{\osi}$ on a Fock
space, which induces via the isomorphism $\Phi$ the basic level-1
representation.  

Fix a $\si_b$-eigenbasis $\{H_i\}_{i\in \{0\}\sqcup\mathfrak{I}_{\rm
    tw} }$ of $\lieh^{(0)}$. It is convenient to define
$H_{01}:=H_{02}:=H_0$ and to assume that the basis is normalized so
that $(H_i|H_{j^*})=\kappa\,\delta_{i,j}$ (compare with \eqref{norm}). 
Put 
\ben
m_{01}:=0,\quad m_{02} := \kappa,\quad m_i:= d_{i^*}\,\kappa, \quad i
\in \mathfrak{I}_{\rm tw},
\een
so that $e^{-2\pi\sqrt{-1}\,d_i}=\eta^{m_i}$ is the eigenvalue
corresponding to the eigen vector $H_i.$ The elements
$$H_{i,\ell}:=H_it^{m_i+\ell\kappa} \, (i\in
\mathfrak{I},\ell\in \Z)$$ 
generate a {\em Heisenberg} Lie subalgebra 
$\lies\subset\lieg^{\osi},$ 
i.e., the following commutation relations hold
\ben
[H_{i,\ell},H_{j,m}] = (m_i+\ell\kappa)\,\delta_{i,j^*}\, \delta_{\ell+m,-1}\,\kappa\, K.
\een
Let us also fix a $\C$-linear basis of $\lies$
\beq\label{s-basis}
H_0:=H_{01},\quad H_{i,\ell},\quad H_{i^*,-\ell-1},\quad K,\quad
(i,\ell)\in I_+,
\eeq
where the index set is defined by
\beq\label{ind-pair}
I_+=\{(i,\ell) \ |\ i\in \mathfrak{I}\setminus{\{(01)\}}, \ell\in\mathbb{Z}_{\geq0}\}.
\eeq
Let $\mathfrak{S}$ be the subgroup of the affine Kac--Moody Lie group
generated by the lifts of the following loops:
\beq\label{h_albe}
h_{\al,\beta} = \exp \Big( \al \log t^{\kappa} + 2\pi\sqrt{-1} \ \beta\Big),
\eeq
where $\al,\be\in \olieh$ are such that 
\ben
\osi(\al)=\al,\quad \osi(\beta)-\beta + \al \in {\Lambda^{(0)}}. 
\een
Let us point out that under the analytical continuation around $t=0$,
the loop $h_{\al,\beta}$ gains the factor $e^{2\pi\sqrt{-1}
  \kappa\al}$. The latter must be 1 because
\ben
\kappa\al = (\al+\osi(\beta)-\beta) + \osi(\al+(\osi(\beta)-\beta)
)+\cdots+ \osi^{\kappa-1}(\al+(\osi(\beta)-\beta)) \in
{\Lambda^{(0)}}. 
\een
It follows that $h_{\al,\beta}$ is single-valued and $\osi$-invariant,
i.e., it defines an element of the affine Kac--Moody loop group acting
on $\lieg^{\osi}$ by conjugation. 
The main result of Kac and Peterson \cite{KP} is the following: \emph{the basic representation of $\lieg^{\osi}$ remains irreducible when restricted to the pair $(\mathfrak{s},
\mathfrak{S})$}.

Let us recall the construction of the representation. 
Let us denote by $$\pi_0:\olieh\to \olieh_0 \text{ and } \pi_*:\olieh\to (\olieh_0)^\perp$$ the orthogonal projections of $\olieh$ onto $\olieh_0:=\C\,H_0$ and $(\olieh_0)^\perp$ respectively. Given $x\in \olieh$, let
\ben
x_0:=\pi_0(x), \quad x_*:=\pi_*(x).
\een
Let $\lies_-\subset\lies$ be the Lie
subalgebra of $\lies$ spanned by the vectors $H_{i^*,-\ell-1}, (i,\ell)\in I_+$. The basic representation can
be realized on the following vector space:
\beq\label{repr}
V_x=S^*(\lies_-)\otimes \C[ e^\omega]e^{x\omega},
\eeq
where $x$ is a complex number and $\omega:=\pi_0(\gamma_b)$. 
The first factor of the tensor product in \eqref{repr} is the symmetric algebra
on $\lies_-$, and the second one is isomorphic to the group algebra of the lattice
$\pi_0 ({\Lambda^{(0)}})=\Z\,\pi_0(\gamma_b).$ We will refer to $|0\rangle
:=1\otimes e^{x\omega}$ as the {\em vacuum vector}. Slightly abusing
the notation we define the operator $$\d_\omega:=\frac{\d}{\d \omega}
-x,$$ acting on $V_x$, so that $\d_\omega\,|0\rangle =0.$

Put
\ben
X_\al(\zeta) = \sum_{n\in \Z} A_{\al,n}\,
\zeta^{-n}=\frac{1}{\kappa}\, \sum_{m=1}^{\kappa} \sum_{n\in \Z} \eta^{-nm}
 ( \osi^m(A_{\al}) t^n) \zeta^{-n},\quad \al \in \oDelta,
\een
where $A_\alpha$ appears in (\ref{A_alpha}). 
Let $E^*_\al(\zeta)$ be the {\em vertex operator}
\beq\label{vop:twisted}
E^*_\al(\zeta)=\exp\Big(
\sum_{(i,\ell)\in I_+} (\al|H_i)  H_{i^*,-\ell-1}
\frac{\zeta^{m_i+\ell\kappa}}{m_i+\ell\kappa}\Big)
\exp\Big(\sum_{(i,\ell)\in I_+} (\al|H_{i^*})  H_{i,\ell} \frac{\zeta^{-m_i-\ell\kappa}}{-m_i-\ell\kappa}\Big).
\eeq
\begin{lemma}\label{lem:C_al}
There are operators $C_\al$, $\al\in \oDelta$, independent of $\zeta$, that commute with 
all basis vectors \eqref{s-basis} of $\mathfrak{s}$ different from $H_0$, such that 
\ben
X_\al(\zeta)
=X^0_\al(\zeta) E^*_\al(\zeta), 
\een
where 
\beq\label{heis:factor}
X^0_\al(\zeta) = \zeta^{\kappa\,|\al_0|^2/2} C_\al \zeta^{\kappa\,
  \al_0}, \quad \alpha_0:=\pi_0(\alpha).
\eeq
\end{lemma}
\proof

After a direct computation we get
\ben
[H_{i,\ell}, X_\al(\zeta)] = (\al|H_i) \zeta^{m_i+\ell\kappa} X_{\al}(\zeta).
\een
It follows that $X_\al(\zeta) = X^0_\al(\zeta) E^*_\al(\zeta)$, where
$X^0_\al(\zeta)$ is an operator commuting with all $H_{i,\ell}\neq H_0$.

After a direct computation we get the following commutation relations:
\ben
h_{\al, \beta}\  (-d)\  h_{\al, \beta}^{-1} & =
&-d +  \kappa\, \al +
\frac{1}{2}|\al|^2 \,\kappa^2K, \\
h_{\al, \beta}\ A_{\gamma,n} \ h_{\al, \beta}^{-1} & = & 
e^{2\pi\sqrt{-1}\,(\beta|\gamma)} A_{\gamma,n+\kappa(\al|\gamma)} +
\delta_{n,0} (\al|A_\gamma)\kappa\, K, 
\een
and $h_{\al,\beta}$ commute with the Heisenberg algebra $\mathfrak{s}$
except for:
\ben
h_{\al, \beta}\ H_0\ h_{\al, \beta}^{-1} & = & 
H_0+(\al|H_0)\kappa K.
\een
Here $h_{\al,\beta}$ are given in (\ref{h_albe}). In order to determine the dependence on $\zeta$ of $X^0_\al(\zeta)$ we
first have to notice that 
\beq\label{vir:0}
-d=\frac{1}{2}|\rho_s |^2 K+\frac{1}{2} H_0^2 +\sum_{(i,\ell)\in I_+} H_{i^*,-l-1} H_{i,\ell},
\eeq
where $H_0=H_{01} = H_{02}.$ Indeed, if we
decompose the basic representation into a direct sum of weight subspaces of
$\mathfrak{s}$, then using the above commutation relations, we get that the LHS of \eqref{vir:0} is an operator that preserves these
weight subspaces while the difference of the LHS and the RHS commutes
with $\mathfrak{s}$ and $\mathfrak{S}$. The formula follows up to the constant
term $\frac{1}{2}|\rho_s |^2 K$, which is fixed by examining the
action of the operator $d\in \lieg$ on the vacuum vector. Using
formula \eqref{tw:iso1}  for $Xt^n=\rho _s$
we get that $\rho_s $ (viewed as an element of $\lieg^{\osi}$) acts on the
vacuum by the scalar $-|\rho_s|^2/\kappa$; then since the RHS of
formula \eqref{tw:iso2} acts by 0 on the vacuum, we get that $d\in
\lieg^{\osi}$ acts by the scalar
\ben 
-|\rho_s|^2/\kappa + \frac{1}{2}|\rho_s|^2 (1/\kappa) = - \frac{1}{2\kappa}|\rho_s|^2. 
\een
Since we have
\ben
[d,X_\al(\zeta)] = -\zeta \d_{\zeta} X_\al(\zeta),\quad
[d,E_\al^*(\zeta)]=-\zeta\d_\zeta E_\al^*(z),
\een
we easily get $-\zeta\d_\zeta X^0_\al = [d,X^0_\al].$ On the other hand,
$X^0_\al(\zeta)$ commutes with $H_{i,\ell}$ for all $i,\ell$, except 
\beq\label{h0-commutator}
[H_0, X^0_\al(\zeta)]= (\al|H_0) X^0_\al.
\eeq
It follows that
\ben
\zeta \d_\zeta X^0_\al =\kappa\Big( X^0_\al \al_0+
\frac{|\al_0|^2}{2} X^0_\al\Big). 
\een
Solving the above equation we get formula \eqref{heis:factor}.
\qed

\begin{lemma}
The operators $C_\al$ in \eqref{heis:factor} satisfy the following commutation relation
\beq\label{mult}
C_\al C_\beta = \epsilon(\al,\beta) B_{\al,\beta}^{-1} C_{\al+\beta} ,
\eeq
where
\ben
B_{\al,\beta} = \kappa^{-(\al|\beta) }\, \prod_{m=1}^{\kappa-1}
(1-\eta^m)^{(\osi^m(\al)|\beta) }. 
\een
\end{lemma}
\proof
 Let us assume first that
$\al\neq -\beta$ are two roots. After a direct
computation we get that  the commutator
$[X_{\al}(\zeta),X_\beta(w)]$ is given by the following formula:
\ben
\frac{1}{\kappa}\, \sum_{m=1}^{\kappa} \left(\prod_{j=1}^{m-1}\upsilon^{-1}(\osi^{j}(\beta))\right)\epsilon(\al,\osi^m(\beta))\,
\delta(\eta^{-m}\zeta,w) \, w \, X_{\al+\osi^m(\beta)}(\zeta),
\een
where $\delta(x,y):=\sum_{n\in \Z} x^n y^{-n-1}$ is the formal delta
function. On the other hand,
\ben
E_\al^*(\zeta)E^*_\beta(w) = \prod_{m=1}^{\kappa}
\Big(1-\eta^m\frac{w}{\zeta}\Big)^{(\osi^m(\al)|\beta) }
:E_\al^*(\zeta)E^*_\beta(w):\, ,
\een
where $:\ :$ is the standard normal ordering in the Heisenberg group -
all annihilation operators $H_{i,\ell}$ must be moved to the
right. Substituting in the above commutator $X_\gamma(\zeta) =
X^0_\gamma(\zeta)E^*_\gamma(\zeta)$ we get that the following two
expressions are equal:
\beq\label{expr:1}
\prod_{m=1}^{\kappa}
\Big(1-\eta^m\frac{w}{\zeta}\Big)^{(\osi^m(\al)|\beta) }
X^0_\al(\zeta)X^0_\beta(w) - 
\prod_{m=1}^{\kappa}
\Big(1-\eta^m\frac{\zeta}{w}\Big)^{(\osi^m(\beta)|\al) }
X^0_\beta(w)X^0_\al(\zeta)
\eeq
and
\beq\label{expr:2}
\frac{1}{\kappa}\, \sum_{m=1}^{\kappa} \left(\prod_{j=1}^{m-1}\upsilon^{-1}(\osi^{j}(\beta))\right) \epsilon(\al,\osi^m(\beta))
\delta(\eta^{-m}\zeta,w)\, w \, X^0_{\al+\osi^m(\beta)}(\zeta).
\eeq
Both formulas have the form $i_{\zeta,w} P_1(\zeta,w) - i_{\zeta,w}
P_2(\zeta,w)$, where $P_1$ and $P_2$ are some rational functions and
$i_{\zeta,w}$ (resp. $i_{w,\zeta}$) means the Laurent series expansion in the region
$|\zeta|>|w|$ (resp. $|w|<|\zeta|$). Since $P_1=P_2$ for the second
expression,  the same must be true for the first one, i.e.,
\ben
\prod_{m=1}^{\kappa}
\Big(1-\eta^m\frac{w}{\zeta}\Big)^{(\osi^m(\al)|\beta) }
X^0_\al(\zeta)X^0_\beta(w) = 
\prod_{m=1}^{\kappa}
\Big(1-\eta^m\frac{\zeta}{w}\Big)^{(\osi^m(\beta)|\al) }
X^0_\beta(w)X^0_\al(\zeta).
\een
Recalling formula \eqref{heis:factor} and \eqref{h0-commutator}, the above equality implies:
\beq\label{com}
C_\al C_\beta = \prod_{m=1}^{\kappa}
(-\eta^m)^{(\al|\osi^m(\beta)) }\ C_\beta C_\al.
\eeq
Using this equality we can easily write \eqref{expr:1} as a sum of
formal delta functions. Comparing with \eqref{expr:2} we get
\eqref{mult}. 
\qed

\begin{lemma}\label{br-star}
Let $\omega_b$, $\omega_i$, $i=(k,p)\in\mathfrak{I}_{\rm tw}$, be the
fundamental weights corresponding to the basis of simple roots $\gamma_b$, $\gamma_i$, $i\in\mathfrak{I}_{\rm tw}$, then 
\ben
(\omega_{i}|\chi \omega_b) = d_i ,\quad 
\pi_0(\gamma_b) = \chi\omega_b,\quad
\pi_*(\gamma_b) = -\sum_{i\in\mathfrak{I}_{\rm tw}}  d_i \, \gamma_{i}.
\een
\end{lemma}
\proof
Let $\{\ep_{k,p}\}_{p=1}^{a_k}$ be the standard basis of $\C^{a_k}$ for any fixed $k=1,2,3$. The root
system of type $A_{a_k-1}$ is given by $\{\ep_{k,p}-\ep_{k,q}\}$ and the
standard choice of simple roots is $\gamma_{k,p}=\ep_{k,p}-\ep_{k,p+1}$, $1\leq
p\leq a_k-1.$ Note that the fundamental weights corresponding to the
basis of simple roots are
\ben
\widetilde{\omega}_{k,p} = 
\Big(1-\frac{p}{a_k}\Big)(\ep_{k,1}+\cdots+\ep_{k,p})- 
\frac{p}{a_k}(\ep_{k,p+1}+\cdots+\ep_{k,a_k}).
\een
It follows that the pairing between the fundamental weights is 
\ben
(\widetilde{\omega}_{k,p} |\widetilde{\omega}_{k,q} ) = {\rm min}(p,q)-pq/a_k. 
\een
In particular, we have 
\beq\label{a-fund}
\widetilde{\omega}_{k,p}  = \Big(1-\frac{p}{a_k}\Big)\gamma_1 + \cdots,
\eeq
where the remaining terms involve only $\gamma_2,\dots,\gamma_{a_k-1}$. 

In our settings, the roots
$\{\gamma_{k,p}\}_{p=1}^{a_k-1}$ give rise to a subroot system of
type $A_{a_k-1}$. Let us denote by $\widetilde{\omega}_{k,p}$ the
corresponding fundamental weights. Note that 
\ben
\omega_{k,p}=\widetilde{\omega}_{k,p} -
(\widetilde{\omega}_{k,p}|\gamma_b)\, \omega_b,
\een
so the first formula of the Lemma follows from \eqref{a-fund} and
\ben
(\gamma_b|\gamma_{k,p})=-\delta_{p,1},\quad (\widetilde{\omega}_{k,p}|\omega_b)=0.
\een 
The other two identities follow easily from the first one. 
\qed

Using formula \eqref{mult} we define $C_\al$ for all $\al$ in the root
lattice ${\Lambda^{(0)}}$; then formula \eqref{com} still
holds. Finally, a similar argument gives us that
\beq\label{c0}
C_\al C_{-\al} = \epsilon(\al,-\al)\,
B_{\al,-\al}^{-1},\mbox {i.e., } C_0=1. 
\eeq

\begin{lemma}\label{const-c}
Let $c_\al$ ($\al\in \oLambda$) be operators defined by 
\beq\label{cal}
C_\al = c_\al 
\exp \Big( (\omega_b|\al) \omega\Big)
\exp \Big(2\pi\sqrt{-1}\, (\rho_b |\al)  \d_\omega
\Big). 
\eeq
Then $[c_\al,c_\be]=0$.
\end{lemma}
\proof
To begin with, note that by definition, the commutator $C_\al C_\beta
C_\al^{-1}C_\beta^{-1}$ is given by the following formula:
\ben
\prod_{m=1}^{\kappa}
(-\eta^m)^{(\al|\osi^m(\beta)) } = e^{\pi \sqrt{-1} \, (\al_0|\beta) }
e^{2\pi \sqrt{-1} \,( (1-\osi)^{-1}\al_*|\beta ) }. 
\een
On the other hand, using \eqref{cal}, the commutator becomes
\beq\label{comm-c}
c_\al c_\be c_\al^{-1} c_\be^{-1} 
\exp \, 2\pi\sqrt{-1}\, 
\Big(
(\rho_b |\al) (\omega_b|\be) -
(\rho_b|\be) (\omega_b|\al)
\Big).
\eeq
Recall that $\osi$ is a composition of 3 matrices $\si^{(0)}_k, k=1,2,3$
whose action on the subspace with basis
$\{\gamma _{k,1},\dots, \gamma _{k,a_k-1}\}$ is
represented by the matrix 
\ben
\si^{(0)}_k=
\begin{bmatrix}
 -1      &          1 & \cdots &0 & 0 \\
 -1      &          0 & \ddots &0 & 0 \\
\vdots & \vdots & \ddots &  \ddots  & \vdots \\
-1       &          0 & \cdots & 0 & 1 \\
-1       &          0 & \cdots & 0 & 0
\end{bmatrix}.
\een
It is easy to check that the
$(p,q)$-th entry 
\beq\label{pq-entry}
\Big[(1-\si^{(0)}_k)^{-1}\Big]_{pq}=\frac{p}{a_k}-\ep_{pq},\quad
\ep_{pq}=
\begin{cases}
0 & \mbox{ if } p\leq q ,\\
1 & \mbox{ if } p>q.
\end{cases}
\eeq
A straightforward computation using formula \eqref{pq-entry} and Lemma
\ref{br-star} yields
\ben
((1-\si^{(0)}_k)^{-1} \gamma_{k,p} |\gamma_{b} )  & = &
-\frac{1}{a_k}, \\
((1-\si^{(0)}_k)^{-1} \gamma_{k,p} |\gamma_{k,q} )  & = & \delta_{p,q}-\delta_{p+1,q}, \\
((1-\si_b)^{-1} (\gamma_b)_* |\gamma_{k,q} )  & = &  \frac{1}{a_k}
\quad ({\rm mod}\ \Z), \\
((1-\si_b)^{-1} (\gamma_b)_*  |\gamma_{b} ) & = & 1-\frac{1}{2}\chi\,. 
\een
Using the above formulas we get
\ben
((1-\osi)^{-1}\pi_*(\al)|\beta)  = 
(\rho_b|\al) \, (\omega_b|\beta) 
-(\rho_b |\beta) \, (\omega_b|\al)   
-\frac{1}{2} (\al_0|\beta_0) \quad ({\rm mod}\ \Z).
\een
For the commutator we get
\ben
C_\al C_\beta C_\al^{-1} C_\beta^{-1} = \exp \Big(2\pi\sqrt{-1}\Big(
(\rho_b|\al) \, (\omega_b|\beta) -
(\rho_b |\beta) \, (\omega_b|\al)   
\Big)\Big).
\een
Comparing with \eqref{comm-c} we get $c_\al c_\be
c_\al^{-1} c_\be^{-1}=1$.
\qed

Lemma \ref{const-c} implies that the operators $c_\al$ can be
represented by scalars, i.e., we can find complex numbers $c_\al$,
$\al\in \oLambda$ such that
\beq\label{mult-ca}
c_\al c_\beta = \epsilon(\al,\beta)\,B^{-1}_{\al,\beta} \,
e^{-2\pi\sqrt{-1} (\rho_b|\beta) (\omega_b|\al)}\
c_{\al+\beta}. 
\eeq
For example we can choose $c_{\al_i}$ arbitrarily for the simple roots
$\al_i$ and then use formula \eqref{mult-ca} to define the remaining
constants.

The level $1$ basic representation can be realized on $V_x$
as follows. Let us represent the Heisenberg 
algebra $\lies$ on $\C[e^\omega]e^{x\omega}$ by letting all generators
act trivially, except for $H_0\mapsto (H_0|\gamma_b)\,\d_\omega$. The
latter is forced by the commutation relation
\ben
[H_0,C_\al]=(\al|H_0)\,C_\al = (\omega_b|\al)\, (H_0|\gamma_b)\, C_\al.
\een
In this way
$V_x$ naturally becomes  an $\lies$-module. Furthermore, put
\beq\label{E0}
E^0_\al(\zeta) = 
\exp\Big( (\omega_b|\al) \omega \Big)
\exp
\Big(\Big( (\omega_b|\al)\,\chi
\, \log\zeta^{\kappa}\  + 2\pi\sqrt{-1} \ (\rho_b |\al)\Big)\d_\omega \Big) 
\eeq
and $E_\al(\zeta) = E_\al^0(\zeta)E_\al^*(\zeta)$, where
$E_\al^*(\zeta)$ is defined by formula \eqref{vop:twisted}. Thus the representation of the Heisenberg algebra $\lies$ on $V_x$ can
be lifted to a representation of the affine Lie algebra  $\lieg^{\sigma_b}$ as follows:
\ben
X_\al(\zeta) & \mapsto  & c_\al \zeta^{\kappa\,|\al_0|^2/2} E_\al (\zeta), \quad \al\in \oDelta \\
K & \mapsto & 1/\kappa, \\
d & \mapsto & -\frac{1}{2\kappa}|\rho_s|^2-\frac{1}{2} H_0^2 -\sum_{(i,\ell)\in I_+} H_{i^*,-\ell-1} H_{i,\ell}.
\een

\subsection{The Kac--Wakimoto hierarchy}\label{sec:KW}

Following Kac--Wakimoto (see \cite{Kac}), we can define an integrable hierarchy in the
Hirota form whose solutions are parametrized by the orbit of the
vacuum vector $|0\rangle$ of the affine
Kac--Moody group. A vector $\tau\in V_x$ belongs to the orbit if and only if
$\Omega_x\, (\tau\otimes \tau )=0$, where  $\Omega_x$ is the operator
representing the following bi-linear Casimir operator:
\ben
\sum_{\al\in \oDelta} \sum_n \frac{A_{\al,n}\otimes A_{-\al,-n}}{(A_\al|A_{-\al})}\,
   + K\otimes d + d\otimes K +
  \frac{H_0\otimes H_0}{\kappa}  + 
\sum_{(i,\ell)\in I_+} \Big( \frac{H_{i,\ell}\otimes H_{i^*,-\ell-1} +
H_{i^*,-\ell-1}\otimes H_{i,\ell} }{\kappa}\Big) ,
\een
On the other hand, we have
\ben
\sum_n \frac{A_{\al,n}\otimes A_{-\al,-n}}{(A_\al|A_{-\al})}\,
   = {\rm Res}_{\zeta =0} \
  \frac{d\zeta}{\zeta} a_{\al}(\zeta)
  E_\al(\zeta)\otimes E_{-\al}(\zeta), 
\een
where the coefficients $a_\al$ can be computed explicitly thanks to
formula \eqref{mult-ca}, i.e., 
\beq\label{a-al}
a_\al(\zeta)= B_{\al,\al} \, \zeta^{\kappa\,|\al_0|^2}\,
e^{2\pi\sqrt{-1} (\rho_b|\al)(\omega_b|\al)}.
\eeq
We identify the symmetric algebra $S^*(\lies_-)$ with the Fock space
$\C[y]$, where $y=(y_{i,\ell})$ is a sequence of formal variables indexed
by $(i,\ell)\in I_+$ as defined in \eqref{ind-pair}, by identifying
$H_{i^*,-\ell-1}=(m_i+\ell\kappa) y_{i,\ell}$. Then (note that $(H_0|\gamma_b)=(\kappa\chi)^{1/2}$)
\ben
H_{i,\ell}=\frac{\d}{\d y_{i,\ell}},\quad H_0= (\kappa\chi)^{1/2} \,
\d_\omega ,\quad K=1/\kappa, 
\een
and 
\ben
d=-\frac{|\rho_s |^2}{2\kappa}  -\frac{\kappa\chi}{2} 
\d_\omega^2 -\sum_{(i,\ell)\in I_+} (m_i+\ell\kappa)y_{i,\ell} \d_{y_{i,\ell}} .
\een
The elements in $V_x$ can also be thought  as
sequences of polynomials in the following way:
\ben
V_x\cong \C[y]^\Z,\quad \sum_{n\in \Z} \tau_n(y)e^{(n+x)\omega}\mapsto
\tau:= (\tau_n(y))_{n\in \Z}.
\een
The above isomorphism turns $\C[y]^\Z$ into a module over the algebra
of differential operators in $e^\omega$:
\ben
(e^\omega\cdot \tau)_n = \tau_{n-1},\quad (\d_\omega\cdot \tau)_n =
n\tau_n. 
\een
The HQEs of the $\si_b$-twisted Kac--Wakimoto hierarchy  will assume the form (\ref{HQE_KW}) stated in Section
\ref{subsec:basic-rep} provided we prove the following identity. 
\begin{lemma}\label{constant}
The following identities hold 
\ben
|\rho_s|^2/\kappa^2 = \frac{1}{12}\sum_{k=1}^3 \frac{a_k^2-1}{a_k} = \frac{1}{2}\operatorname{tr}\Big(
\frac{1}{4}+\theta\, \theta^{\rm T}\Big),
\een 
where $\theta$ is the Hodge grading operator \eqref{eq:theta}.
\end{lemma}
\proof
Since $\tau=|0\rangle$ is a solution to the hierarchy, we must have
\ben
|\rho_s|^2/\kappa^2 = \sum_{\al: \ (\omega_b|\al)=0} a_\al(\zeta).
\een
Let $\al\in \Delta^{(0)}$ be such that $(\omega_b|\al)=0$, then
formula \eqref{a-al} reduces simply to 
\ben
a_\al(\zeta) = B_{\al,\al} = \kappa^{-2} \prod_{m=1}^{\kappa-1}
(1-\eta^m)^{(\si_b^m(\al)|\al)}. 
\een
Recall the notation in the proof of Lemma \ref{br-star}. We claim that
$\al$ must belong to one of the root subsystems $\Delta_k^{(0)}$ of
type $A_{a_k-1}$ corresponding to the legs of the Dynkin diagram
for some $k$. Indeed, let us write $\al$ as a linear combination
$\sum_{k,p} c_{k,p}\gamma_{k,p}$ for some integers
$c_{k,p}$. If this linear combination involves a simple root
$\gamma_{k,p}$ for some $k$, then using
reflections $s_{k,p}$ with $p>1$ we can 
transform $\al$ to a cycle $\al'$ such that the decomposition of $\al'$
as a sum of simple roots will involve $\gamma_{k,1}$. Moreover, 
we still have $(\omega_b|\al')=0$. In other words, we may assume that $c_{k,1}\neq
0$ as long as $c_{k,p}\neq 0$ for some $p$. However, since
$(\al|\gamma_b) = -\sum_k c_{k,1}$ and the coefficients
$c_{k,p}$ have the same sign (depending on whether $\al$ is a
positive or a negative root) we get that there is precisely one $k$
for which $c_{k,1}\neq 0$. 

Assume that $\al\in \Delta_k^{(0)}$, then since $\si_b$ is a product
of the Coxeter transformations $\si_{k'}=\cdots
s_{k',2}s_{k',1}$, in the above formula for $a_\al$ only $\si_k$
contributes and since the order of $\si_k$ is $a_k$, after a short
computation we get
\ben
a_\al(\zeta) = a_k^{-2} \prod_{m=1}^{a_k-1}
(1-\eta_k^m)^{(\si_k^m(\al)|\al)},\quad \eta_k =
e^{2\pi\sqrt{-1}/a_k}. 
\een
These are precisely the coefficients of the principal Kac-Wakimoto
hierarchy  of type $A_{a_k-1}$. Let $\rho_k$ be the sum of the fundamental weights of
$\Delta_k^{(0)}$. It is well known that $|\rho_k|^2 =
(a_k-1)a_k(a_k+1)/12$. 
According to \cite{FGM} the sum
\ben
|\rho_s|^2/\kappa^2 =\sum_{\al\in \Delta^{(0)}_k} a_{\al}(\zeta) = |\rho_k|^2/a_k^2
= \frac{1}{12} \sum_{k=1}^3 \Big(a_k -
\frac{1}{a_k}\Big). 
\een 
It remains only to notice (using $\theta^{\rm T}=-\theta$) that
\ben
\frac{1}{2} \operatorname{tr}\Big(
\frac{1}{4}+\theta\, \theta^{\rm T}\Big) & = & \frac{1}{2} \operatorname{tr}\Big(
\frac{1}{2}+\theta\Big)\Big(\frac{1}{2}-\theta\Big)  = \frac{1}{2} \sum_{i\in\mathfrak{I}_{\rm tw}}
d_i(1-d_i) 
= \frac{1}{12} \sum_{k=1}^3 \Big(a_k -
\frac{1}{a_k}\Big). 
\een
\qed

\subsection{Formal discrete Laplace transform}
Let $\al\in \Delta^{(0)}$ and $\tilde{\al}\in \Delta$ be as in Section \ref{cal-per}. 
We would like to compare the vertex operators $E_\al(\zeta)$ and
$\tGamma^{\tilde{\al}}(\la):=e^{(\tf_{\tilde{\al}}(\la;z))\sphat}$, where $(-)\sphat\ $ is the quantization operation defined in Section \ref{sec:can_quantize} and
\ben
\tf_{\tilde{\al}}(\la;z) = \sum_{n\in \Z} \, \tI^{(n)}_{\tilde{\al}}(\la)\, (-z)^n, 
\een
see (\ref{eqn:I^n}). Using the formulas for the calibrated periods from Section
\ref{cal-per} we get
\ben
\tGamma^{\tilde{\al}}(\la) = U_{\tilde{\al}}(\la) \, \tGamma_0^{\tilde{\al}}(\la)\, \tGamma_*^{\tilde{\al}}(\lambda), 
\een
where (we dropped the superscript and set $\omega_b:=\omega_b^{(0)}$)
\ben
U_{\tilde{\al}}(\la) &=& \exp\Big(\sum_{\ell=1}^\infty \Big( (\omega_b|\al)
\chi (\log\la-C_\ell) +2\pi\sqrt{-1} (n+(\rho_b|\al)  )\Big) \frac{\la^\ell}{\ell!}\,
q_\ell^{01}/\sqrt{\hbar} \Big),\\
\tGamma_0^{\tilde{\al}}(\la) & = &
\exp\Big(\Big( (\omega_b|\al)\,\chi
\, (\log\, \la - C_0) +2\pi\sqrt{-1} (n+(\rho_b|\al)  )\Big) \,
\frac{q_0^{01}}{\sqrt{\hbar}} \Big)\times 
\exp\Big(-(\omega_b|\al)\, \sqrt{\hbar} \frac{\d}{\d q_0^{01}} \Big),\\
\tGamma_*^{\tilde{\al}}(\la) & = &
\exp\Big( \sum_{(i,\ell)\in I_+} (\al|H_i)  \,
\zeta^{m_i+\ell\kappa}\, y_{i,\ell} \Big)\ 
\exp\Big( \sum_{(i,\ell)\in I_+} (\al|H_{i^*})  \,
\frac{\zeta^{-m_i-\ell\kappa}}{-m_i-\ell\kappa} \, \frac{\d}{\d y_{i,\ell}} \Big),
\een
where $\la=\zeta^{\kappa}/\kappa$, and we use the change of
variables 
\beqa
\label{dv:change1}
y_{02,\ell} & = &\frac{1}{\sqrt{\hbar}}\, \frac{\kappa^{d_{02}}}{\sqrt{\kappa\chi}}\, 
\frac{ q_\ell^{02} }{m_{02}(m_{02}+\kappa)\cdots(m_{02}+\ell\kappa)}, \\
\label{dv:change2}
y_{i,\ell} & = & \frac{1}{\sqrt{\hbar}}\, 
\frac{  \kappa^{d_i} }{\sqrt{  \kappa a_{i} } }\, 
\frac{ q_\ell^i }{m_i(m_i+\kappa)\cdots(m_i+\ell\kappa)}, \quad (i,\ell)\in\mathfrak{I}_{\rm tw}\times\mathbb{Z}_{\geq0}.
\eeqa
Comparing with \eqref{vop:twisted} and \eqref{E0} we get that
$\tGamma^{\tilde{\al}}_*(\la)=E^*_\al(\zeta)$ and that $\tGamma^{\tilde{\al}}_0(\la)$ is a Laplace
transform of $E_\al^0(\zeta)$. We make the last statement precise as follows. Put 
\ben
\widehat{V}:=\C_\hbar[\![y,x,q^{01}_1+1,q^{01}_2,\dots ]\!]^\Z.
\een
The space $\widehat{V}$ contains a completion of the basic
representation $V_x$. It has also some additional variables $q_\ell^{01}$,
$\ell\geq 1$ which will be treated as parameters.  Just like before, we
identify the elements of $\widehat{V}$ with
formal Fourier series
\ben
f=(f_n)_{n\in \Z} \mapsto \sum_{n\in \Z}\, f_n\, e^{(n+x)\omega}.
\een
Given $f(\hbar;\q)\in \C_\hbar[\![\q]\!]$ satisfying the condition
\beq\label{lt-cond}
f(\hbar;\q)\Big|_{q_0^{01} = x\sqrt{\hbar}}\in
\C_\hbar[\![\q]\!]\quad \forall x\in \C,
\eeq
define the formal Laplace transform of $f$ depending on a
parameter $C$ $(C\neq 0)$
\ben
\cF_C(f(q_0^{01},\dots)) := \sum_{n\in \Z}
f ((x+n)\sqrt{\hbar} ,\dots)\, 
e^{(n+x)\omega} \, C^{\frac{1}{2} n^2}\quad \in \widehat{V},
\een
where the dots stand for the remaining $\q$-variables on which $f$
depends. 
It is easy to check that 
\beq\label{var_ft:1}
\cF_C\ \circ \ q_0^{01}/\sqrt{\hbar} = \frac{\d}{\d \omega}\ \circ\ \cF_C
\eeq
and
\beq\label{var_ft:2}
\cF_C\ \circ\  e^{-m\sqrt{\hbar} \d/\d q_0^{01} }=
e^{m\omega}\, C^{\frac{1}{2}m^2+m\d_\omega}\
\circ\  \cF_C,
\eeq
where recall that $\d_\omega =  \frac{\d}{\d \omega} - x.$ 
\begin{lemma}\label{vop-lt}
Let $C=\kappa^\chi\, e^{\chi C_0}$, then 
\ben
\ E_\al^0(\zeta)\ \cF_C = \cF_C\  e^{-AB-\frac{1}{2}B^2\log C}\,
e^{Ax}\, \tGamma_0^{\tilde{\al}},
\een
where
\ben
A=(\omega_b|\al)\,\chi
\, (\log\, \la - C_0) +2\pi\sqrt{-1} (n+(\rho_b|\al)),\quad
B=(\omega_b|\al).
\een
\end{lemma}
\proof
Using \eqref{var_ft:1} and \eqref{var_ft:2} we get that the vertex
operators in $q_0^{01}$ transform as follows:
\ben
\cF_C\ e^{A\,q_0^{01}/\sqrt{\hbar}}\,e^{-B\sqrt{\hbar}\d/\d q_0^{01}}
= e^{AB+\frac{1}{2}B^2\log C}\, e^{Ax}\, e^{B\omega}e^{(A+B\log
  C)\d_\omega}\ \cF_C.
\een
On the other hand, after a straightforward computation, we get  
\ben
e^{AB+\frac{1}{2}B^2\log C} =
\zeta^{\kappa|\al_0|^2}e^{-\frac{|\al_0|^2}{2\chi}\,(2\chi (C_0+\log\kappa)- \log
  C)}e^{2\pi\sqrt{-1}(\omega_b|\al)(\rho_b|\al)}   
\een
and 
\beq\label{n}
A+B\log C = (\omega_b|\al)
\Big(\chi \log \zeta^\kappa +\log C-
\chi(C_0+\log \kappa) \Big) +2\pi\sqrt{-1} (n+(\rho_b|\al)).
\eeq
Furthermore, note that since the operator
$e^{2\pi\sqrt{-1} \d_\omega}$ acts as the identity on
$\widehat{V}$, the integer $n$ in \eqref{n} may be set to 0.
Finally, it remains only to compare with \eqref{E0} and to recall our
assumption 
\beq\label{C}
\log C = \chi (C_0+\log\kappa).
\eeq
\qed

\subsection{Integrable hierarchies for the affine cusp polynomials}\label{sec:hqe}
For every root $\al\in \oDelta\subset H^{(0)}$ we fix an arbitrary lift $\tilde{\al}\in
\Delta\subset \lieh$ (cf. Section \ref{cal-per}). The subset of affine roots
obtained in this way will be denoted by  $\Delta'$. 
Following the construction of Givental and Milanov in \cite{GM} we
introduce the following Casimir-like operator 
\ben
\widetilde{\Omega}_{\Delta'}(\la) & = & 
-\frac{\la^2}{2}\left(\sum_{m=1}^N\, 
:(\tphi_m(\la)\otimes_\liea 1-1\otimes_\liea\tphi_j(\la)) (\tphi^m(\la)\otimes_\liea
1-1\otimes_\liea\tphi^m(\la)) :\, \right) +\\
&&
+
\sum_{\tilde{\al}\in \Delta'}
\tilde{b}_{\tilde{\al}}(\la)\tGamma^{\tilde{\al}}(\la)\otimes_\liea
\tGamma^{-\tilde{\al}}(\la) 
-\frac{1}{2}\operatorname{tr}\left(\frac{1}{4}+\theta\, \theta^{\rm T}\right),
\een
where the notation is as follows. Let $\{{\tilde{\al}}_m\}_{m=1}^N$ and
$\{{\tilde{\al}}^m\}_{m=1}^N$ be two sets of vectors in $\lieh$ such that under the
projection $\tI^{(0)}(1): \lieh\to H^{(0)}$ they project to
bases dual with respect to the intersection form $(\cdot |\cdot)$, i.e.,
$(\tilde{\al}_j|{\tilde{\al}}^m) = \delta_{j,m}$. Then
\ben
\tphi_m(\la) = (\d_\la\,\tf_{{\tilde{\al}}_m}(\la;z))\sphat\ ,\quad 
\tphi^m(\la) = (\d_\la\,\tf_{{\tilde{\al}}^m}(\la;z))\sphat\ ,\quad 1\leq m\leq N.
\een
The tensor product is over the polynomial algebra
$\liea := \C_\hbar[q_1^{01},q_2^{01},\dots],$
which in particular means that almost all terms that involve $\log \la$
cancel. 

The first sum in the definition of $\widetilde{\Omega}_{\Delta'}$ is monodromy
invariant around $\la=\infty$ and hence it expands in only integral
powers of $\la$. In fact one can check that the corresponding
coefficients give rise to a representation of the Virasoro algebra,
which can be identified with an instance of the so called {\em coset
  Virasoro construction}\footnote{We are thankful to B. Bakalov for this
remark.}. After a straightforward computation using the formulas for the
periods from Section  \ref{cal-per}, we get the following formula for 
the coefficient in front of $\la^{-2}$ (i.e., the $L_0$-Virasoro operator)
\ben
\frac{\chi}{2\hbar}\, (q_0^{01}\otimes_\liea 1-1\otimes_\liea q_0^{01} )^2
+\sum_{(i,\ell)\in I_+}   \left( \frac{m_i}{\kappa}+\ell\right) (q_\ell^i\otimes_\liea
1-1\otimes_\liea q_\ell^i )(\d_{q_\ell^i}\otimes_\liea 1-1\otimes_\liea \d_{q_\ell^i}).
\een 
The coefficient $\tilde{b}_{\tilde{\al}}$ are defined in 
terms of the vertex operators $\tGamma^{\tilde{\al}}(\la)$ as follows
\beq\label{coef:b2}
\tilde{b}_{\tilde{\al}}^{-1}(\la) =\lim_{\mu\to \la} 
\ \Big(1-\frac{\mu}{\la}\Big)^2
\, \widetilde{B}_{{\tilde{\al}},-{\tilde{\al}}}(\la,\mu), \quad \tilde{\al}, \tilde{\be}\in \Delta,
\eeq
where $\widetilde{B}_{\tilde{\al},\tilde{\beta}}(\la,\mu)$ is the phase factor from the composition of
the following two vertex operators:
\ben
\tGamma^{\tilde{\al}}(\la)\tGamma^{\tilde{\be}}(\mu) = \widetilde{B}_{{\tilde{\al}},{\tilde{\be}}}(\la,\mu)\,
:\tGamma^{\tilde{\al}}(\la)\tGamma^{\tilde{\be}}(\mu):\, .
\een
After a straightforward computation as in Section \ref{kac-peterson}, we get
\beq\label{phase-fac-calib}
\widetilde{B}_{{\tilde{\al}},{\tilde{\be}}}(\la,\mu) = 
\mu^{-(\al_0|\be_0)} e^{C_0
  (\al_0|\be_0)-2\pi\sqrt{-1}(\omega_b|\al)\, (\rho_b|\be)}
\prod_{m=1}^\kappa \Big( 1-\eta^m(\mu/\la)^{1/\kappa}\Big)^{(\si_b^m(\al)|\be)}.
\eeq

We are interested in the following system of Hirota quadratic
equations: for every integer $n\in \Z$
\beq\label{eth-ade}
{\rm Res}_{\la=\infty}\ \left.\frac{d\la}{\la} \Big(
  \widetilde{\Omega}_{\Delta'}(\la)\, (\tau\otimes_\liea \tau)
  \Big)\right|_{q_0^{01}\otimes 1-1\otimes q_0^{01}=n\sqrt{\hbar}}=0
\eeq
where $\tau\in \C_\hbar[\![q_0,q_1+1,q_2\dots]\!].$ The operator
$\widetilde{\Omega}_{\Delta'}(\la)$ is multivalued near $\la=\infty$: the analytic
continuation around $\la=\infty$ corresponds to a monodromy
transformation of each cycles ${\tilde{\al}}\in \Delta'$ of the type
${\tilde{\al}}\mapsto \si_b({\tilde{\al}})+ n_{\tilde{\al}} \varphi,$ where $n_{\tilde{\al}}\in \Z$. Using
Proposition \ref{per-tor} we get that the analytic continuation
transforms $\widetilde{\Omega}_{\Delta'}(\la)$ by permuting the cycles ${\tilde{\al}}$ and
multiplying each vertex operator term by $e^{2\pi\sqrt{-1}n_{\tilde{\al}}
  (q_0^{01}\otimes 1-1\otimes q_0^{01})}$. Therefore the 1-form in
\eqref{eth-ade} is invariant with respect to the analytic continuation
near $\la=\infty$. Moreover, for the same
  reason the equations \eqref{eth-ade} are independent of the choice
  of a lift $\Delta'$ of $\oDelta$. 
\begin{remark}
The Hirota quadratic equations \eqref{eth-ade} are a straightforward
generalization of the construction of Givental and Milanov \cite{GM}
(see also \cite{FGM}, where the coefficients $\tilde{b}_{\tilde{\al}}$ were interpreted
in terms of the vertex operators)
of integrable hierarchies for simple singularities. 
\end{remark}
The following is the main result of this section.

\begin{theorem}\label{t2}
If $\tau$ is a solution to the Hirota quadratic equations
\eqref{eth-ade}, then $\cF_C(\tau)$ with $C=\kappa^\chi Q$ is a
tau-function of the $\si_b$-twisted Kac--Wakimoto hierarchy.
\end{theorem} 
\proof
We just have to find the
Laplace transform of the Hirota quadratic equations (\ref{HQE_KW}) of the
Kac--Wakimoto hierarchy. Let $\al\in \Delta^{(0)}$ and $\tilde{\al}\in \Delta$ be as in Section \ref{cal-per}. Using Lemma \ref{vop-lt} we get
\ben
\Big(a_\al(\zeta)\,  E_\al(\zeta)\otimes E_{-\al}(\zeta)\Big)\ \Big(\cF_C\otimes \cF_C \Big)
= \Big(\cF_C\otimes \cF_C \Big)\, \Big(b_{\tilde{\al}}(\la)\, \tGamma^{\tilde{\al}}(\la)\otimes_\liea \tGamma^{-{\tilde{\al}}}(\la) \Big),
\een
where the coefficient $b_{\tilde{\al}}$ is given by 
\ben
a_\al(\zeta)\, \zeta^{-2\kappa|\al_0|^2}e^{\frac{|\al_0|^2}{\chi}\, \log
  C}e^{-4\pi\sqrt{-1}(\omega_b|\al)(\rho_b|\al)} .
\een
Recalling formula \eqref{a-al} and $\la=\zeta^\kappa/\kappa$ we get 
\beq\label{coef:b1}
b_{\tilde{\al}}(\la)=B_{\al,\al}\, \la^{-|\al_0|^2}\, 
e^{|\al_0|^2\, C_0}\, 
e^{-2\pi\sqrt{-1}(\omega_b|\al)(\rho_b|\al)} .
\eeq
Using \eqref{coef:b2} and \eqref{phase-fac-calib}, it is not hard to verify that $b_{\tilde{\al}}(\la)=\tilde{b}_{\tilde{\al}}(\la)$. 

In other words, $\cF_C(\tau)$ is a solution to the Kac--Wakimoto
hierarchy if $\tau$ satisfies the following equations:
\ben
{\rm Res}_{\la=\infty}\ \frac{d\la}{\la} \Big(
(\cF_C\otimes \cF_C)  \widetilde{\Omega}_{\Delta'}(\la)\, (\tau\otimes_\liea \tau)
  \Big)=0.
\een
Comparing the coefficients in front of $e^{(n'+x)\omega}\otimes
e^{(n''+x)\omega}$ we get \eqref{eth-ade} with $n=n'-n''.$
\qed

\section{The main Theorem} \label{sec:eth-anc}

\subsection{Vertex operators}
The symplectic  loop space formalism in GW theory was
introduced by Givental \cite{G3}. We apply this natural framework to
describe and investigate further the Hirota quadratic equations
\eqref{eth-ade}.  In this section, we again adopt the notation that $\al,\be$ are in the affine root system $\Delta$.

Recall the series \eqref{falpha}. We are interested in the vertex
operators 
\beq\label{vop-al}
\Gamma^\alpha(t,\la) = :e^{\widehat{\f}^\alpha(t,\la)}:,\quad \al\in \Delta,
\eeq
and their {\em phase factors} $B_{\alpha,\beta}(t,\la,\mu)$ defined by
\ben
\Gamma^\alpha(t,\la)\Gamma^\beta(t,\mu) =
B_{\alpha,\beta}(t,\la,\mu):\Gamma^\alpha(t,\la)\Gamma^\beta(t,\mu):
\,\quad \al,\be\in\Delta,
\een
where $:\cdot:$ is the usual normal ordering -- move all differentiation
operators to the right of the multiplication operators. Note that 
\beq\label{phase-factor}
B_{\al,\be}(t,\la,\mu) := e^{
  \Omega(\f_\al(t,\la;z)_+,\f_\be(t,\mu;z)_-)}.
\eeq
The action of the vertex operators on the Fock space is not well
defined in general. We would like to recall the conjugation laws from
\cite{G1} and to make sense of the vertex operator action on the Fock
space. 
\subsubsection{Vertex operators at infinity}\label{vop-inf}
Let us fix $t\in M$ and expand the vertex operators $\Gamma^\alpha(t,\la)$ in a
neighborhood of $\la=\infty.$ By definition (see (\ref{eqn:calib_lim})) we have
$f_\alpha(t,\la;z)=S_t\widetilde{\f}_\alpha(\la;z)$. Using formula
\eqref{S:fock}, it is easy to prove that  
\beq\label{S:vop}
\tGamma^\al(\la)\, \widehat{S}_t^{-1} = e^{\frac{1}{2}\,
  W(\tf_\al(\la)_+,\tf_\al(\la)_+)}
\widehat{S}_t^{-1} \, \Gamma^\al(t,\la).
\eeq
In particular, using the formal $\la^{-1}$-adic
topology we get that the vertex operator $\Gamma^\al(t,\la)$ defines a linear map  
$\C_{\hbar}[\![\q]\!] \to K_{\hbar}[\![\q]\!] ,$
where $K$ is an appropriate field extension of the field
$\C(\!(\la^{-1})\!).$ 

Let us explain the relation between the phase factors. Recall formula \eqref{phase-fac-calib}, 
the RHS is interpreted as an element in
$\C(\!(\la^{-1/\kappa})\!)(\!(\mu^{-1/\kappa})\!)$ by taking the Laurent
series expansion with respect to $\la$ at $\la=\infty$.
\begin{proposition}\label{intertw}
The following formula holds:
\ben
B_{\al,\be}(t,\la,\mu)  =  \widetilde{B}_{\al,\be}(\mu,\la) 
e^{W_t(\tf_\al(\mu)_+,\tf_\be(\la)_+)}.
\een
\end{proposition}
\proof
Conjugating the identity $\tGamma^\al(\la)\tGamma^\be(\mu) =
\widetilde{B}_{\al,\be}(\la,\mu)\, :\tGamma^\al(\la)\tGamma^\be(\mu):$ by
$\widehat{S}_t$ and using formula \eqref{S:vop}
we get that
\ben
e^{\frac{1}{2}\Big(W_t(\tf_\al(\la)_+,\tf_\al(\la)_+) +
  W_t(\tf_\be(\mu)_+,\tf_\be(\mu)_+) \Big)}\,
B_{\al,\be}(t,\la,\mu) 
\een
coincides with 
\ben
e^{\frac{1}{2}W_t(\tf_\al(\la)_++\tf_\be(\mu)_+,\tf_\al(\la)_++\tf_\be(\mu)_+)
}
\widetilde{B}_{\al,\be}(\la,\mu).
\een
The quadratic form $W$ is symmetric, so comparing the above identities
yields the desired formula.
\qed

\subsubsection{Vertex operators at a critical value}\label{sec:vop-kdv}
Let us assume now that $\la$ is near one of the critical values
$u_j(t)$ and that $\be$ is a cycle vanishing over
$\la=u_j(t), 1\leq j\leq N+1$. According to Lemma \ref{vanishing_a1}  we have
$\f_\be(t,\la;z) = \Psi_tR_t(z)\f_{A_1}(u_j,\la;z)$.  
Using Lemma \ref{lrfock} it is easy to prove (see \cite{G1},
Section 7) that 
\beq\label{R:vop}
\Gamma^\be(t,\la)\, \widehat{\Psi}_t\widehat{R}_t = e^{\frac{1}{2}\,
  V_t(\f_\be(t,\la)_-, \f_\be(t,\la)_-)}
\widehat{\Psi}_t\widehat{R}_t \, \Gamma_{A_1}^{\pm}(u_j,\la),
\eeq
where 
$
\Gamma_{A_1}^{\pm}(u_j,\la) = :e^{\pm \widehat{\f}_{A_1}(u_j,\la)}:
$ 
is the vertex operator of the $A_1$-singularity, $V_t$ is the second
order differential operator defined in Lemma 
\ref{lrfock}, and
\ben
 V_t(\f_\be(t,\la)_-, \f_\be(t,\la)_-) = \sum_{\ell,m=0}^\infty (I_\be^{(-\ell)}(t,\la),V_{\ell m}\,I_\be^{(-m)}(t,\la)).
\een
In this
case, the action of the vertex operators is 
well-defined on the subspace spanned by the tame asymptotical
functions and it yields a linear map 
\ben
\Gamma^\be(t,\la): \C_{\hbar}[\![\q]\!]_{\rm tame} \to K_{\hbar}[\![\q]\!],
\een
where $K=\C(\!((\la-u_j)^{1/2})\!).$
Furthermore, the phase factor $B_{\al,\be}(t,\la,\mu)$ is well defined if
$\be$ is a vanishing cycle, since it can be interpreted as an element
in $\C(\!((\mu-u_j)^{1/2})\!)(\!((\la-u_j)^{1/2})\!).$ Finally,
similarly to Proposition \ref{intertw}, we have 
\beq\label{intertw-R}
B_{\be,\be}(t,\la,\mu) = B_{A_1}(u_j,\la,\mu)
e^{-V_t(\f_\be(t,\la)_-,\f_\be(t,\mu)_-)},
\eeq
where  $B_{A_1}(u_j,\la,\mu)$ is the phase factor of the product
$\Gamma_{A_1}^{\pm}(u_j,\la) \Gamma_{A_1}^{\pm}(u_j,\mu)$. A straightforward computation gives
\beq\label{phase-a1}
B_{A_1}(u_j,\la,\mu) =
\left(
\frac{\sqrt{\la-u_j}-\sqrt{\mu-u_j}}{\sqrt{\la-u_j}+\sqrt{\mu-u_j}} 
\right)^2,
\eeq
where the RHS should be expanded into a Laurent series with respect to
$\mu$ at $\mu=u_j$.

\subsection{From descendants to ancestors}\label{sec:desc-anc}
Following our construction of the HQEs from Section \ref{sec:hqe} we
would like to introduce an integrable hierarchy for the ancestor
potential $\mathcal{A}_t$. Let us introduce 
the Heisenberg fields $$\phi_\beta(t,\la)=\d_\la
\widehat{\f}^\beta(t,\la), \quad \be\in\Delta',$$ 
and the corresponding Casimir operator
\ben
\Omega_{\Delta'}(t,\la) & = & 
-\frac{\la^2}{2}\Big(\sum_{m=1}^N\, 
:(\phi_{\be_m}(t,\la)\otimes_\liea 1-1\otimes_\liea\phi_{\be_m}(t,\la)) (\phi^{\be_m}(t,\la)\otimes_\liea
1-1\otimes_\liea\phi^{\be_m}(t,\la)) :\, \Big) +\\
&&
+
\sum_{\beta\in \Delta'}
b_\beta(t,\la)\Gamma^\beta(t,\la)\otimes_\liea
\Gamma^{-\beta}(t,\la) 
-\frac{1}{2}\operatorname{tr}\Big(\frac{1}{4}+\theta\, \theta^{\rm T}\Big),
\een
where $\{\be_m\}$ and
$\{\be^m\}$ are chosen as $\{\tilde{\al}_m\}$ and $\{\tilde{\al}^m\}$ as in Section \ref{sec:hqe}, and the coefficients
$b_\beta(t,\la)$ are defined by
\beq\label{coef-be}
b_\beta(t,\la)^{-1} = \lim_{\mu\to \la} \ \Big(1-\frac{\mu}{\la}\Big)^2
\, B_{\beta,-\beta}(t,\la,\mu).
\eeq
Finally, we need also to discretize the HQEs corresponding to the
above Casimir operator so that we offset the problem of
multivaluedness. Note that, for the toroidal cycle $\varphi$ in Section \ref{sec:tor_cycle}, according to Proposition \ref{per-tor} the
vector $\f^\varphi(t,\la;z)$ has only negative powers of $z$, so the
quantization $\widehat{\f}^\varphi(t,\la)$ is a linear function in
$\q$. 
\begin{lemma}\label{discr}
Let $\varphi$ be the toroidal cycle. Then the equation
\beq\label{discr-cond}
\widehat{\f}^\varphi(t,\la)\otimes 1- 1\otimes
\widehat{\f}^\varphi(t,\la)=2\pi\sqrt{-1}\, n
\eeq
is equivalent to
\begin{align}\label{discr-anc1}
[S_t^{-1}\q(z)]_{0,01}\otimes 1-1\otimes
[S_t^{-1}\q(z)]_{0,01} & = n\sqrt{\hbar}  \\
\label{discr-anc2}
[S_t^{-1}\q(z)]_{\ell,01}\otimes 1-1\otimes
[S_t^{-1}\q(z)]_{\ell,01}  & = 0, \quad \forall \ell\geq 1,
\end{align}
where $[S_t^{-1}\q(z)]_{\ell,i}$ denotes the coefficient of $S_t^{-1}\q(z)$ in front of $\phi_i
z^\ell$. 
\end{lemma}
\proof
Note that 
\ben
\widetilde{\f}^\varphi(\la;z)=2\pi\sqrt{-1}\ \sum_{\ell=0}^\infty \,
\frac{\la^\ell}{\ell!}\, \phi_{02}\,(-z)^{-\ell-1}.
\een
The equations \eqref{discr-anc1}--\eqref{discr-anc2} can be written
equivalently as 
\ben
\Omega(\widetilde{\f}^\varphi(\la;z), S_t^{-1}\q(z)) = 2\pi\sqrt{-1}\, n\sqrt{\hbar}.
\een
It remains only to recall that $S_t$ is a symplectic transformation
and that 
\ben
\f^\varphi(t,\la;z)=S_t \widetilde{\f}^\varphi(\la;z). 
\een
\qed

We will be interested in the following HQEs: for every integer $n\in
\Z$
\beq\label{eth-anc}
\operatorname{Res}_{\la=\infty} \left.\frac{d\la}{\la}\, \Big(
\Omega_{\Delta'}(t,\la)\, (\tau\otimes \tau) \Big) \right|_{
\widehat{\f}^\varphi(t,\la)\otimes 1- 1\otimes
\widehat{\f}^\varphi(t,\la)=2\pi\sqrt{-1}\, n } = 0,
\eeq
where $\tau$ belongs to an appropriate Fock space and we have to
require also that the discretization is well defined. For our purposes
the HQEs \eqref{eth-anc} will be on the Fock space $\C_\hbar[\![
q_0+t,q_1+1,q_2,\dots]\!]$. On the other hand the operator
$\widehat{S}_t^{-1}$ gives rise to an isomorphism
\ben
\widehat{S}_t^{-1}: \C_\hbar[\![ q_0+t,q_1+\one,q_2,\dots]\!]\to \C_\hbar[\![ q_0,q_1+1,q_2,\dots]\!].
\een
which allows us to identify the HQEs \eqref{eth-ade} and \eqref{eth-anc}.
\begin{proposition}\label{anc-desc}
A function $\tau$ is a solution to the HQEs \eqref{eth-anc} iff
$\widehat{S}_t^{-1}\tau$ is a solution to the HQEs \eqref{eth-ade}.
\end{proposition}
\proof
Using Proposition \ref{intertw} we get that
\ben
\widetilde{\Omega}_{\Delta'}(\la)\, (\widehat{S}_t^{-1}\otimes \widehat{S}_t^{-1}) = 
(\widehat{S}_t^{-1}\otimes \widehat{S}_t^{-1})\, \Omega_{\Delta'}(t,\la).
\een
It remains only to notice that the discretization in both HQEs are compatible
with the action of $\widehat{S}_t$, which follows easily from Lemma
\ref{lsfock} and Lemma \ref{discr}.\qed

\subsection{The integrable hierarchy for $A_1$-singularity}
It was conjectured by Witten \cite{W1} and first proved by Kontsevich
\cite{Ko1} that the total descendant potential of a point is a
tau-function of the KdV hierarchy. The latter can be written in two
different ways: via the Kac-Wakimoto construction and as a reduction
of the KP hierarchy. We will need both realizations, so let us recall
them. 
\subsubsection{The Kac--Wakimoto construction of KdV}
The Casimir operator (cf. Section \ref{sec:desc-anc}) for the $A_1$-singularity $f(x)=x^2/2+u$ takes the form
\ben
\Omega_{A_1}(u,\la) & = &
-\frac{\la^2}{4}\, :\phi^{V\otimes V}_{\be}(u,\la) \phi^{V\otimes
  V}_{\be}(u,\la): + \\
&&
 +b_{\be}(u,\la) \Big(
\Gamma_{A_1}^\be(u,\la)\otimes \Gamma_{A_1}^{-\be}(u,\la) + 
\Gamma_{A_1}^{-\be}(u,\la)\otimes \Gamma_{A_1}^\be(u,\la) \Big)-\frac{1}{8},
\een
where the coefficient 
\ben
b_{\be}(u,\la) =\lim_{\mu\to \la}
\Big(1-\frac{\mu}{\la}\Big)^{-2}\, B_{\be,\be}(u,\mu,\la) = \frac{\la^2}{16(\la-u)^2}.
\een
We denoted by $V$ the Fock space $\C_\hbar[\![\q]\!]$, and 
\ben
\phi^{V\otimes V}_{\be}(u,\la) := \phi_{\be}(u,\la)\otimes 1 - 1\otimes \phi_{\be}(u,\la).
\een
Witten's conjecture (Kontsevich's theorem) can be stated as follows:
\beq\label{kdv-kw}
\operatorname{Res}_{\la=\infty} \Omega_{A_1}(0,\la)\, (\cD_{\rm
  pt}\otimes \cD_{\rm pt})=0. 
\eeq
To compare the above equation with the principal Kac-Wakimoto
hierarchy of type $A_1$, note that 
\ben
\Gamma_{A_1}^\be (u,\la) = \exp\Big(2\sum_{n=0}^\infty
\frac{(2(\la-u))^{n+1/2}}{(2n+1)!!}\,\frac{q_n}{\sqrt{\hbar}}\Big) 
\exp\Big(- 2\sum_{n=0}^\infty
\frac{(2n-1)!! }{(2(\la-u))^{n+1/2}}\sqrt{\hbar}\d_n\Big),
\een
and that the coefficient in front of $\la^{-2}$ in 
$ \frac{1}{4}:\phi^{V\otimes V}_{\be}(0,\la) \phi^{V\otimes
V}_{\be}(0,\la):$ is precisely
\ben
\sum_{n=0}^\infty \Big(n+\frac{1}{2}\Big)(q_n\otimes 1 -1\otimes
q_n)(\d_n\otimes 1-1\otimes \d_n),
\een
where $\d_n:=\d/\d q_n.$ It follows that the above equations coincide
with the Kac--Wakimoto form of the KdV hierarchy up to the rescaling
$q_n=t_{2n+1} (2n+1)!!$.

On the other hand, the total descendant potential $\cD_{\rm pt}$
satisfies the string equation, which can be stated as follows (see
\cite{G3}): $e^{(u/z)\sphat} \cD_{\rm pt}= \cD_{\rm pt}$. Using that 
\ben
\Omega_{A_1}(0,\la)\, \Big(e^{(u/z)\sphat} \otimes
e^{(u/z)\sphat}\Big)  =  
\Big(e^{(u/z)\sphat} \otimes
e^{(u/z)\sphat}\Big) \, \Omega_{A_1}(u,\la)
\een
we get that $\cD_{\rm pt}$ satisfies also the following HQEs:
\beq\label{kdv-1}
\operatorname{Res}_{\la=\infty} \Omega_{A_1}(u,\la)\, (\cD_{\rm
  pt}\otimes \cD_{\rm pt})  =0. 
\eeq
\subsubsection{The KdV hierarchy as a reduction of KP}
According to Givental \cite{G1} the KdV hierarchy \eqref{kdv-kw} can
be written also as
\ben
\operatorname{Res}_{\la=0}\, \Big(\sum_{\pm} \frac{d\la}{\pm\sqrt{\la}}\,
\Gamma_{A_1}^{\pm\be/2}(0,\la)\otimes \Gamma_{A_1}^{\mp\be/2}(0,\la)\Big)\, 
(\cD_{\rm pt}\otimes \cD_{\rm pt}) = 0.
\een
Using again the string equation and Proposition \ref{intertw} we get that
$\cD_{\rm pt}$ satisfies also the following HQEs:
\beq\label{kdv-2}
\operatorname{Res}_{\la=u}\, \Big(\sum_{\pm} \frac{d\la}{\pm\sqrt{\la-u}}\,
\Gamma_{A_1}^{\pm\be/2}(u,\la)\otimes \Gamma_{A_1}^{\mp\be/2}(u,\la)\Big)\, 
(\cD_{\rm pt}\otimes \cD_{\rm pt}) = 0.
\eeq

\subsection{The phase factors}\label{sec:phase-factors}

In this section we will prove Proposition \ref{phase-mon}, that the phase factors $B_{\al,\be}(t,\la,\mu)$ (see \eqref{phase-factor}) are multivalued
analytic function and that the analytic continuation is compatible with the
monodromy action on the cycles $\alpha$ and $\beta$. 
To begin with put
\ben
B^\infty_{\al,\be}(t,\la,\mu)=\exp\Omega^\infty_{\al,\be}(t,\la,\mu),
\een
where
\beq\label{phase-series}
\Omega^\infty_{\al,\be}(t,\la,\mu):=\iota_{\la^{-1}}\iota_{\mu^{-1}}
\sum_{n=0}^\infty
(-1)^{n+1}(I^{(n)}_\al(t,\la),I^{(-n-1)}_\be(t,\mu)),
\eeq
where $\iota_{\la^{-1}}$ (resp. $\iota_{\mu^{-1}}$) is the Laurent
series expansion at $\la=\infty$ (resp. $\mu=\infty$). 
The differential of \eqref{phase-series} with respect to $t$ is 
\ben
\widetilde{\cW}_{\al,\be}(\la,\mu) := 
I^{(0)}_{\al}(t,\la)\bullet I^{(0)}_\be(t,\mu) = \sum_{i\in \lieI} 
(I^{(0)}_{\al}(t,\la),\d_i\bullet I^{(0)}_\be(t,\mu))\,dt_i,
\een
which will be interpreted as a 1-form on $M$ depending on the parameters $\la$ and $\mu$. 
Furthermore, for each $t\in M$, put $r(t)=\operatorname{max}_j |u_j(t)|$, where $\{u_j(t)\}_{j=1}^{N+1}$ is the set of all critical values of $F(x,t).$ 
In other words, $r(t)$ is the radius of the smallest disk (with center at $0$) that contains all critical values of $F(x,t)$. 
Let
\ben
D^+_\infty=\{(t,\la,\mu)\in M\times \C^2\ :\ |\la-\mu|<|\mu|-r(t)<|\la|-r(t) \}.
\een
Note that since $|\la-\mu|\geq 0$ we have $|\la|>r(t)$ and $|\mu|>r(t)$ for all $(t,\la,\mu)\in D_\infty^+$, which implies that the 
Laurent series expansions of $I^{(0)}_\al(t,\la)$ and $I^{(0)}_\be(t,\mu)$ at respectively $\la=\infty$ and $\mu=\infty$ are convergent.
The first inequality in the definition of $D_\infty^+$ guarantees that the line segment $[\la,\mu]$ is outside the disk $|x|\leq r(t)$. In
particular, in order to specify a branch of $\cW_{\al,\be}(\la,\mu)$ it is enough to specify the branches of the period vectors only at the 
point $(t,\la)$, the branch of the periods at $(t,\mu)$ is determined via the line segment $[\la,\mu]$.  
\begin{proposition}\label{conv}
The series \eqref{phase-series} is convergent for all $(t,\la,\mu)\in
D^+_\infty$. 
\end{proposition}
\proof
Using Proposition \ref{intertw} we can write \eqref{phase-series} as a
sum of two formal series 
\beq\label{log-phase}
\Omega^\infty_{\al,\be}(t,\la,\mu) =
\widetilde{\Omega}^\infty_{\al,\be}(\la,\mu)+W_t(\widetilde{\f}_\al(\la)_+,\widetilde{\f}_\be(\mu)_+ ),
\eeq
where $\widetilde{\Omega}^\infty_{\al,\be}$ is the Laurent series expansion of $\log \widetilde{B}_{\al,\be} $ in the domain $|\la|>|\mu|$. Since the series 
 $\widetilde{\Omega}^\infty_{\al,\be}$ is convergent for $|\la|>|\mu|>|\la-\mu|$, it is
enough to prove the proposition for the second series on the RHS of \eqref{log-phase}. Recalling the
definition of $W_t$ and using the fact that modulo $Q$ the series
$S_t(z)=e^{\frac{1}{z} t\,\cup}$, where $t\,\cup$ means the classical orbifold cup
product multiplication by $t$, we get that 
\ben
\lim_{\operatorname{Re}(t_{02})\to -\infty}\, \lim_{t\to (0,\cdots,0,t_{02})}\ (W_t-
t_{02}\,P) = 0,
\een
On the other hand, since
\ben
d W_t(\widetilde{\f}_\al(\la)_+,\widetilde{\f}_\be(\mu)_+ ) = 
d \Omega_{\al,\be}(t,\la,\mu)  = I^{(0)}_\al(t,\la)\bullet I^{(0)}_\be(t,\mu),
\een
the series 
\beq\label{W-phase}
 W_t(\widetilde{\f}_\al(\la)_+,\widetilde{\f}_\be(\mu)_+ )-t_{02}(\al_0|\be_0)/\chi,
\eeq
viewed as a formal Laurent series in $\la^{-1}$ and $\mu^{-1}$ can be identified with the improper integral
\beq\label{phase:int-inf}
\lim_{\ep \to \infty}\int_\ep^t \Big( I^{(0)}_\al(t',\la)\bullet I^{(0)}_\be(t',\mu) - dt'_{02}(\al_0|\be_0)/\chi\Big),
\eeq
where $\ep\in M$ and the limit is taken along a straight segment,
s.t., $\ep_i\to 0$ for $i\neq 02$ and $\operatorname{Re}(\ep_{02})\to -\infty$. More precisely, if we take the Laurent series 
expansion of the integrand at $\la=\infty$ and $\mu=\infty$ and integrate termwise, we get 
\eqref{W-phase}. It remains only to notice that the integrand extends holomorphically at the limiting point $\ep=\infty$ (because 
we removed the singularity), so the termwise integration preserves the convergence.
\qed

The proof of Proposition \ref{conv} yields slightly more. Namely, we proved that the 2nd summand on the RHS of \eqref{log-phase} 
is a convergent Laurent series in $\la^{-1}$ and $\mu^{-1}$ and that the corresponding limit is a multi-valued analytic function on 
\ben
D_\infty:=\{(t,\la,\mu)\in M\times \C^2 \ :\  |\la-\mu|<\operatorname{min}(|\la|-r(t),|\mu|-r(t))\}.
\een
On the other hand, the phase factor $\widetilde{B}_{\al,\be}(\la,\mu)$ is also a multivalued
analytic function on $D_\infty$ except for a possible pole along $\la=\mu$. Hence we have the following corollary (of the proof).
\begin{corollary}\label{ac-1}
The series $B^\infty_{\al,\be}(t,\la,\mu)$ extends analytically to a multivalued analytic function on $D_\infty$ except for a possible pole along 
the diagonal $\la=\mu$.
\end{corollary}
Using the analytic extension of $B^\infty_{\al,\be}(t,\la,\mu)$ we define a multi-valued function with values in the 
space $\C\{\!\{\xi\}\!\} $ of convergent Laurent series at $\xi=0$ in the following way:
\ben
B_{\al,\be}:(M\times \C)_\infty \to \C\{\!\{\xi\}\!\},\quad (t,\la)\mapsto \iota_{\mu-\la} B^\infty_{\al,\be}(t,\la,\mu),
\een
where $\xi=\mu-\la$, $\iota_{\mu-\la}$ is the Laurent series expansion at $\mu=\la$, and 
\ben
(M\times \C)_\infty:=\{(t,\la)\in M\times \C\ :\ |\la|>r(t)\}.
\een 
It is convenient to introduce the 1-form $\cW_{\al,\be}(\xi) := \widetilde{\cW}_{\al,\be}(0,\xi)$. Following \cite{G1} we call $\cW_{\al,\be}(\xi)$ 
the {\em  phase form}. Note that if $C\subset (M\times\C)_\infty$ is a path from $(t,\la)$ to $(t',\la')$, then 
\beq \label{ac-2}
B_{\al,\be}(t',\la') = B_{\al,\be}(t,\la)\, e^{\int_C \cW_{\al,\be}(\xi)}.
\eeq
Therefore we can uniquely extend the function $B_{\al,\be}$ to a function on $(M\times \C)'$, so that
formula \eqref{ac-2} holds for all paths $C\subset (M\times \C)'$. Finally, for every $(t,\la)\in (M\times\C)'$ and $\mu$ sufficiently close to $\la$
we define 
\ben
B_{\al,\be}(t,\la,\mu) = \left. B_{\al,\be}(t,\la)\right|_{\xi=\mu-\la},\quad \Omega_{\alpha,\beta}(t,\la,\mu):=\log B_{\al,\be}(t,\la,\mu).
\een 
Note that $B_{\al,\be}(t,\la,\mu) =B_{\al,\be}^\infty(t,\la,\mu)$ if $(t,\la,\mu)\in D^+_\infty$.

Let $t_0\in M$ be a generic point, so that all critical points of
$F(x,t_0)$ are of type $A_1$ and the absolute values of the corresponding critical values are
pairwise distinct. Let $u_j(t_0)$ be a critical value of $F(x,t_0)$ 
with a maximal absolute value, i.e., $|u_j(t_0)|=r(t_0)$. There exists a real number $\epsilon_0>0$, s.t., if $|x|<\epsilon_0$,
then $r(t_0+x\one)=|u_j(t_0)+x|.$  We fix 
$t=t_0+x_0\one,$ $\la$, and $\mu$, s.t., the line segment $[\mu-u_j(t_0),x_0]$ is contained inside the disk $\{ |x|<\epsilon_0\}\subset \C$  
and the line segment $[t_0,t]\times \{(\la,\mu)\} \subset D^+_\infty$ . For example, fix $\mu$, s.t., $|\mu|>u_j(t_0)$ and 
$|\mu-u_j(t_0)|<\ep_0$ and put $x_0=\frac{1}{2}(\mu-u_j(t_0))$, then we can find $\la$ such that all requirements are fulfilled. 
\begin{lemma}\label{int:phase-inf}
If $\be\in  H_2(X_{t,\mu};\Z)$ is a cycle vanishing over $t=t_0+(\mu-u_j(t_0))\one$, then  
\beq\label{integral}
\Omega_{\al,\be}(t,\la,\mu) = \lim_{\ep\to 0} \int^{t}_{t_0+(\ep+\mu-u_j(t_0))\one} \widetilde{\mathcal{W}}_{\al,\be}(\la,\mu),
\eeq
where the integration is along a straight segment.
\end{lemma}
\proof
By definition $\Omega_{\al,\be}(t,\la,\mu)$ is the  Laurent
series expansion near $\la=\infty$ of the series 
\beq\label{series}
\sum_{n=0}^\infty (-1)^{n+1}(I_\al^{(n)}(t,\la),I^{(-n-1)}_\be(t,\mu)),
\eeq
while the RHS of \eqref{integral} is 
\beq\label{integral-2}
\lim_{\ep\to 0}\int^{x_0}_{\ep+\mu-u_j}  (I^{(0)}_\al(t_0,\la-x),I^{(0)}_\be(t_0,\mu-x))dx.
\eeq
Using integration by parts $(n+1)$ times and the fact that the periods $I^{(-p-1)}_\be(t_0,\mu-x)$ vanish at $x=\mu-u_j$,
we get that the integral \eqref{integral-2} coincides with
\beqa\notag 
\sum_{p=0}^n (-1)^{p+1}(I^{(p)}_\al(t_0',\la), I^{(-p-1)}_\be(t_0',\mu)) + \\ 
\label{integral-3}
\lim_{\ep\to 0} \, (-1)^{n+1}
\int^{x_0}_{\ep+\mu-u_j}  (I^{(n+1)}_\al(t_0,\la-x),I^{(-n-1)}_\be(t_0,\mu-x))dx.
\eeqa
The Laurent series expansion of $I^{(n+1)}_\al(t_0,\la-x) = I^{(n+1)}_\al(t_0+x\one,\la)$ in $\la^{-1}$ has radius of convergence 
$r(t_0+x\one)$. Hence, it is uniformly convergent for all 
$x$ that vary along a compact subset of the open subset in $\C$ defined by the inequality
\ben
\{x\in \C\ :\ |\la|>r(t_0+x\one)\}.
\een
On the other hand, according to our choice of $x_0,\la,$ and $\mu$, the point $(t_0+x\one,\la,\mu)\in D_\infty^+$ for all $x$ on the integration path. In particular, $|\la|>|\mu|>r(t_0+x\one)$, which means that the integration
path is entirely contained in the above open subset. Hence the integral \eqref{integral-3} has a convergent Laurent
series in $\la^{-1}$. Moreover, the leading order term of the expansion is $\la^{-e}$ for some rational number $e>n.$ 
This proves that the Laurent series expansions (in $\la^{-1}$) of the integral \eqref{integral-2} and of the series \eqref{series} coincide.
\qed

\bigskip
Our next goal is to prove that the analytic continuation of the phase factor $B_{\al,\be}(t,\la,\mu)$ 
is compatible with the monodromy representation in the following sense. Recall the monodromy representation 
(cf. Section \ref{sec:periods})
\ben
\rho: \pi_1((M\times \C)')\rightarrow \operatorname{GL}(\lieh).
\een  
Let $U\subset (M\times \C)'$ be an open subdomain and $f_{\al,\be}(t,\la)$ be
a (vector-valued) function depending bi-linearly on $(\al,\be)\in \lieh\times\lieh$
and analytic in a neighborhood of some point $(t_0,\la_0)\in
U$. We say that $f_{\al,\be}$ is {\em multi-valued analytic} on $U$ if it can be
extended analytically along any path in $U$. Furthermore, we say that 
$f_{\al,\be}$ is {\em compatible} with the monodromy representation $\rho$,
if for every closed loop $C$ in $U$, the analytic continuation of
$f_{\al,\be}(t,\la)$ along $C$ coincides with
$f_{w(\al),w(\be)}(t,\la)$, where $w=\rho(C)$ is the corresponding
monodromy transformation. 

Recall that (see Corollary \ref{ac-1}) the Laurent series $\Omega^\infty_{\be,\al}(t,\la,\mu) $ extends analytically to a multi-valued
analytic function $\Omega_{\be,\al}(t,\la,\mu) $ defined for all $(t,\la,\mu)\in D_\infty$, s.t., $\la\neq \mu$.
\begin{lemma}\label{locality}
Let $\al$ and $\be$ be cycles in the vanishing cohomology, s.t., $(\al|\be)=0$ then 
\ben
\Omega_{\al,\be}(t,\la,\mu)-\Omega_{\be,\al}(t,\mu,\la) = 2\pi\sqrt{-1}\operatorname{SF}(\al,\be)\quad 
\forall (t,\la,\mu)\in D^+_\infty,
\een
where $\operatorname{SF}$ is the bi-linear form \eqref{seifert-form}.
\end{lemma}
\proof
Since the difference
\ben
\Omega_{\be,\al}(t,\la,\mu) - \widetilde{\Omega}_{\be,\al}(\la,\mu),\quad  
\mbox{where}\quad 
\widetilde{\Omega}_{\be,\al}(\la,\mu) := \log \widetilde{B}_{\be,\al}(\la,\mu),
\een
has a convergent Laurent series expansion in $D_\infty$ and it is invariant under switching 
$(\be,\la)\leftrightarrow (\al,\mu)$, it is enough to prove the statement for 
$\widetilde{\Omega}_{\al,\be}(\la,\mu)$ where $(\la,\mu)$ is a point in the open subset 
\ben
\{|\la-\mu|<\operatorname{min}(|\la|,|\mu|)\}\subset \C^2.
\een
Recalling formula \eqref{phase-fac-calib}, the rest of the proof is a straightforward computation (see also the proof
of Lemma \ref{const-c}, where some of the computations were already done). 
\qed
\begin{remark}
If we omit the condition $(\al|\be)=0$ in Lemma \ref{locality}, then the identity is true only up to an integer 
multiple of $2\pi\sqrt{-1} (\al|\be)$. The ambiguity comes from the fact that the phase factor 
$\widetilde{\Omega}_{\al,\be}(\la,\mu)$ has a logarithmic singularity along $\la=\mu$ of the type 
$(\al|\be)\log(\la-\mu)$.
\end{remark}

\begin{proposition}\label{mon-comp}
The phase factor $B_{\al,\be}(t,\la)$ is compatible with the
monodromy representation in the domain $(M\times \C)'$.
\end{proposition}
\proof
By definition we have to prove that if $C'\subset (M\times \C)'$ is an arbitrary loop based at 
$(t,\la)$ and $\mu$ is sufficiently close to $\la$, then
\ben
B_{w(\al),w(\be)}(t,\la,\mu) = B_{\al,\be}(t,\la,\mu)\, e^{\int_C \widetilde{\cW}_{\al,\be}(\la,\mu)},
\een
where $w=\rho(C')$ and $C\subset M$ is the path parametrized by 
\ben
t'+(\la-\la')\one,\quad (t',\la')\in C'.
\een
We may assume that $(t,\la,\mu)\in D^+_\infty$, because by definition the value of $B_{\al,\be}$ at any other point differs by
an integral along the path of the phase form $\widetilde{\cW}_{\al,\be}(\la,\mu).$ Under this assumption the above equality is equivalent to
\beq\label{mon-inv}
\Omega^\infty_{w(\al),w(\be)}(t,\la,\mu) = \Omega^\infty_{\al,\be}(t,\la,\mu)+\int_C \widetilde{\cW}_{\al,\be}(\la,\mu)\quad 
(\operatorname{mod}\ 2\pi\sqrt{-1}\, \Z).
\eeq
We first prove  a special case of the above formula. Namely, let us choose a generic point $t_0\in M$, s.t.,
the absolute values of the critical values of $F(x,t_0)$ are pairwise distinct and let  $u_j(t_0)$ be the critical value with maximal absolute value
(here the notation is the same as in Lemma \ref{int:phase-inf}). We will assume that $t=t_0+x_0\one$ is sufficiently close to 
$t_0+(\mu-u_j(t_0))\one$ and that $C$ is a closed loop of the type $t_0+x\one$, where the parameter $x$ varies
along a small closed loop based at $x_0\in \C$ going around $\mu-u_j(t_0)$, so that the line segment $[\la-x,\mu-x]$ 
moves around $u_j$. Let us denote by $\gamma\in H_2(X_{t,\la};\Z)$ the vanishing cycle
vanishing over $(t_0,u_j(t_0))$, then we have the following
decompositions:
\ben
\al=\al'+\frac{(\al|\gamma)}{2}\, \gamma,\quad
\be=\be'+\frac{(\be|\gamma)}{2}\, \gamma,
\een
where $\al'$ and $\be'$ are cycles invariant w.r.t. the local monodromy around the point $(t_0,u_j(t_0))$. 
After a straightforward computation we get
\ben
\Omega_{w(\al),w(\be)}(t,\la,\mu) - \Omega_{\al,\be}(t,\la,\mu) = 
-(\al|\gamma)\Omega_{\gamma,\be'}(t,\la,\mu) - (\be|\gamma)\Omega_{\al',\gamma}(t,\la,\mu),
\een
while
$
\int_C \widetilde{\cW}_{\al,\be}(\la,\mu) 
$
is
\beq\label{integral-C}
\frac{1}{2}(\be|\gamma)\int_C \widetilde{\cW}_{\al',\gamma}(\la,\mu) + 
\frac{1}{2}(\al|\gamma)\int_C \widetilde{\cW}_{\gamma,\be'}(\la,\mu) +
\frac{1}{4}\, (\al|\gamma)\, (\be|\gamma)\, \int_C \widetilde{\cW}_{\gamma,\gamma}(\la,\mu),
\eeq
where we used that $\int_C \widetilde{\cW}_{\al',\be'}(\la,\mu)=0,$ because the
periods $I^{(0)}_{\al'}(t_0,\la-x)$ and $I^{(0)}_{\be'}(t_0,\mu-x)$ are holomorphic respectively at $x=\la-u_j$ and $x=\mu-u_j,$
which means that the phase form is holomorphic inside the loop $C$. 
The last integral in the above formula is easy to compute because only the singular terms of  $I^{(0)}_\gamma(t_0,\la-x)$ and $I^{(0)}_\gamma(t_0,\mu-x)$ contribute, i.e., 
\ben
\int_C \widetilde{\cW}_{\gamma,\gamma}(\la,\mu)  = 2\oint\frac{dx}{\sqrt{(\la-u_j(t_0)-x)(\mu-u_j(t_0)-x)} }=
  4\pi\sqrt{-1}.  
\een
According to Lemma \ref{int:phase-inf} 
\ben
\Omega_{\al',\gamma}(t,\la,\mu) = \int_{t_0+(\mu-u_j(t_0))\one}^t \widetilde{\cW}_{\al',\gamma}(\la,\mu)
\een 
and the integral on the RHS has a convergent Laurent series expansion in $\la-u_j(t)$ and $(\mu-u_j(t))^{1/2}$, which allows
us to evaluate the integral
\ben
\int_C \widetilde{\cW}_{\al',\gamma}(\la,\mu)  = -2 \int_{t_0+(\mu-u_j(t_0))\one}^t \widetilde{\cW}_{\al',\gamma}(\la,\mu) =
-2 \Omega^\infty_{\al',\gamma}(t,\la,\mu) =-2 \Omega_{\al',\gamma}(t,\la,\mu) .
\een
It remains only to evaluate the 2nd integral in \eqref{integral-C}. We have
\ben
\int_C \widetilde{\cW}_{\gamma,\be'}(\la,\mu) =  \int_C \widetilde{\cW}_{\be',\gamma}(\mu,\la) = 
-2\Omega_{\be',\gamma}(t,\mu,\la),  
\een
where the 2nd identity is derived just like above when $|\mu|>|\la|$,
and then we use analytic continuation to extend the formula for
$|\mu|<|\la|$ as well. Recalling Lemma \ref{locality}, we get
\ben
\Omega_{\be',\gamma}(t,\mu,\la) =\Omega_{\gamma,\be'}(t,\la,\mu) + 2\pi\sqrt{-1}
\operatorname{SF}(\be',\gamma).
\een
Using that $\be'=\be-(\be|\gamma)\gamma/2$ and that 
$\operatorname{SF}(\gamma,\gamma)=1$, we finally get
\ben
\int_C \widetilde{\cW}_{\gamma,\be'}(\la,\mu) =
-2\Omega_{\gamma,\be'}(t,\la,\mu) - 4\pi\sqrt{-1}
\operatorname{SF}(\be,\gamma) + 2\pi\sqrt{-1} (\be|\gamma).
\een
Since $\operatorname{SF}(\be,\gamma) \in \Z$, the proof of formula \eqref{mon-inv} in the special case is complete. 

The general case follows easily, because the fundamental group  $\pi_1((M\times \C)')$ is generated by loops like the above one. 
Indeed, we already know that the affine cusp polynomial $f(x)$ has a real Morsification $F(x,t'_0)$, i.e., all critical points of 
$F(x,t'_0)$ are real and the corresponding critical values are real as well. In particular, we can find a small deformation 
$F(x,t_0)$ of the real Morsification, s.t., the critical values $u_j$ are vertices of a convex polygon.  
The fundamental group $\pi_1((M\times \C)')$ 
is generated by simple loops in $\{t_0\}\times \C$ that go around the vertices of the polygon. Let us pick one of these loops and let $(t_0,u_j(t_0))$ be the corresponding vertex of the polygon. Since the translations of the 
type  $t_0\mapsto  t_0+c\one$, $c\in \C$, do not change the homotopy class of the loop, we can find a representative 
(namely, pick $c$, s.t., the $|u_j(t_0)+c|>|u_j(t_0)+c|$ for all other vertices $(t_0,u_j(t_0))$) of the homotopy class, which has the 
special form from above. 
\qed

\begin{proposition}\label{phase-mon}
There exists a generic point $t_0\in M$ (i.e. $F(x,t_0)$ is a Morse function) and a critical value $u_j(t_0)$, s.t., 
\beq\label{int-phase}
B_{\al,\be}(t,\la,\mu) = \lim_{\ep\to 0}\exp\Big( 
-\int_t^{t_0+(\ep+\mu-u_j(t_0))\one} \widetilde{\cW}_{\al,\be}(\la,\mu) 
\Big),
\eeq
where the integration is along any path avoiding the poles of the 1-form $\widetilde{\cW}_{\al,\be}(\la,\mu)$, s.t., the cycle 
$\be\in H_2(X_{t,\mu},\Z)$ vanishes along it. 
\end{proposition}
\proof
Let us assume that $t_0$ is a generic point and that $u_j(t_0)$ is the critical value with maximal absolute value. 
It is enough to prove the statement for an arbitrary point $(t,\la,\mu)\in D^+_\infty$, because by definition
the value of $B_{\al,\be}(t',\la,\mu)$ at any other point $(t',\la,\mu)$ differs by an integral of 
$\widetilde{\cW}_{\al,\be}(\la,\mu) $ along a path connecting $t$ and $t'$, while the RHS of \eqref{int-phase} clearly has the
same property. Let $(t,\la,\mu)\in D_\infty^+$ be a point such that Lemma \ref{int:phase-inf} holds and let $C_\ep''$ be the straight segment 
$[t, t_0+(\ep+\mu-u_j(t_0))\one]$. Put $C'=(C_\ep'')^{-1}\circ C_\ep$ and $w=\rho(C')$, where $C_\ep$ is the integration path 
(from $t$ to $t_0+(\ep+\mu-u_j(t_0))\one $), then by definition the cycle $w(\be)\in H_2(X_{t,\mu};\Z)$ is the vanishing cycle along 
the line segment $[t,t_0+(\mu-u_j(t_0))\one]$. According to Lemma \ref{int:phase-inf},
formula \eqref{int-phase} holds for $C''$ and $B_{w(\al),w(\be)}$. Therefore, we need to prove that   
\beq\label{gp-coh}
-\int_{C'} \widetilde{\mathcal{W}}_{\al,\be}(\la,\mu) = \Omega_{\al,\be}(t,\la,\mu)-\Omega_{w(\al),w(\be)}(t,\la,\mu) 
\quad (\operatorname{mod}\ 
2\pi\sqrt{-1} \Z),
\eeq
which follows from Proposition \ref{mon-comp}.
\qed

\subsection{The ancestor solution}

Now we are in a position to prove 
\begin{theorem}\label{thm:ancestor}
The total ancestor potential $\cA_t(\hbar;\q)$ is a solution to the HQEs \eqref{eth-anc}.
\end{theorem}

To begin with, put $\q'=\q\otimes 1$, $\q''=1\otimes \q$, and let us
assume that the discretization condition \eqref{discr-cond} is
satisfied for some integer $n$. The tameness of $\cA(\hbar;\q)$
implies that the LHS of \eqref{eth-anc} (for $\tau=\cA(\hbar;\q)$) is a formal
series in $\q'$ and $\q''$ with coefficients
formal Laurent series in $\sqrt{\hbar}$, whose coefficients are
polynomial expressions of the period vectors $I^{(n)}_\al(t,\la)$. In
particular, the residue in \eqref{eth-anc} can be computed via the
residue theorem, i.e., we have to compute the residues at all critical
points and at $\la=0$ and prove that their sum is 0.

Let $u_j(t)$ be one of the critical points of $F$, where $t\in M$ is a
generic point such that all critical values are pairwise
different. Furthermore, we assume that $\la$ is near $u_j(t)$ and that
a path in $(M\times \C)'$ from the reference point $(0,1)$ to
$(t,\la)$ is fixed in such a way that the vanishing cycle $\be$,
vanishing over $\la=u_j(t)$, belongs to the subset $\Delta'$ of affine roots defined in Section \ref{sec:hqe}. 

\subsubsection{The Virasoro term}  \label{vir-term}
Let us compute 
\beq\label{res-i}
-\operatorname{Res}_{\la=u_j(t)} \frac{\la}{2} d\la\, 
\sum_{m=1}^N:\phi^{V\otimes V}_{\be_m}(t,\la) \phi^{V\otimes
  V}_{\be^m}(t,\la): \ \cA_t^{\otimes 2}, 
\eeq
where $\phi^{V\otimes V}_\al:=\phi_\al\otimes 1-1\otimes
\phi_\al$. Put $\be_m=\al_m+(\be_m|\be)\be/2$ and
$\be^m=\al^m+(\be^m|\be)\be/2$, where $(\al_m|\be)=(\al^m|\be)=0.$
The above operator can be written as the sum of
\ben
\sum_{m=1}^N:\phi^{V\otimes V}_{\al_m}(t,\la) \phi^{V\otimes
  V}_{\al^m}(t,\la): + \Big(\sum_{m=1}^N (\be_m|\be)(\be^m|\be)\Big)\frac{1}{4}
:\phi^{V\otimes V}_{\be}(t,\la) \phi^{V\otimes V}_{\be}(t,\la):
\een
and
\beq\label{mixed}
\sum_{m=1}^N \frac{1}{2}\Big(
(\be_m|\be)
:\phi^{V\otimes V}_{\be}(t,\la) \phi^{V\otimes V}_{\al^m}(t,\la): +
(\be^m|\be)
:\phi^{V\otimes V}_{\be}(t,\la) \phi^{V\otimes V}_{\al_m}(t,\la):
\Big)
\eeq
 The Picard-Lefschetz formula implies that the periods
$I^{(n)}_{\al_m}(t,\la)$ and $I^{(n)}_{\al^m}(t,\la)$ are invariant
with respect to the local monodromy around $\la=u_j(t)$, so they must be
holomorphic in a neighborhood of $\la=u_j(t)$. 
The operator $\phi_\varphi^{V\otimes V}(t,\la)$, where $\varphi$ is
the toroidal cycle,  vanishes after we impose the discretization
condition \eqref{discr-cond}. On the other hand, since
$\sum_m(\be_m|\be)(\be^m|\al)=(\be|\al)$, the cycles 
\ben
-\be+ \sum_{m=1}^N (\be_m|\be)\be^m\quad \mbox{and}\quad
-\be+ \sum_{m=1}^N (\be^m|\be)\be_m
\een 
are in the kernel of the intersection form, so they must be
proportional to $\varphi$. Hence the operator \eqref{mixed} vanishes
after the discretization condition \eqref{discr-cond} is imposed. The residue
\eqref{res-i} turns into
\ben
-\operatorname{Res}_{\la=u_j(t)} \frac{\la}{4} d\la\,
:\phi^{V\otimes V}_{\be}(t,\la) \phi^{V\otimes V}_{\be}(t,\la):\ \cA_t(\hbar;\q')\cA_t(\hbar;\q'').
\een
To compute the above residue, note that the expression
\ben
:\phi^{V\otimes V}_{\be}(t,\la) \phi^{V\otimes V}_{\be}(t,\la):
(\widehat{\Psi}_t\widehat{R}_t)^{\otimes 2}  
\een
can be written as
\ben
(\widehat{\Psi}_t\widehat{R}_t)^{\otimes 2} 
:\phi^{V\otimes V}_{A_1}(u_j,\la) \phi^{V\otimes V}_{A_1}(u_j,\la): +
2V_t(\phi_\be(t,\la)_-, \phi_\be(t,\la)_-).
\een
Let us compute 
\ben
-\operatorname{Res}_{\la=u_j(t)} \frac{\la}{4}
d\la\,2V_t(\phi_\be(t,\la)_-, \phi_\be(t,\la)_-) = 
-\operatorname{Res}_{\la=u_j(t)} \frac{\la}{2}
d\la\,(V_{00}(t) I^{(0)}_\be(t,\la),I^{(0)}_\be(t,\la)),
\een
where we used the fact that only the leading term (w.r.t. $z$) of
$\phi_\be(t,\la;z)_-=-I^{(0)}_\be(t,\la)z^{-1}+\cdots$ will contribute
because the remaining ones have a zero at $\la=u_j(t)$ of order at least $\frac{1}{2}$. Furthermore, the Laurent series expansion of
$I^{(0)}_\be$  at $\la=u_j(t)$ has the form
\ben
I^{(0)}_\be(t,\la) = 2(2(\la-u_j))^{-1/2} e_j +\cdots,\quad
e_j=du_j/\sqrt{\Delta_j},
\een
where the dots stand for terms that have at $\la=u_j$ a zero of order
at least $\frac{1}{2}.$  These terms do not contribute to the residue,
so we get 
\ben
-\operatorname{Res}_{\la=u_j(t)} \frac{\la}{2}
d\la\,(V_{00}(t)e_j,e_j)\frac{2}{\la-u_j(t)} = u_j(t)\,(R_1(t)e_j,e_j).
\een
We get the following formula for the residue \eqref{res-i}:
\ben
(\widehat{\Psi}_t\widehat{R}_t)^{\otimes 2} \Big(u_jR_1^{jj}-
\operatorname{Res}_{\la=u_j} \frac{\la}{4} d\la\,
:\phi^{V\otimes V}_{A_1}(u_j,\la) \phi^{V\otimes V}_{A_1}(u_j,\la): 
\Big) \prod_{m=1}^{N+1} \cD_{\rm
  pt}(\hbar\Delta_m;\leftexp{m}{\q})^{\otimes 2},
\een
where $R_1^{jj}=(R_1e_j,e_j)$ is the $j$-th diagonal entry of $R_1$.

\subsubsection{The $A_1$-subroot system}\label{a1-term}
The vanishing cycles $\{-\be,\be\}$ form a subroot system of type
$A_1$. Let us compute the residue of the corresponding vertex operator
terms, i.e.,
\beq\label{a1-res}
\operatorname{Res}_{\la=u_j(t)}\,\frac{d\la}{\la}\,\Big(\sum_{\pm}
b_{\pm\be}(t,\la)\Gamma^{\pm\be}(t,\la)\otimes \Gamma^{\mp\be}(t,\la)
\Big)\,\cA_t^{\otimes 2} .
\eeq
We have $b_\be(t,\la) = b_{-\be}(t,\la)$ and 
\ben
b_\be(t,\la)
\Gamma^{\pm\be}(t,\la)\otimes \Gamma^{\mp\be}(t,\la)
(\widehat{\Psi}_t\widehat{R}_t)^{\otimes 2} = 
(\widehat{\Psi}_t\widehat{R}_t)^{\otimes 2}
b_{A_1}(u_j,\la)
\Gamma_{A_1}^{\pm\be}(u_j,\la)\otimes \Gamma_{A_1}^{\mp\be}(u_j,\la),
\een
where we used formula \eqref{R:vop} together with the identity 
\ben
b_\be(t,\la) e^{V_t(\f_\be(t,\la)_-,\f_\be(t,\la)_-)} = b_{A_1}(u_j,\la),
\een
which follows immediately from \eqref{intertw-R}.
Using that $\cA_t=\widehat{\Psi}_t\widehat{R}_t\prod_j \cD_{\rm
  pt}^{(j)}$, where  the factors 
$\cD_{\rm  pt}^{(j)}= \cD_{\rm pt}(\hbar\Delta_j;\leftexp{j}{\q})$ are solutions to
KdV, we can compute the residue \eqref{a1-res} via the Kac-Wakimoto
form of the KdV hierarchy \eqref{kdv-1}. After a short computation we
get that the residue \eqref{a1-res} is
\ben
(\widehat{\Psi}_t\widehat{R}_t)^{\otimes 2} \Big( 
\frac{1}{8} + \operatorname{Res}_{\la=u_j} \frac{\la}{4} d\la\,
:\phi^{V\otimes V}_{A_1}(u_j,\la) \phi^{V\otimes V}_{A_1}(u_j,\la): 
\Big) \prod_{m=1}^{N+1} \cD_{\rm
  pt}(\hbar\Delta_m;\leftexp{m}{\q})^{\otimes 2}.
\een

\subsubsection{The $A_2$-subroot subsystem}\label{a2-term}
Let $\al\in \Delta'$ be a cycle such that $(\al|\be)=1.$ We claim that
the expression 
\beq\label{vop-a2}
\Big(b_\al(t,\la)\Gamma^{\al}(t,\la)\otimes \Gamma^{-\al}(t,\la) +
b_{\al-\be}(t,\la) \Gamma^{\al-\be}(t,\la)\otimes
\Gamma^{-\al+\be}(t,\la) \Big)\cA_t^{\otimes 2}
\eeq
is analytic near $\la=u_j$. Using the decompositions
\ben
\al=\al'+\be/2,\quad \al-\be=\al'-\be/2,
\een
where $(\al'|\be)=0$, the above expression can be written as
\ben
\Gamma^{\al'}\otimes \Gamma^{-\al'}\Big(
a'\Gamma^{\be/2}\otimes \Gamma^{-\be/2} + 
a''\Gamma^{-\be/2}\otimes \Gamma^{\be/2}\Big)\ \cA_t^{\otimes 2},
\een
where the coefficients $a'$ and $a''$ are given by
\ben
a'(t,\la)&=&
\lim_{\mu\to \la} 
\Big(1-\frac{\mu}{\la}\Big)^{-2}
B_{\al,\al}(t,\la,\mu)B^{u_j}_{\al',-\be}(t,\la,\mu),\\
a''(t,\la)&=&\lim_{\mu\to \la}
\Big(1-\frac{\mu}{\la}\Big)^{-2}
B_{\al-\be,\al-\be}(t,\la,\mu)B^{u_j}_{\al',\be}(t,\la,\mu),
\een
where the phase factor $B^{u_j}_{\al',\be}=\exp \Omega^{u_j}_{\al',\be}$ with 
\ben
\Omega^{u_j}_{\al',\be}(t,\la,\mu)= \iota_{\la-u_j}\iota_{\mu-u_j} \sum_{n=0}^\infty (-1)^{n+1}(I^{(n)}_{\al'}(t,\la), I^{(-n-1)}_\be(t,\mu)),
\een
where $\iota_{\la-u_j}$ (resp. $\iota_{\mu-u_j}$) is the Laurent series expansion at $\la=u_j$ (resp. $\mu=u_j$). Since the Laurent series
expansions are convergent for $\la$ and $\mu$ sufficiently close to $u_j$, integration by parts yields
\ben
\Omega^{u_j}_{\al',\be}(t,\la,\mu)=\lim_{\ep\to 0} \int_{L_\ep} \cW_{\al',\be}(\mu-\la),
\een
where $L_\ep$ is the straight segment $[ t+(\ep+\mu-\la-u_j)\one, t-\la\one].$
On the other hand we have
\ben
\Gamma^{\pm\be/2}\otimes\Gamma^{\mp\be/2}(\widehat{\Psi}_t\widehat{R}_t)^{\otimes
  2} =
(\widehat{\Psi}_t\widehat{R}_t)^{\otimes 2}\,
e^{V_t(\f_{\be/2}(t,\la)_-,\f_{\be/2}(t,\la)_-)}
\Gamma_{A_1}^{\pm\be/2}\otimes\Gamma_{A_1}^{\mp\be/2}.
\een
The exponential factor can be
expressed in terms of the phase factors as follows (cf. Section
\ref{sec:vop-kdv}):
\ben
e^{V_t(\f_{\be/2}(t,\la)_-,\f_{\be/2}(t,\la)_-)}=
\frac{1}{2\sqrt{\la-u_j}}\, \lim_{\mu\to \la} (\la-\mu)^{1/2}B^{u_j}_{\be/2,-\be/2}(t,\la,\mu),
\een
where the limit is taken in the region $|\la|>|\mu|$. Recalling the
KP-reduction HQEs of KdV \eqref{kdv-2} we get that if the coefficients
\ben
c'(t,\la)=
\la^2\,
\lim_{\mu\to \la} 
(\la-\mu)^{-3/2}
B_{\al,\al}(t,\la,\mu)B^{u_j}_{\al',-\be}(t,\la,\mu) B^{u_j}_{\be/2,-\be/2}(t,\la,\mu)
\een 
and 
\ben
c''(t,\la)=
\la^2\,
\lim_{\mu\to \la}
(\la-\mu)^{-3/2}
B_{\al-\be,\al-\be}(t,\la,\mu)B^{u_j}_{\al',\be}(t,\la,\mu)
B^{u_j}_{\be/2,-\be/2}(t,\la,\mu)
\een
are analytic near $\la=u_j$, and $c'/c''=-1$, then the expression
\eqref{vop-a2} is analytic near $\la=u_j$. 

Let us prove the analyticity of $c'$. The argument for $c''$ is
similar. Let us choose a small $\ep\in \C$
and a generic point $t_0\in M$ on the discriminant, so that Proposition \ref{phase-mon} holds. Furthermore, we fix
2 paths $C'_\ep,$ and $C''_\ep$ in $M'=M\setminus{\{{\rm  discr}\}}$ from  $t_0+(\mu-\la+\ep)\one$ to $t-\la\one$ such that the
parallel transport transforms the cycle $\varphi$ vanishing
over $t_0$ respectively into $\al$, and $\al-\be$.  The phase
factors in the definition of $c'$ can be written in terms of integrals
along the path as follows 
\ben
B_{\al,\al}(t,\la,\mu) & = & \lim_{\ep\to 0} \exp\Big( \int_{C'_\ep}
\cW_{\al,\al}(\mu-\la)\Big),\\
B_{\al',-\be}(t,\la,\mu) & = & \lim_{\ep\to 0} \exp\Big( \int_{L_\ep}
\cW_{\al',-\be}(\mu-\la)\Big),\\
B_{\be/2,-\be/2}(t,\la,\mu) & = & \lim_{\ep\to 0} \exp\Big( \int_{L_\ep}
\cW_{\be/2,-\be/2}(\mu-\la)\Big).
\een
Using these formulas, we can express the coefficient $c'(t,\la)$ as the
limit $\ep\to 0$ of the following expression:
\ben
\la^2\, \lim_{\mu\to \la}(\la-\mu)^{-3/2}\exp\Big( 
\int_{C'_\ep} \cW_{\al,\al}(\mu-\la) -
\int_{L_\ep} \cW_{\al,\al}(\mu-\la) +
\int_{L_\ep} \cW_{\al',\al'}(\mu-\la) 
\Big).
\een
Let us examine the dependence on the parameters $t,\la$, and
$\xi:=\mu-\la$. 
The difference of the first two integrals in the above formula does not
depend on $\la$, because the paths $C'_\ep$ and $L_\ep$ have
the same ending point, while the starting points are independent of $\la$. After passing to the limit the
difference contributes a constant independent of $\la$, and
$\mu$. The last integral is analytic near $\la=u_j$, because the cycle
$\al'$ is invariant with respect to the local monodromy, which means that the
period vector $I^{(0)}_{\al'}(t',\xi)$  and respectively the phase
form $\cW_{\al',\al'}(\xi)$ are analytic for $t'$
sufficiently close to $t-u_j\one$ and $|\xi|\ll 1.$   This proves the
analyticity of $c'$. 

It remains only to prove that $c'/c''=-1.$ Using the above path
integrals, we can write $\log(c'/c'')$ in the following way:
\ben
\int_{C'_\ep} \cW_{\al,\al}-
\int_{L_\ep} \cW_{\al,\al} -
\int_{C''_\ep}\cW_{\al-\be,\al-\be} +
\int_{L_\ep} \cW_{\al-\be,\al-\be} +
\int_{\gamma_\ep} \cW_{\al,\al}-
\int_{\gamma_\ep} \cW_{\al,\al},
\een
where $\gamma_\ep$ is a small loop in $M'$ based at the starting point of $L_\ep$ (i.e. $t+(\ep+\mu-\la-u_j)\one$) that goes
counterclockwise around the discriminant and the branch of the phase form is determined by its value at the
point $t-\la\one$ (which belongs to the integration paths of the first 4 integrals and it is connected via the line segment $L_\ep$ to the contour of the 
last 2 ones). The above expression coincides with 
\ben
\oint_{(C''_\ep)^{-1}\circ L_\ep\circ \gamma_\ep\circ L_\ep^{-1}\circ C'_\ep} \cW_{\al,\al} -
\oint_{\gamma_\ep} \cW_{\al,\al}.
\een
By definition the cycle $\al$ is invariant along the integration contour of the first integral, 
so the first integral is an integer multiple of $2\pi\sqrt{-1}$. We get
\ben
c'/c''= \lim_{\xi\to 0} \lim_{\ep\to 0} \, \exp \Big(
-\oint_{\gamma_{\ep}} 
\cW_{\al,\al}(\xi)\Big).
\een
The limit here is easy to compute because the integral involves only 
local information. Using again the decomposition $\al=\al'+\be/2$ and Lemma  \ref{vanishing_a1} we get 
\ben
I^{(0)}_{\be}(t',\xi) = 2(2(\xi-u))^{-1/2}\, \frac{du}{\sqrt{\Delta}}+\cdots,
\een  
where the dots stand for higher order terms. On the other hand, the
period vector $I_{\al'}^{(0)}(t',\xi)$ is analytic for $(t',\xi)$
sufficiently close $(t,u_j).$ Expanding the phase form into a Laurent
series about $\xi=u$ we get  
\ben
\lim_{\ep\to 0} \oint_{\gamma_{\ep}} \cW_{\al,\al}(\xi) = 
\frac{1}{4} \oint_{\gamma_{\ep}} \cW_{\be,\be}(\xi) = 
\frac{1}{4} \oint \frac{2 du}{\sqrt{(-u)(\xi-u)}} = \pi\sqrt{-1},
\een
i.e., $c'/c''=-1.$

\subsubsection{Proof of Theorem \ref{thm:ancestor}}
The 1-form 
\ben
\frac{d\la}{\la} \Omega_{\Delta'}(t,\la)\, \cA_t(\hbar;\q')\cA_t(\hbar;\q'')
\een
has poles only at $\la=0,\infty$, and the critical values $u_j$,
$1\leq j\leq N+1$. Let $u_j$ be one of the critical values and $\be$
be the cycle vanishing over $\la=u_j$. Note that non-trivial
contributions to the residue at $\la=u_j$ come only from
vertex operator terms corresponding to vanishing cycles that have
non-zero intersection with $\be$. Recalling our computations in Sections
\ref{vir-term}, \ref{a1-term}, and \ref{a2-term}, we get that the residue at
$\la=u_j$ is $(1/8 + u_jR_1^{jj})\cA_t^{\otimes 2},$ while the residue
at $\la=0$ is 
$-\frac{1}{2}\operatorname{tr}\left(\frac{1}{4}+\theta\theta^T\right) \cA_t^{\otimes 2}$.
In order to prove that the
residue at $\la=\infty$ is 0, we just need to check that 
\ben
\sum_{j=1}^{N+1} u_jR_1^{jj} = 
\frac{1}{2}\operatorname{tr}\left(\theta\theta^T\right).
\een
The above identity is well-known from the theory of Frobenius
manifolds (see \cite{GM,He}). Hence the ancestor potential
$\cA_t(\hbar;\q)$ is a solution to the HQEs \eqref{eth-anc}. Theorem \ref{thm:ancestor} is thus proved.
\qed

\proof[Proof of Theorem \ref{t1}]
Given Theorem \ref{thm:ancestor}, Proposition \ref{anc-desc} implies that the total
descendant potential $\mathcal{D}_{\bf a}(\hbar; \q)$ is a solution to the HQEs
\eqref{eth-ade}. Theorem \ref{t1} then follows from Theorem \ref{t2}. 
\qed

\section{An example: $\mathbb{P}^1_{2,2,2}$}\label{ex}
In this section we consider the example ${\bf a}=\{2,2,2\}$, namely $\mathbb{P}^1_{\bf a}=\mathbb{P}^1_{2,2,2}$. In this case $\Delta^{(0)}$ is the root system of type $D_4$. It is convenient to denote the indexes in the
index set $\mathfrak{I}_{\rm tw}=\{(1,1),(2,1),(3,1)\}$ simply by $1,2,3.$ There are 12 positive roots 
\ben
&&
\gamma_i\ (1\leq i\leq 3), \quad \gamma_b, \quad \gamma_b+\gamma_i\ (1\leq i\leq 3), \quad
\gamma_b+\gamma_i+\gamma_j\ (1\leq i<j\leq 3), \\
&&
\gamma_b+\gamma_1+\gamma_2+\gamma_3,\quad  2\gamma_b+\gamma_1+\gamma_2+\gamma_3,
\een
where $\gamma_b$ is the simple root corresponding to the branching node of the Dynkin diagram and 
$\gamma_i \, (1\leq i\leq 3)$ are the remaining simple roots. The fundamental weight is
$\omega_b=2\gamma_b+\gamma_1+\gamma_2+\gamma_3$. The eigenbasis 
for $\si_b$ used in our construction is 
\ben
H_i:=-(\kappa/2)^{1/2}\gamma_i \quad (1\leq i\leq 3),\quad H_0:=(\kappa/2)^{1/2}\omega_b,
\een
and we have $m_i=\frac{\kappa}{2}$, $d_i=\frac{1}{2}$, $1\leq i\leq 3$, where $\kappa=4$. 

Let us write the HQEs for $\tau=(\tau_n(y))_{n\in \Z}$. We have
\ben
a_\al(\zeta) = \frac{1}{4}2^{(\si_b(\al)|\al)} \zeta^{\kappa|\al_0|^2}e^{2\pi\sqrt{-1}(\rho_b|\al)(\omega_b|\al)}
\een
and 
\ben
\Big(E_\al(\zeta)\tau\Big)_0 = \zeta^{-\kappa|\al_0|^2}e^{-2\pi\sqrt{-1}(\rho_b|\al)(\omega_b|\al)} 
E^*_\al(\zeta)\tau_{-(\omega_b|\al)},
\een
where the subscript $0$ on the LHS means the $0$-th component of the
corresponding vector in our Fock space. Recall that the HQEs give rise to a system of PDEs in the 
following way. First we make a substitution 
\ben
\y':=\y\otimes 1 = \x+\t,\quad \y'':=1\otimes \y=\x-\t,
\een
which implies that
\ben
\y'-\y'' = 2\t,\quad \frac{\d}{\d \y'} - \frac{\d}{\d \y''} = \frac{\d}{\d\t},
\een
and that 
\ben
\operatorname{Res}_{\zeta=0}\Big( a_\al(\zeta)E_\al(\zeta)\tau \otimes E_{-\al}(\zeta)\tau\Big)_{0,0} 
\een
is the coefficient in front of $\zeta^{0}$ in the following expression
\ben
2^{(\si_b(\al)|\al)-2} e^{-2\pi\sqrt{-1}(\rho_b|\al)(\omega_b|\al)}
\Big(\zeta^{-\kappa|\al_0|^2}
e^{\sum_{i,\ell} 2(\al|H_i)\zeta^{m_i+\ell\kappa}t_{i,\ell}} \Big) \\
\Big(e^{-\sum_{i,\ell} (\al|H_{i^*})\frac{\zeta^{-m_i-\ell\kappa}}{m_i+\ell\kappa}\d_{x_{i,\ell}}}\tau_{-(\omega_b|\al)}(\x+\t)\Big)
\Big(e^{\sum_{i,\ell} (\al|H_{i^*})\frac{\zeta^{-m_i-\ell\kappa}}{m_i+\ell\kappa}\d_{x_{i,\ell}}}\tau_{(\omega_b|\al)}(\x-\t)\Big).
\een
By definition the HQEs are
\ben
&&
\operatorname{Res}_{\zeta=0}\sum_{\al\in \Delta^{(0)}} 
\Big( a_\al(\zeta)E_\al(\zeta)\tau \otimes E_{-\al}(\zeta)\tau\Big)_{m,n} = \\
&&
\Big(\frac{3}{8} + \frac{1}{4}(m-n)^2+2\sum_{i,\ell} (d_{i^*}+\ell)t_{i,\ell}\d_{t_{i,\ell}} \Big) \tau_m(\x+\t)\tau_n(\x-\t).
\een
Comparing the coefficients in front of the various monomials in $\t$ we obtain a system of PDEs whose equations 
are some quadratic polynomials in the partial derivatives of $\tau$. Let us specialize to the case $m=n=0$.
In order to get non-trivial equations we have to compare coefficients in front of monomials that are invariant under
the involution $\t\mapsto -\t$. The simplest case is $\t^0$, which corresponds to the identity
\ben
\sum_{\al\in \Delta^{(0)}:(\omega_b|\al)=0} 2^{(\si_b(\al)|\al)-2} = \frac{3}{8}.
\een
Comparing the coefficients in front of the monomial $t_{02,0}^2$, we get 
\ben
4\frac{\d^2}{\d x_{02,0}^2}\, \log \tau(\x) =
8\kappa \frac{\tau_{-2}(\x)\tau_2(\x) }{\tau^2(\x)}-
\left. 4(2/\kappa)^{1/2}\frac{\d^3}{\d t_{1,0}\d t_{2,0}\d
    t_{3,0}}\Big(\frac{\tau_{-1}(\x+\t)\tau_1(\x-\t)}{\tau^2(\x)}\Big)\right|_{\t=0} .
\een
Recalling the substitution
\eqref{dv:change1}--\eqref{dv:change2}, which in this case is 
\ben
y_{02,0} & = & \frac{1}{\sqrt{\hbar}}\,
\frac{\sqrt{2}}{\kappa\sqrt{\kappa}} q_0^{02},\\
y_{i,0} & = &  
\frac{1}{\sqrt{\hbar}}\,
\frac{\sqrt{2}}{\kappa}\, q_0^i,\quad 1\leq i\leq 3,
\een
we get
\ben
\hbar \frac{\d^2}{\d (q_0^{02})^2} \log \tau(\q) =
\frac{4}{\kappa^2}\frac{\tau_{-2}(\q)\tau_2(\q) }{\tau^2(\q)} -
\frac{\hbar^{3/2}}{\kappa^{1/2}}
\left.
\d_1\d_2\d_3\Big(\frac{\tau_{-1}(\q+\t)\tau_1(\q-\t)}{\tau^2(\q)}\Big)\right|_{\t=0}  ,
\een
where for brevity we put $\d_i:=\d/\d t_0^i$. To get a differential
equation for the total descendant potential we just have to substitute
\ben
\tau_{\pm 2}(\q) = C^2\cD(\hbar;\q\pm 2\sqrt{\hbar}),\quad 
\tau_{\pm 1}(\q) = C^{1/2}\cD(\hbar;\q\pm \sqrt{\hbar})
\een
where $C=\kappa^{1/2}Q$. 

Let us use the above equation to compute the genus-$0$ primary potential $F$. 
Put $q_k^i=0,$ $\forall k>0$, and compare the
leading terms of the genus expansion. We get the following PDE for $F$:
\ben
F_{02,02} = 4Q^4 e^{4F_{01,01}} + 
Q e^{F_{01,01}} 
\Big(8F_{01,1}F_{01,2}F_{01,3} + 4(F_{01,1} F_{2,3}+F_{01,2}F_{1,3}+F_{01,3}F_{1,2})\Big),
\een
where $F_{i,j}:=\d^2F/\d q_0^i\d q_0^j$. To simplify the notation, let us put $t_i:=q_0^i$. String equation gives
\ben
F_{01,01} = t_{02},\quad F_{01,i} = \frac{1}{2}t_i,
\een
so from the above equation we get the following relation
\beq\label{top-sol}
F_{02,02} = 4Q^4 e^{4t_{02}} + 
Qe^{t_{02}} 
\Big(t_1t_2t_3 + 2(t_1F_{2,3}+t_2F_{1,3}+t_3F_{1,2})\Big).
\eeq
Equation \eqref{top-sol} allows us to compute the potential $F$ recursively, by the degree of the Novikov variable $Q$. 
Indeed, it is easy to see that up to degree-1 terms, $F$ is given by
\ben
\frac{1}{2} t_{01}^2t_{02} + \frac{1}{4}t_{01}(t_1^2+t_2^2+t_3^2) + \frac{1}{96}(t_1^4+t_2^4+t_3^4) 
+ Qe^{t_{02}}  t_1t_2t_3.
\een
Comparing the degree-2 terms in \eqref{top-sol} we get that the degree-2
term of $F$ must be $\frac{1}{2}(t_1^2+t_2^2+t_3^2)Q^2e^{2t_{02}}$. Arguing in the same way we get that 
$F$ does not have degree-3 terms, while the degree 4 term must be $\frac{1}{4} Q^4e^{4t_{02}}.$ 
The potential $F$ takes the form
\ben
F(t) &=& \frac{1}{2} t_{01}^2t_{02} + \frac{1}{4}t_{01}(t_1^2+t_2^2+t_3^2) + \frac{1}{96}(t_1^4+t_2^4+t_3^4) 
+ Qe^{t_{02}}  t_1t_2t_3 + \\
&&
+\frac{1}{2}Q^2e^{2t_{02}} (t_1^2+t_2^2+t_3^2)+\frac{1}{4} Q^4e^{4t_{02}}.
\een
The above formula agrees with the computation of P. Rossi \cite[Example 3.2]{Ro} based on Symplectic Field Theory.

\appendix

\section{An alternative proof of higher genus reconstruction}
In this subsection, we use the degree of virtual fundamental cycle and tautological relations to give a simple proof for Teleman's higher genus reconstruction theorem for the target $\mathbb{P}^1_{\bf a}$, see Proposition \ref{anc-recons} below. This proof does not require the semi-simple assumption. 
  
We first recall the $g$-reduction property introduced in \cite{FSZ}, which is a consequence of results by Ionel \cite{Io}, and by Faber and Pandharipande \cite{FP}:
\begin{lemma}[\cite{Io,FP}]\label{lem-g}
 If $M(\psi, \kappa)$ is a polynomial of $\psi$-classes and $\kappa$-classes with $\deg M\geq g$ for $g\geq1$ or $\deg M\geq1$ for $g=0$, then $M(\psi, \kappa)$ can be presented as a linear combination of dual graphs on the boundary of $\overline{\mathcal{M}}_{g,n}$.
\end{lemma}
Our second tool is the Getzler's relation in \cite{Get}. It is a linear relation between codimension two cycles in $H_*(\overline{\mathcal{M}}_{1,4},\mathbb{Q})$. Here we briefly introduce this relation for our purpose. The dual graph
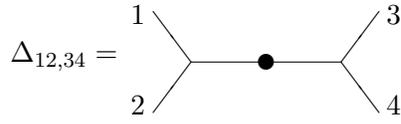
\begin{center}
\begin{picture}(50,20)
    \put(-24,9){$\Delta_{12,34}=$}

    \put(10,9){\circle*{2}}

	\put(0,9){\line(-3,4){5}}
    \put(0,9){\line(-3,-4){5}}
	\put(0,9){\line(1,0){9}}
	\put(11,9){\line(1,0){9}}
    \put(20,9){\line(3,4){5}}
    \put(20,9){\line(3,-4){5}}

    \put(-8,2){2}
    \put(-8,14){1}
    \put(26,2){4}
    \put(26,14){3}
\end{picture}
\end{center}
 represents a codimension-two stratum in $\overline{\mathcal{M}}_{1,4}$: A filled circle represents a genus-1 component, other vertices represent genus-0 components. An edge connecting two vertices represents a node, a tail (or half-edge) represents a marked point on the component of the corresponding vertex. $\Delta_{2,2}$ is defined to be the $S_4$-invariant of the codimension-two stratum in $\overline{\mathcal{M}}_{1,4}$,
\begin{equation*}
\Delta_{2,2}=\Delta_{12,34}+\Delta_{13,24}+\Delta_{14,23}.
\end{equation*}
We denote $\delta_{2,2}=[\Delta_{2,2}]$ the corresponding cycle in $H_4(\overline{\mathcal{M}}_{1,4},\mathbb{Q})$. We list the corresponding unordered dual graph for other strata below,  see \cite{Get} for more details.
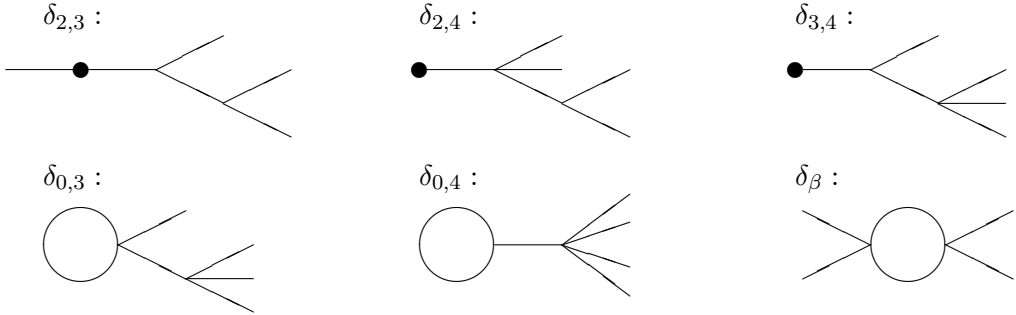
\begin{center}
\begin{picture}(50,20)

    \put(-35,9){\circle*{2}}
	\put(-40,15){$\delta_{2,3}:$}

	\put(-45,9){\line(1,0){9}}
	\put(-34,9){\line(1,0){9}}
    \put(-25,9){\line(2,1){9}}
    \put(-25,9){\line(2,-1){9}}
    \put(-16,4.5){\line(2,1){9}}
    \put(-16,4.5){\line(2,-1){9}}

    \put(10,9){\circle*{2}}
	\put(10,15){$\delta_{2,4}:$}


	\put(11,9){\line(1,0){9}}
    \put(20,9){\line(2,1){9}}
	\put(20,9){\line(1,0){9}}
    \put(20,9){\line(2,-1){9}}
    \put(29,4.5){\line(2,1){9}}
    \put(29,4.5){\line(2,-1){9}}

    \put(60,9){\circle*{2}}
	\put(60,15){$\delta_{3,4}:$}

	\put(79,4.5){\line(1,0){9}}
	\put(61,9){\line(1,0){9}}
    \put(70,9){\line(2,1){9}}
    \put(70,9){\line(2,-1){9}}
    \put(79,4.5){\line(2,1){9}}
    \put(79,4.5){\line(2,-1){9}}
\end{picture}
\end{center}

\begin{center}
\begin{picture}(50,20)

	\put(-40,17){$\delta_{0,3}:$}
	\put(10,17){$\delta_{0,4}:$}
	\put(60,17){$\delta_{\beta}:$}

    \put(-35,9){\circle{10}}


	\put(-21,4.5){\line(1,0){9}}
    \put(-30,9){\line(2,1){9}}
    \put(-30,9){\line(2,-1){9}}
    \put(-21,4.5){\line(2,1){9}}
    \put(-21,4.5){\line(2,-1){9}}

    \put(15,9){\circle{10}}


    \put(29,9){\line(4,3){9}}
    \put(29,9){\line(3,1){9}}
	\put(20,9){\line(1,0){9}}
    \put(29,9){\line(4,-3){9}}
    \put(29,9){\line(3,-1){9}}


    \put(75,9){\circle{10}}

    \put(80,9){\line(2,1){9}}
    \put(80,9){\line(2,-1){9}}
    \put(70,9){\line(-2,1){9}}
    \put(70,9){\line(-2,-1){9}}
\end{picture}
\end{center}

In \cite{Get}, Getzler found the following identity:
\begin{equation}\label{eq:Getzler}
  12\delta_{2,2}+4\delta_{2,3}-2\delta_{2,4}+6\delta_{3,4}+\delta_{0,3}+\delta_{0,4}-2\delta_{\beta}=0\in H_4(\overline{\mathcal{M}}_{1,4},\mathbb{Q}).
\end{equation}

Now we prove the following higher genus reconstruction result.
\begin{proposition}\label{anc-recons}
The total ancestor potential $\mathcal{A}_{\bf a}(\hbar;\t)$ is uniquely determined by the quantum cohomology of $\mathbb{P}^1_{\bf a}$ when ${\bf a}\neq\{1,1,1\}$ and $\chi>0$.
\end{proposition}
\proof
We consider the ancestor correlator $\langle \phi_{1}\,\bar{\psi}_1^{k_1},\dots,\phi_{n}\bar{\psi}_n^{k_n}\rangle_{g,n,d}$ in \eqref{anc-cor}.
According to the degree formula \eqref{virt-deg}, if the correlator is nonzero, then
\beq\label{dim-virt}
\frac{1}{2}\sum_{j=1}^n\deg\phi_{j}+\sum_{j=1}^{n}k_j=(3-\frac{1}{2}\dim\mathbb{P}^1_{\bf a})(g-1)+\chi\cdot d+n.
\eeq
Now if $\sum_{j=1}^{n}k_j\geq g$ for $g\geq1$ or $\sum_{j=1}^{n}k_j\geq1$ for $g=0$, then we can apply Lemma \ref{lem-g} of $g$-reduction to rewrite the ancestor correlator as a linear combination of intersection numbers over the corresponding homology cycles of some dual graphs, each of the dual graph lives on the boundary of $\overline{\mathcal{M}}_{g,n}$. The splitting axiom in GW theory allows us to reconstruct the ancestor correlator in \eqref{anc-cor} using intersection numbers over each component of the boundaries. We can keep doing this process until on each component, the $g$-reduction property does not hold. In other words, all the ancestor correlators are determined completely by those \eqref{anc-cor} which satisfies $\sum_{j=1}^{n}k_j\leq g-1$ for $g\geq1$ or $\sum_{j=1}^{n}k_j=0$ for $g=0$. On the other hand, since $\deg\phi_{j}\leq2$, $\chi>0$ and $\dim\mathbb{P}^1_{\bf a}=1$, the formula \eqref{dim-virt} implies such intersection numbers must vanish unless $g=0$ and all $k_{j}=0$, or $g=1, d=0$, all $k_{j}=0$ and all $\deg\phi_{j}=2$. 

In order to finish the proof, it only remains to consider genus $1$ correlator $\langle P\rangle_{1,1,0}$. If ${\bf a}\neq\{1,1,1\}$, then according to Rossi's computation \cite{Ro}, we can always find a twisted sector $\phi_{i}\in H$, such that
\beq\label{4-pt}
\langle \phi_{i}, \phi_{i}, \phi_{i^*}, \phi_{i^*}\rangle_{0,4,0}\neq0.
\eeq
Consider the integration of the cohomology cycle $\Lambda_{1,4,0}(\phi_{i}, \phi_{i}, \phi_{i^*}, \phi_{i^*})$ over the Getzler's relation \eqref{eq:Getzler}, with four fixed insertions $\phi_{i}, \phi_{i}, \phi_{i^*}, \phi_{i^*}$. By the splitting axiom in GW theory, it is not hard to see that the integration vanishes on those homology classes with a genus-$1$ component except that
$$\int_{\delta_{3,4}}\Lambda_{1,4,0}(\phi_{i}, \phi_{i}, \phi_{i^*}, \phi_{i^*})$$
is a multiplication by a nonzero scalar and $\langle P\rangle_{1,1,0}$, because of \eqref{4-pt}. Thus the equality \eqref{eq:Getzler} implies $\langle P\rangle_{1,1,0}$ is reconstructed from genus-0 correlators. 
\qed
\begin{remark}
The technique above only uses properties of cohomology field theories and tautological relations over the moduli space of stable curves. So it also works for the reconstruction of the ancestor potential in \eqref{ancestor}. It also works for elliptic orbifold projective curves $\mathbb{P}^1_{\bf a}$, where $\chi=0$, see \cite{KS}. The genus-1 correlator $\langle P\rangle_{1,1,0}$ in GW theory can be calculated directly using virtual cycle or virtual localization, see \cite{Ts}.
\end{remark}

\bibliographystyle{amsalpha}

\begin{thebibliography}{FKRW}

\bibitem{AL}
M. J.~Ablowitz, J. F.~Ladik.
\textit{Nonlinear differential-difference equations}. 
J. Math. Phys. 16 (1975), 598--603. 

\bibitem{AbGV}
D.~Abramovich, T.~Graber, A.~Vistoli.
\textit{Gromov-Witten theory of Deligne-Mumford stacks.}
Amer. J. Math. 130 (2008), no. 5, 1337--1398.

\bibitem{AGV}
V.I.~Arnol'd, S.M.~Gusein-Zade, A.N.~Varchenko.
\textit{Singularities of differentiable maps. Vol. II. Monodromy and asymptotics of integrals.}  
Monographs in Mathematics, 83. Birkh\"auser Boston, Inc., Boston, MA, 1988

\bibitem{BM}
B.~Bakalov, T.~Milanov.
\textit{$\mathcal{W}$-constraints for the total descendant potential of a simple singularity.}
Compositio Math. 149 (2013), no. 5, 840--888.

\bibitem{BCFK}
A.~Bertram, I.~Ciocan-Fontanine, B.~Kim.
\textit{Two proofs of a conjecture of Hori and Vafa.}
Duke Math. J. 126 (2005), no. 1, 101--136.

\bibitem{Br}
A.~Brini.
\textit{The local Gromov-Witten theory of $\mathbb{CP}^1$ and integrable hierarchies}. 
Comm. Math. Phys. 313 (2012), no. 3, 571--605.

\bibitem{BPS}
A.~Buryak, H.~Posthuma, S.~Shadrin.
\textit{On deformations of quasi-Miura transformations and the Dubrovin--Zhang bracket.}  
J. Geom. Phys. 62 (2012), no. 7, 1639--1651. 

\bibitem{BPS2}
A.~Buryak, H.~Posthuma, S.~Shadrin.
\textit{A polynomial bracket for the Dubrovin-Zhang hierarchies.}
J. Differential Geom. 92 (2012), no. 1, 153--185. 

\bibitem{Cadman}
C.~Cadman.
\textit{Using stacks to impose tangency conditions on curves.}
Amer. J. Math. 129 (2007), no. 2, 405--427. 

\bibitem{Ca}
G.~Carlet.
\textit{The extended bigraded Toda hierarchy}. 
J. Phys. A 39 (2006), no. 30, 9411--9435.

\bibitem{CDZ}
G.~Carlet, B.~Dubrovin, Y.~Zhang.
\textit{The extended Toda hierarchy}. 
 Mosc. Math. J. 4 (2004), no. 2, 313--332, 534. 
 
\bibitem{CvdL}
G.~Carlet, J.~van de Leur.
\textit{Hirota equations for the extended bigraded Toda hierarchy and the total descendent potential of $\mathbb{CP}^1$ orbifolds}.
J. Phys. A: Math. Theor. 46 (2013), 405205.
 
 
\bibitem{CR}
W.~Chen, Y.~Ruan.
\textit{Orbifold Gromov-Witten theory.} 
In: ``Orbifolds in mathematics and physics''. 25--85.
Contemp. Math., 310, Amer. Math. Soc., Providence, RI, 2002. 

\bibitem{CCIT1}
T.~Coates, A.~Corti, H.~Iritani, H.-H.~Tseng.
\textit{Computing Genus-Zero Twisted Gromov-Witten Invariants}.
Duke Math. J. 147 (2009), no. 3, 377--438.


\bibitem{CCIT2}
T.~Coates, A.~Corti, H.~Iritani, H.-H.~Tseng.
\textit{A Mirror Theorem for Toric Stacks}.
arXiv:1310.4163.


\bibitem{Du} 
B.~Dubrovin.
\textit{Geometry of 2D topological field theories.} 
In: ``Integrable Systems and Quantum Groups''. 120--348.
Lecture Notes in Math. 1620, Springer-Verlag, New York, 1996.

\bibitem{DZ}
B.~Dubrovin, Y.~Zhang.
\textit{Extended Affine Weyl groups and Frobenius manifolds}.
Compos. Math. 111 (1998), no. 2, 167--219.

\bibitem{DZ2}
B.~Dubrovin, Y.~Zhang.
\textit{Normal forms of hierarchies of integrable PDEs, Frobenius manifolds and Gromov-Witten invariants}.
a new 2005 version of arXiv:math/0108160v1, 295 pp.

\bibitem{DZ3}
B.~Dubrovin, Y.~Zhang.
\textit{Virasoro Symmetries of the Extended Toda Hierarchy.} 
Comm. Math. Phys. 250 (2004), no. 1, 161–193.

\bibitem{EY}
T. Eguchi, S.-K. Yang.
\textit{The topological $\C P^1$-model and the large-$N$ matrix integral.} 
Mod. Phys. Lett. A9 (1994), 2893--2902.

\bibitem{FP}
C.~Faber, R.~Pandharipande. 
\textit{Relative maps and tautological classes}. 
J. Eur. Math. Soc. 7 (2005), 13--49.

\bibitem{FSZ}
C.~Faber, S.~Shadrin, D.~Zvonkine. 
\textit{Tautological relations and the $r$-spin Witten conjecture.} 
Ann. Sci. c. Norm. Supr. (4) 43 (2010), no. 4, 621--658.

\bibitem{FJR}
H.~Fan, T.~Jarvis, Y.~Ruan.
\textit{Quantum singularity theory for $A_{(r-1)}$ and $r$-spin theory}. 
Ann. Inst. Fourier (Grenoble) 61 (2011), no. 7, 2781--2802.

\bibitem{FJR2}
H.~Fan, T.~Jarvis, Y.~Ruan.
\textit{The Witten equation, mirror symmetry, and quantum singularity theory}. 
Ann. of Math. (2) 178 (2013), no. 1, 1--106. 

\bibitem{FGM}
E.~Frenkel, A.~Givental, T.~Milanov. 
\textit{Soliton equations, vertex operators, and simple singularities}.
Funct. Anal. Other Math. 3 (2010), no. 1, 47--63.

\bibitem{FK}
I.~Frenkel, V.~Kac.
\emph{Basic representations of affine Lie algebras and the dual resonance models}.
Invent. Math. 62 (1980/81), no. 1, 23--66.

\bibitem{GGI}
S.~Galkin, V.~Golyshev, H.~Iritani.
\textit{Gamma classes and quantum cohomology of Fano manifolds: Gamma conjectures}.
arXiv:1404.6407.

\bibitem{Get}
E.~Getzler. 
\textit{Intersection theory on $\overline{\mathcal{M}}_{1,4}$ and elliptic Gromov-Witten invariants.}
J. Amer. Math. Soc. 10 (1997), no. 4, 973--998.

\bibitem{Ge}
E.~Getzler.
\textit{The Toda conjecture}.
In: ``Symplectic geometry and mirror symmetry (Seoul, 2000)'', 51--79, World Sci. Publ., River Edge, NJ, 2001.

\bibitem{G0}
A.~Givental.
\textit{A mirror theorem for toric complete intersections.}
In: ``Topological field theory, primitive forms and related topics (Kyoto, 1996)'', 141--175, Progr. Math., 160, Birkh\"auser Boston, Boston, MA, 1998. 

\bibitem{G1}
A.~Givental. 
\textit{$A_{n-1}$ singularities and $n$KdV Hierarchies}. 
Mosc. Math. J. 3 (2003), no. 2, 475--505.

\bibitem{G2} 
A.~Givental.
{\em Semi-simple Frobenius structures at higher genus.}
Internat. Math. Res. Notices 2001, no. 23, 1265--1286.

\bibitem{G3}
A.~Givental.
\textit{Gromov--Witten invariants and quantization of quadratic Hamiltonians}. 
Mosc. Math. J. 1 (2001), no. 4, 551--568.

\bibitem{GM} 
A.~Givental, T.~Milanov. 
\textit{Simple singularities and integrable hierarchies}. 
In: ``The breadth of symplectic and Poisson geometry''. 173--201.
Progr. Math., 232, Birkh\"auser Boston, Boston, MA, 2005

\bibitem{Gri}
P.~Griffiths.
\textit{On the periods of certain rational integrals: I.} 
Ann. of Math. (2) 90 (1969), no. 3, 460--495.

\bibitem{He}
C.~Hertling.
\textit{Frobenius manifolds and moduli spaces for singularities}. 
Cambridge Tracts in Math., 151, Cambridge Univ. Press, Cambridge, 2002


\bibitem{Io}
E.~Ionel.
\textit{Topological Recursive Relations in $H^{2g}(\mathcal{M}_{g,n})$}. 
Invent. Math. 148 (2002), no. 3, 627--658.


\bibitem{I}
H.~Iritani.
\textit{Convergence of quantum cohomology by quantum Lefschetz.}
 J. Reine Angew. Math. 610 (2007), 29--69.

\bibitem{Ir}
H.~Iritani.
\textit{An integral structure in quantum cohomology and mirror symmetry for toric orbifolds.}
Adv. Math. 222 (2009), no. 3, 1016--1079.

\bibitem{Iri}
H.~Iritani.
\textit{Quantum cohomology and periods.}
Ann. Inst. Fourier 61 (2011), 7, 2909--2958.

\bibitem{IST}
Y.~Ishibashi, Y.~Shiraishi, A.~Takahashi.
\textit{Primitive forms for affine cusp polynomials.}
arXiv:1211.1128.

\bibitem{JK}
D.~Joe, B.~Kim.
\textit{Equivariant mirrors and the Virasoro conjecture for flag manifolds.}
 Int. Math. Res. Not. 2003, no. 15, 859--882.

\bibitem{J}
P.~Johnson.
\textit{Equivariant Gromov-Witten theory of one dimensional stacks}.
arXiv:0903.1068.

\bibitem{Kac} 
V.~Kac. 
\textit{Infinite dimensional Lie algebras.} 
Third edition. Cambridge University Press, 1990.

\bibitem{Kac-v} 
V.~Kac. 
\textit{Vertex algebras for beginners.}
Second edition. University Lecture Series, 10. American Mathematical Society, Providence, RI, 1998.

\bibitem{KP} 
V.~Kac, D.~Peterson.
\textit{$112$ constructions of the basic representation of the loop group of $E_8$.} 
In: ``Symposium on anomalies, geometry, topology (Chicago, Ill., 1985)''. 276--298.
World Sci. Publishing, Singapore, 1985.

\bibitem{KKP}
L.~Katzarkov, M.~Kontsevich, T.~Pantev.
\textit{Hodge theoretic aspects of mirror symmetry.}
arXiv:0806.0107.

\bibitem{Ko1}
M.~Kontsevich, 
\textit{Intersection theory on the moduli space of curves and the matrix Airy function}. 
Comm. Math. Phys. 147 (1992), no. 1, 1--23.

\bibitem{KS}
M.~Krawitz, Y.~Shen.
\textit{Landau-Ginzburg/Calabi-Yau correspondence of all genera for elliptic orbifold $\mathbb{P}^1$}.
arXiv:1106.6270 [math.AG]

\bibitem{LRZ}
S.-Q.~Liu, Y.~Ruan, Y.~Zhang.
\textit{BCFG Drinfeld-Sokolov Hierarchies and FJRW-Theory}.
arXiv:1312.7227.


\bibitem{M}
T.~Milanov. 
\textit{Gromov--Witten theory of\/ $\mathbb{C}P^1$ and integrable hierarchies}. 
arXiv:math-ph/0605001.

\bibitem{M2}
T.~Milanov.
\emph{Hirota Quadratic Equations for the Extended Toda Hierarchy.} 
Duke Math. J. 138 (2007), no. 1, 161--178.

\bibitem{MT}
T.~Milanov, H.-H.~Tseng. 
\textit{The spaces of Laurent polynomials, Gromov--Witten theory of\/ $\mathbb{P}^1$-orbifolds, and integrable hierarchies}.
J. Reine Angew. Math. 622 (2008), 189--235.

\bibitem{OP}
A.~Okounkov, R.~Pandharipande. 
\textit{The equivariant Gromov-Witten theory of $\mathbf{P}^1$}. 
Ann. of Math. (2) 163 (2006), no. 2, 561--605.

\bibitem{Ro}
P.~Rossi.
\textit{Gromov--Witten theory of orbi-curves, the space of tri-polynomials, and symplectic field theory of Seifert fibrations}.
Math. Ann. 348 (2010), no. 2, 265--287.

\bibitem{Sa}
K.~Saito.
\textit{On Periods of Primitive Integrals, I.}
Preprint RIMS (1982). 

\bibitem{SaT}
K.~Saito, A.~Takahashi.
\textit{From primitive forms to Frobenius manifolds.}
In: ``From Hodge theory to integrability and TQFT tt*-geometry''. 31--48,
Proc. Sympos. Pure Math., 78, Amer. Math. Soc., Providence, RI, 2008

\bibitem{Ste}
R.~Stekolshchik.
\textit{Notes on Coxeter transformations and the McKay correspondence.}
Springer Monographs in Mathematics, Springer-Verlag Berlin Heidelberg, 2008.

\bibitem{Stb}
R.~Steinberg.
\textit{Finite subgroups of $SU_2$, Dynkin diagrams and affine Coxeter elements.} 
Pacific J. Math. 118 (1985), no. 2, 2587--598.

\bibitem{Ta}
A.~Takahashi.
\textit{Mirror symmetry between orbifold projective lines and cusp singularities.}
to appear in Adv. Stud. Pure Math.

\bibitem{Te}
C.~Teleman.
\textit{The structure of 2D semi-simple field theories.} 
Invent. Math. 188 (2012), no. 3, 525--588.

\bibitem{tseng}
H.-H.~Tseng.
\textit{Orbifold quantum Riemann-Roch, Lefschetz and Serre.}
Geom. Topol. 14 (2010), no. 1, 1--81. 

\bibitem{Ts}
H.-H.~Tseng. 
\textit{On degree-0 elliptic orbifold Gromov-Witten invariants.}
Internat. Math. Res. Notices 2011, no. 11, 2444--2468.



\bibitem{W1}
E.~Witten,
\textit{Two-dimensional gravity and intersection theory on moduli space}. 
In: ``Surveys in differential geometry,'' 243--310, 
Lehigh Univ., Bethlehem, PA, 1991.

\bibitem{W2}
E.~Witten,
\textit{Algebraic geometry associated with matrix models of  two-dimensional gravity}. 
Topological methods in modern mathematics (Stony Brook, NY, 1991), Publish or Perish, Houston, TX, 1993, 235--269.

\bibitem{Z}
Y. Zhang.
\textit{On the $\C P^1$ topological sigma model and the Toda lattice hierarchy.}
J. Geom. Phys. 40 (2002), no. 3-4, 215--232.

\end{thebibliography}

\end{document}